\PassOptionsToPackage{unicode}{hyperref}
\PassOptionsToPackage{hyphens}{url}
\PassOptionsToPackage{dvipsnames,svgnames,x11names}{xcolor}
\documentclass[
  english,
  11pt]{article}
\usepackage{xcolor}
\usepackage[margin=1in]{geometry}
\usepackage{amsmath,amssymb}
\usepackage{iftex}
\ifPDFTeX
  \usepackage[T1]{fontenc}
  \usepackage[utf8]{inputenc}
  \usepackage{textcomp} % provide euro and other symbols
\else % if luatex or xetex
  \usepackage{unicode-math} % this also loads fontspec
  \defaultfontfeatures{Scale=MatchLowercase}
  \defaultfontfeatures[\rmfamily]{Ligatures=TeX,Scale=1}
\fi
\usepackage{lmodern}
\ifPDFTeX\else
\fi
\IfFileExists{upquote.sty}{\usepackage{upquote}}{}
\IfFileExists{microtype.sty}{% use microtype if available
  \usepackage[]{microtype}
  \UseMicrotypeSet[protrusion]{basicmath} % disable protrusion for tt fonts
}{}
\makeatletter
\@ifundefined{KOMAClassName}{% if non-KOMA class
  \IfFileExists{parskip.sty}{%
    \usepackage{parskip}
  }{% else
    \setlength{\parindent}{0pt}
    \setlength{\parskip}{6pt plus 2pt minus 1pt}}
}{% if KOMA class
  \KOMAoptions{parskip=half}}
\makeatother
\makeatletter
\ifx\paragraph\undefined\else
  \let\oldparagraph\paragraph
  \renewcommand{\paragraph}{
    \@ifstar
      \xxxParagraphStar
      \xxxParagraphNoStar
  }
  \newcommand{\xxxParagraphStar}[1]{\oldparagraph*{#1}\mbox{}}
  \newcommand{\xxxParagraphNoStar}[1]{\oldparagraph{#1}\mbox{}}
\fi
\ifx\subparagraph\undefined\else
  \let\oldsubparagraph\subparagraph
  \renewcommand{\subparagraph}{
    \@ifstar
      \xxxSubParagraphStar
      \xxxSubParagraphNoStar
  }
  \newcommand{\xxxSubParagraphStar}[1]{\oldsubparagraph*{#1}\mbox{}}
  \newcommand{\xxxSubParagraphNoStar}[1]{\oldsubparagraph{#1}\mbox{}}
\fi
\makeatother

\usepackage{longtable,booktabs,array}
\usepackage{calc} % for calculating minipage widths
\usepackage{etoolbox}
\makeatletter
\patchcmd\longtable{\par}{\if@noskipsec\mbox{}\fi\par}{}{}
\makeatother
\IfFileExists{footnotehyper.sty}{\usepackage{footnotehyper}}{\usepackage{footnote}}
\makesavenoteenv{longtable}
\usepackage{graphicx}
\makeatletter
\newsavebox\pandoc@box
\newcommand*\pandocbounded[1]{% scales image to fit in text height/width
  \sbox\pandoc@box{#1}%
  \Gscale@div\@tempa{\textheight}{\dimexpr\ht\pandoc@box+\dp\pandoc@box\relax}%
  \Gscale@div\@tempb{\linewidth}{\wd\pandoc@box}%
  \ifdim\@tempb\p@<\@tempa\p@\let\@tempa\@tempb\fi% select the smaller of both
  \ifdim\@tempa\p@<\p@\scalebox{\@tempa}{\usebox\pandoc@box}%
  \else\usebox{\pandoc@box}%
  \fi%
}
\def\fps@figure{htbp}
\makeatother

\ifLuaTeX
\usepackage[bidi=basic,shorthands=off]{babel}
\else
\usepackage[bidi=default,shorthands=off]{babel}
\fi
\ifLuaTeX
  \usepackage{selnolig} % disable illegal ligatures
\fi

\providecommand{\tightlist}{%
  \setlength{\itemsep}{0pt}\setlength{\parskip}{0pt}}

\usepackage{amsmath,amssymb,bm,booktabs,graphicx,float,appendix,needspace}
\usepackage{etoolbox}
\makeatletter
\newcommand{\nmcontinue}{\@afterindentfalse\@afterheading}
\AfterEndEnvironment{equation}{\nmcontinue}
\AfterEndEnvironment{equation*}{\nmcontinue}
\AfterEndEnvironment{displaymath}{\nmcontinue}
\AfterEndEnvironment{align}{\nmcontinue}
\AfterEndEnvironment{align*}{\nmcontinue}
\AfterEndEnvironment{gather}{\nmcontinue}
\AfterEndEnvironment{gather*}{\nmcontinue}
\AfterEndEnvironment{multline}{\nmcontinue}
\makeatother
\usepackage{fancyhdr}
\makeatletter
\AtBeginDocument{%
  \@ifundefined{c@theorem}{}{%
    \@ifundefined{c@proposition}{}{\let\c@proposition\c@theorem}%
    \@ifundefined{c@corollary}{}{\let\c@corollary\c@theorem}%
    \@ifundefined{c@lemma}{}{\let\c@lemma\c@theorem}%
    \@ifundefined{c@example}{}{\let\c@example\c@theorem}%
    \@ifundefined{c@remark}{}{\let\c@remark\c@theorem}}%
  \@ifundefined{c@theorem}{}{\counterwithout{theorem}{section}}%
  \@ifundefined{c@proposition}{}{\counterwithout{proposition}{section}}%
  \@ifundefined{c@corollary}{}{\counterwithout{corollary}{section}}%
  \@ifundefined{c@lemma}{}{\counterwithout{lemma}{section}}%
  \@ifundefined{c@example}{}{\counterwithout{example}{section}}%
  \@ifundefined{c@definition}{}{\counterwithout{definition}{section}}%
  \@ifundefined{c@remark}{}{\counterwithout{remark}{section}}%
  \@ifundefined{c@refremark}{}{\counterwithout{refremark}{section}}}
\makeatother
\makeatletter
\@ifpackageloaded{tcolorbox}{}{\usepackage[skins,breakable]{tcolorbox}}
\@ifpackageloaded{fontawesome5}{}{\usepackage{fontawesome5}}
\definecolor{quarto-callout-color}{HTML}{909090}
\definecolor{quarto-callout-note-color}{HTML}{0758E5}
\definecolor{quarto-callout-important-color}{HTML}{CC1914}
\definecolor{quarto-callout-warning-color}{HTML}{EB9113}
\definecolor{quarto-callout-tip-color}{HTML}{00A047}
\definecolor{quarto-callout-caution-color}{HTML}{FC5300}
\definecolor{quarto-callout-color-frame}{HTML}{acacac}
\definecolor{quarto-callout-note-color-frame}{HTML}{4582ec}
\definecolor{quarto-callout-important-color-frame}{HTML}{d9534f}
\definecolor{quarto-callout-warning-color-frame}{HTML}{f0ad4e}
\definecolor{quarto-callout-tip-color-frame}{HTML}{02b875}
\definecolor{quarto-callout-caution-color-frame}{HTML}{fd7e14}
\makeatother
\makeatletter
\@ifpackageloaded{caption}{}{\usepackage{caption}}
\AtBeginDocument{%
\ifdefined\contentsname
  \renewcommand*\contentsname{Table of contents}
\else
  \newcommand\contentsname{Table of contents}
\fi
\ifdefined\listfigurename
  \renewcommand*\listfigurename{List of Figures}
\else
  \newcommand\listfigurename{List of Figures}
\fi
\ifdefined\listtablename
  \renewcommand*\listtablename{List of Tables}
\else
  \newcommand\listtablename{List of Tables}
\fi
\ifdefined\figurename
  \renewcommand*\figurename{Figure}
\else
  \newcommand\figurename{Figure}
\fi
\ifdefined\tablename
  \renewcommand*\tablename{Table}
\else
  \newcommand\tablename{Table}
\fi
}
\@ifpackageloaded{float}{}{\usepackage{float}}
\floatstyle{ruled}
\@ifundefined{c@chapter}{\newfloat{codelisting}{h}{lop}}{\newfloat{codelisting}{h}{lop}[chapter]}
\floatname{codelisting}{Listing}

\usepackage{amsthm}
\theoremstyle{plain}
\newtheorem{theorem}{Theorem}[section]
\theoremstyle{plain}
\newtheorem{lemma}{Lemma}[section]
\theoremstyle{plain}
\newtheorem{proposition}{Proposition}[section]
\theoremstyle{plain}
\newtheorem{corollary}{Corollary}[section]
\theoremstyle{remark}
\AtBeginDocument{}

\newtheorem{refremark}{Remark}[section]

\makeatother
\makeatletter
\@ifpackageloaded{caption}{}{\usepackage{caption}}
\@ifpackageloaded{subcaption}{}{\usepackage{subcaption}}
\makeatother
\usepackage{bookmark}
\IfFileExists{xurl.sty}{\usepackage{xurl}}{} % add URL line breaks if available
\hypersetup{
  pdftitle={Rao's score tests on correlation matrices},
  pdflang={en},
  pdfkeywords={correlation matrix, density power
divergence, equicorrelation, multivariate normal distribution, Rao's
score test, robustness, uncorrelatedness},
  colorlinks=true,
  linkcolor={blue},
  filecolor={Maroon},
  citecolor={Blue},
  urlcolor={Blue},
  pdfcreator={LaTeX via pandoc}}

\title{Rao's score tests on correlation matrices}
\author{Nirian Martín\textsuperscript{*} and Elena Almaraz\\[0.5em] \normalsize Department of Statistics and Operations Research and\\ \normalsize Interdisciplinary Mathematics Institute (IMI),\\ \normalsize Faculty of Mathematics, Complutense University of Madrid,\\ \normalsize Plaza de Ciencias 3, 28040 Madrid, Spain\\[0.3em] \normalsize \texttt{nimartin@ucm.es}, \texttt{ealmaraz@ucm.es}\\[0.3em] \small $^{*}$Corresponding author}
\date{December 2, 2026}
\begin{document}
\maketitle
\begin{abstract}
The Rao's score tests on correlation matrices of multivariate normal
distributions have received less attention than the likelihood ratio and
Wald tests, and their explicit expressions were not available for a
general correlation matrix. We derive the Rao's score test-statistics,
and their robust generalization, the Rao's \(\beta\)-score
test-statistics, for a specified correlation matrix, for fixed and
unknown equicorrelation, and for complete uncorrelatedness, in fixed
dimension. They are based on the restricted minimum density power
divergence estimators, in which a tuning parameter \(\beta\geq0\)
balances robustness and efficiency, the case \(\beta=0\) corresponding
to the classical Rao's score tests. All the test-statistics share a
simple trace expression in the correlation matrix estimated under the
null hypothesis, with an asymptotic chi-square distribution. We show
that this short form is valid for every \(\beta\geq0\), and that at
\(\beta=0\) the test-statistics are asymptotically equivalent to the
classical tests. The asymptotic power under contiguous alternatives,
with an explicit asymptotic relative efficiency, and the influence
functions of the estimators and of the tests are also obtained. A
simulation study and an application to the heptathlon data of the 1988
Olympic Games illustrate the robustness gained and provide a practical
rule for the choice of the tuning parameter.


\end{abstract}

\section{Introduction}\label{sec-Sec1}

Testing hypotheses on the correlation matrix of a multivariate normal
distribution is a classical problem. Bartlett (1954) gave the likelihood
ratio test of complete independence, Lawley (1963) a large-sample test
for the equality of all the correlation coefficients, and Jennrich
(1970) an asymptotic chi-square test for a specified correlation matrix,
whose test-statistic is algebraically equal to a Wald test-statistic
(Neudecker and Satorra, 1996). None of them is obtained as a Rao's score
test-statistic.

The Rao's score tests, introduced by C. R. Rao (1948), became popular in
Econometrics with an alternative version and refinements through
restrictions, called Lagrange multiplier tests, presented by Aitchison
and Silvey (1958) and Silvey (1959). In between both publications, Wald
(1943) had proposed another test with the same asymptotic distribution
as the likelihood ratio test and focussed on maximum likelihood
estimators (MLEs) too, but these estimators were not exactly the same.
While the likelihood ratio and Wald tests consider the whole parameter
space, under the null and alternative hypotheses, the parameter space of
the MLEs for the Rao's score test is only restricted to the null
hypothesis, so that the Rao's score test-statistic only requires the
estimators restricted to the null hypothesis.

In the same way as the classical likelihood ratio tests devoted to
correlation matrices, the Rao's score tests presented in Section
\ref{sec-Sec3} are based on the sample Pearson correlation matrix, which
is well known to be very susceptible to outliers or anomalous
observations in the data. The Rao's \(\beta\)-score tests (Basu et al.,
2022b) replace the MLEs by the minimum density power divergence
estimators (Basu et al., 1998), indexed by a single parameter
\(\beta>0\) which controls the trade-off between robustness and
efficiency. Values of \(\beta\) near zero retain efficiency but lose
robustness, and the classical Rao's score tests are recovered as
\(\beta\) tends to zero from the right, so that the classical tests
belong to the class rather than compete with it.

The contributions of the paper are the following. The Rao's
\(\beta\)-score test-statistics are derived for four hypotheses on the
correlation matrix of a multivariate normal distribution of any fixed
dimension \(p\), namely a specified correlation matrix, a fixed
equicorrelation, complete uncorrelatedness and an equicorrelation with
an unknown common coefficient, together with the restricted minimum
density power divergence estimators on which they are based, in closed
form whenever this is possible. The four test-statistics are shown to
share the form
\(\mathrm{trace}((\boldsymbol{R}_{0}^{-1}\boldsymbol{M}^{\ast})^{2})\),
where \(\boldsymbol{M}^{\ast}\) is an adjusted matrix defined by the
estimating equations of the nuisance parameters themselves, which
reduces four derivations to a single one and makes the comparison with
the classical tests immediate. On the inferential side, the regularity
conditions are established under which the restricted estimating
equations have a consistent sequence of roots, and this root is
identified; for \(\beta>0\) it is not the global minimizer of the
objective function, which is unbounded from below for small sample
sizes, a fact that does not seem to have been noticed for this family of
estimators. The asymptotic distribution of the test-statistics is
obtained under contiguous alternatives, and from it the asymptotic
relative efficiency with respect to the classical Rao's score tests in
closed form for general \(p\), which quantifies the price of robustness
and shows that it grows with the dimension. Finally, the influence
function of the restricted estimators and the second-order, level and
power influence functions of the tests are derived, all of them bounded
if and only if \(\beta>0\), with the gross-error sensitivity of the test
of complete uncorrelatedness in closed form. A simulation study examines
the finite-sample size and power of the test-statistics, verifies
numerically the asymptotic relative efficiency and shows that the value
of the tuning parameter which is needed depends on the type of
contamination, which gives a practical rule for its choice. Finally, the
theory is illustrated with a real data set in which a single observation
is responsible for most of the evidence reported by the classical
test-statistic.

For the particular case of a theoretical identity correlation matrix,
that is, for the test of independence, Leung and Drton (2018) derived
the corresponding classical Rao's score test by differentiating with
respect to the variance-covariance matrix, under the assumption that
this matrix is fixed except for the variances, which are estimated under
independence. Their technique gives the correct expression in that case,
but it relies on the restricted estimators of the variances being used
and on the particular form of the null hypothesis, and it need not give
the correct expression for a general theoretical correlation matrix. The
derivation of Appendix \ref{sec-ProofThTestGeneral}, in terms of a
change of the parameters rather than a transformation of the data, see
(\ref{eq-R-dif-par}), covers the general case, and the test-statistic of
Leung and Drton (2018) is recovered as the limit of ours as \(\beta\)
tends to zero from the right. As far as we are aware, for a general
theoretical correlation matrix neither the classical Rao's score
test-statistic nor the robust Rao's \(\beta\)-score test-statistic had
been derived before the first version of this paper (Martín, 2020). The
matrices \(\boldsymbol{J}_{\beta}(\boldsymbol{\theta})\) and
\(\boldsymbol{K}_{\beta}(\boldsymbol{\theta})\) of the minimum density
power divergence estimators of the multivariate normal distribution
obtained there have been used by Basu et al.~(2022a).

The rest of the paper is organized as follows. Section \ref{sec-Sec2}
covers the technical issues related to Rao's tests of any parameter
associated with the multivariate normal distribution. Notational aspects
are described in Subsection \ref{sec-Not}, and Subsections \ref{sec-BR1}
and \ref{sec-BR2} provide some basic results, valid for the specific
tests of the current paper as well as for any additional Rao's score
test associated with multivariate normal distributions, either for
correlation or for central or dispersion parameters, the second of them
justifying the form of the test-statistic used throughout. In Section
\ref{sec-Sec3}, Subsection \ref{sec-estim} explains how to compute the
restricted estimators, and Subsections \ref{sec-test1}, \ref{sec-test2},
\ref{sec-test3} and \ref{sec-test4} introduce each of the four
correlation tests. Section \ref{sec-Asym} establishes the regularity
conditions and the choice of the root of the estimating equations, the
asymptotic power under contiguous alternatives and the influence
functions of the restricted estimators and of the tests. Section
\ref{sec-Sim} reports a simulation study and Section \ref{sec-Appl}
illustrates the tests with a real data set, and the discussion of
Section \ref{sec-Discuss} relates the new tests with the classical tests
of Lawley (1963) and Jennrich (1970). The results appearing in Sections
\ref{sec-Sec2}, \ref{sec-Sec3} and \ref{sec-Asym} are proven in the
Appendix.

\section{Preliminary results}\label{sec-Sec2}

\subsection{Notational aspects}\label{sec-Not}

Let
\(\boldsymbol{X}=(X_{1},\ldots,X_{p})^{\top}\sim\mathcal{N}_{p} (\boldsymbol{\mu},\boldsymbol{\Sigma})\)
be a \(p\)-variate normal distribution, with
\(\boldsymbol{\mu}=(\mu_{1},\ldots,\mu_{p})^{\top}\),
\(\boldsymbol{\Sigma }=(\sigma_{ij})\), such that
\(\sigma_{ii}=\sigma_{i}^{2}=Var[X_{i}]\),
\(\sigma_{ij}=Cov[X_{i},Y_{j}]\). Let us consider the parameter vector
as
\(\boldsymbol{\theta}=(\boldsymbol{\mu}^{\top},\mathrm{vech}^{\top} (\boldsymbol{\Sigma}))^{\top}\),
where \(\mathrm{vech}\) is the so called vech operation for a symmetric
matrix. The vech of \(\boldsymbol{\Sigma}\) (vech for vector half) is
the \(p(p+1)/2\) dimensional vector obtained by stacking the unique part
of each column that lies on or below the diagonal of
\(\boldsymbol{\Sigma}\) into a single vector (for more details, see
Henderson \& Searle (1979)). For example, when \(p=2\) then
\(\mathrm{vech}(\boldsymbol{\Sigma })=(\sigma_{11},\sigma_{12},\sigma_{22})^{\top}=(\sigma_{1}^{2},\sigma _{12},\sigma_{2}^{2})^{\top}\).
In this setting, the density function is expressed as

\[
f_{\theta}(\boldsymbol{x})=\frac{1}{(2\pi)^{\frac{p}{2}}\left\vert
\boldsymbol{\Sigma}\right\vert ^{\frac{1}{2}}}\exp\left\{  -\frac{1}
{2}(\boldsymbol{x}-\boldsymbol{\mu})^{\top}\boldsymbol{\Sigma}^{-1}
(\boldsymbol{x}-\boldsymbol{\mu})\right\}  ,
\]

where \(\boldsymbol{x}=(x_{1},\ldots,x_{p})^{\top}\) is any point in the
support \(\mathbb{R}^{p}\) and
\(\boldsymbol{\theta}\in\mathbb{R}^{p}\times\mathbb{R}_{+}^{\frac{p(p+1)}{2}}\).

Let us consider the parameter vector as

\[
\begin{aligned}
\boldsymbol{\theta}  &  =(\boldsymbol{\theta}_{1}^{\top},\boldsymbol{\theta}
_{2}^{\top})^{\top}, \\
\boldsymbol{\theta}_{1}  &  =\boldsymbol{\mu}, \\
\boldsymbol{\theta}_{2}  &  =\mathrm{vech}(\boldsymbol{\Sigma}),
\end{aligned}
\]

but now the variance-covariance components are reordered and in addition
the correlation matrix is considered according to

\[
\boldsymbol{\Sigma}=\left(  \sigma_{ij}\right)  =\boldsymbol{\Lambda}
^{1/2}\boldsymbol{R\Lambda}^{1/2}=\left(  \sigma_{i}\rho_{ij}\sigma
_{j}\right)  ,
\]

where

\[
\begin{aligned}
\boldsymbol{\Lambda}  &  =\mathrm{diag}\{\sigma_{j}^{2}\}_{j=1}^{p}, \\
\boldsymbol{R}  &  =(\rho_{ij}).
\end{aligned}
\]

The first parameter vector variable change is

\[
\begin{aligned}
\boldsymbol{\phi}  &  =(\boldsymbol{\phi}_{1}^{\top},\boldsymbol{\phi}_{2}
^{\top})^{\top}, \\
\boldsymbol{\phi}_{1}  &  =\boldsymbol{\mu}=\boldsymbol{\theta}_{1}, \\
\boldsymbol{\phi}_{2}  &  =((\boldsymbol{\Lambda1}_{p})^{\top},\mathrm{vecl}
^{\top}(\boldsymbol{\Sigma}))^{\top}=\boldsymbol{M}^{\top}\boldsymbol{\theta}_{2},
\end{aligned}
\]

where \(\mathrm{vecl}(\cdot)\) denotes the vectorization operator of the
lower off-diagonal elements of matrix \(\cdot\) (unlike the
\(\mathrm{vech}(\cdot)\) operator, this operator excludes the diagonal
elements) and

\[
\begin{aligned}
\boldsymbol{M}  &  =\left(  \boldsymbol{P},\boldsymbol{Q}\right)  , \\
\boldsymbol{P}  &  =(\boldsymbol{P}_{\cdot1},..,\boldsymbol{P}_{\cdot p}), \\
\boldsymbol{P}_{\cdot i}  &  =\boldsymbol{e}_{(i-1)(p+1)-\frac{i(i-1)}{2}
+1}, \\
\boldsymbol{Q}  &  =(\boldsymbol{Q}_{\cdot1},..,\boldsymbol{Q}_{\cdot
(p-1)p/2}), \\
\boldsymbol{Q}_{\cdot s}  &  =\boldsymbol{e}_{(i-1)(p+1)-\frac{i(i-1)}
{2}+1+j-i},
\end{aligned}
\]

\(i<j\) in lexicographical order. Notice that \(\boldsymbol{M}\) is a
permutation matrix and so it is orthogonal,
i.e.~\(\boldsymbol{M}^{-1}=\boldsymbol{M}^{\top}\). Let

\[
\boldsymbol{\phi}_{2}=(\boldsymbol{\phi}_{2,1}^{\top},\boldsymbol{\phi}_{2,2}
^{\top})^{\top}=((\boldsymbol{\Lambda1}_{p})^{\top},\mathrm{vecl}^{\top}
(\boldsymbol{\Sigma}))^{\top},
\]

\[
\boldsymbol{\theta}_{2}=\mathrm{vech}(\boldsymbol{\Sigma}
),\;\boldsymbol{\Sigma}=\boldsymbol{\Sigma}(\boldsymbol{\Lambda}
)=\boldsymbol{\Lambda}^{1/2}\boldsymbol{R\Lambda}^{1/2}
\]

be a partition of \(\boldsymbol{\phi}_{2}\), then from previous
expressions it is concluded that

\[
\begin{aligned}
\boldsymbol{\phi}_{2,1}  &  =\boldsymbol{P}^{\top}\boldsymbol{\theta}_{2}, \\
\boldsymbol{\phi}_{2,2}  &  =\boldsymbol{Q}^{\top}\boldsymbol{\theta}_{2}.
\end{aligned}
\]

Ferrari and Yang (2010) used the previous parameter scheme for
estimation though the so-called \(Lq\)-estimators~of
\(\boldsymbol{\phi}\). The second parameter change is

\[
\begin{aligned}
\boldsymbol{\eta}  &  =(\boldsymbol{\eta}_{1}^{\top},\boldsymbol{\eta}_{2}
^{\top})^{\top}, \\
\boldsymbol{\eta}_{2}  &  =(\boldsymbol{\eta}_{2,1}^{\top},\boldsymbol{\eta
}_{2,2}^{\top})^{\top}, \\
\boldsymbol{\eta}_{2,1}  &  =\boldsymbol{\phi}_{2,1}=\boldsymbol{\Lambda1}
_{p}, \\
\boldsymbol{\eta}_{2,2}  &  =\mathrm{vecl}(\boldsymbol{R})=\mathrm{vecl}
(\boldsymbol{\Lambda}^{-1/2}\boldsymbol{\Sigma\Lambda}^{-1/2}).
\end{aligned}
\]

For example, if \(p=3\) then
\(\boldsymbol{\theta}=(\mu_{1},\mu_{2},\mu _{3},\sigma_{1}^{2},\sigma_{12},\sigma_{13},\sigma_{2}^{2},\sigma_{23} ,\sigma_{3}^{2})^{\top}\),
\(\boldsymbol{\phi}=(\mu_{1},\mu_{2},\mu_{3},\sigma_{1}^{2},\sigma_{2}^{2},\sigma_{3}^{2},\allowbreak\sigma_{12},\sigma_{13} ,\sigma_{23})^{\top}\)
and
\(\boldsymbol{\eta}=(\mu_{1},\mu_{2},\mu_{3},\sigma_{1}^{2},\sigma_{2}^{2},\sigma_{3}^{2},\rho_{12},\rho_{13},\rho_{23})^{\top}\).

\subsection{\texorpdfstring{Basic results for the partition of parameter
vector
\(\boldsymbol{\theta}\)}{Basic results for the partition of parameter vector \textbackslash boldsymbol\{\textbackslash theta\}}}\label{sec-BR1}

\begin{proposition}[]\protect\hypertarget{prp-Props}{}\label{prp-Props}

The expression of the vectorial score function,
\(\boldsymbol{s} _{\boldsymbol{\theta}}(\boldsymbol{x})=(\boldsymbol{s}_{\boldsymbol{\mu}} ^{\top}(\boldsymbol{x}),\boldsymbol{s}_{\mathrm{vech}(\boldsymbol{\Sigma})} ^{\top}(\boldsymbol{x}))^{\top}\),
is given as follows

\begin{equation}\protect\phantomsection\label{eq-scMu}{
\begin{aligned}
\boldsymbol{s}_{\boldsymbol{\mu}}(\boldsymbol{x})  &  =\frac{\partial
}{\partial\boldsymbol{\mu}}\log f_{\boldsymbol{\theta}}(\boldsymbol{x}
) \\
&  =\boldsymbol{\Sigma}^{-1}(\boldsymbol{x}-\boldsymbol{\mu}),
\end{aligned}
}\end{equation}

\begin{equation}\protect\phantomsection\label{eq-scSig}{
\begin{aligned}
\boldsymbol{s}_{\mathrm{vech}(\boldsymbol{\Sigma})}(\boldsymbol{x})  &
=-\frac{1}{2}\mathrm{vech}\left(  \frac{\partial}{\partial\boldsymbol{\Sigma}
}\left(  \log\left\vert \boldsymbol{\Sigma}\right\vert +(\boldsymbol{x}
-\boldsymbol{\mu})^{\top}\boldsymbol{\Sigma}^{-1}(\boldsymbol{x}-\boldsymbol{\mu
})\right)  \right) \\
&  =-\frac{1}{2}\boldsymbol{G}_{p}^{\top}\mathrm{vec}(\boldsymbol{\Sigma}
^{-1})+\frac{1}{2}\boldsymbol{G}_{p}^{\top}\left[  \left(  \boldsymbol{\Sigma
}^{-1}(\boldsymbol{x}-\boldsymbol{\mu})\right)  \otimes\left(  (\boldsymbol{x}
-\boldsymbol{\mu})^{\top}\boldsymbol{\Sigma}^{-1}\right)  \right]  ,
\end{aligned}
}\end{equation}

with \(\boldsymbol{x}\in\mathbb{R}^{p}\) being any point of the support
and \(\boldsymbol{G}_{p}\) the so-called ``duplication matrix'' of order
\(p\), i.e., the unique \(p^{2}\times\frac{p(p+1)}{2}\)~matrix such that
\(\mathrm{vec} (\boldsymbol{\Sigma})=\boldsymbol{G}_{p}\mathrm{vech}(\boldsymbol{\Sigma})\).

\end{proposition}

\begin{theorem}[]\protect\hypertarget{thm-Th-J-BS1}{}\label{thm-Th-J-BS1}

The expression of
\(\boldsymbol{J}_{\beta}(\boldsymbol{\theta })=E_{\boldsymbol{\theta}}[\boldsymbol{s}_{\boldsymbol{\theta}}(\boldsymbol{X} )\boldsymbol{s}_{\boldsymbol{\theta}}^{\top}(\boldsymbol{X})f_{\boldsymbol{\theta }}^{\beta}]\),
is given as follows

\begin{equation}\protect\phantomsection\label{eq-eqJ}{
\begin{aligned}
\boldsymbol{J}_{\beta}(\boldsymbol{\theta})=
\begin{pmatrix}
\boldsymbol{J}_{\beta}(\boldsymbol{\mu}) & \boldsymbol{0}_{p\times
\frac{p(p+1)}{2}}\\
\boldsymbol{0}_{\frac{p(p+1)}{2}\times p} & \boldsymbol{J}_{\beta
}(\mathrm{vech}(\boldsymbol{\Sigma}))
\end{pmatrix}
,
\end{aligned}
}\end{equation}

where

\begin{equation}\protect\phantomsection\label{eq-eqJ11}{
\begin{aligned}
\boldsymbol{J}_{\beta}(\boldsymbol{\mu})  &  =\frac{(\beta+1)^{-\frac{p}{2}
-1}}{(2\pi)^{\frac{\beta p}{2}}\left\vert \boldsymbol{\Sigma}\right\vert
^{\frac{\beta}{2}}}\boldsymbol{\Sigma}^{-1},
\end{aligned}
}\end{equation}

\begin{equation}\protect\phantomsection\label{eq-eqJ22}{
\begin{aligned}
\boldsymbol{J}_{\beta}(\mathrm{vech}(\boldsymbol{\Sigma}))  &  =\frac
{(\beta+1)^{-\frac{p}{2}-2}}{4(2\pi)^{\frac{\beta p}{2}}\left\vert
\boldsymbol{\Sigma}\right\vert ^{\frac{\beta}{2}}}\left[  \beta^{2}
\boldsymbol{C}_{\boldsymbol{\theta}}+2\boldsymbol{G}_{p}^{\top}\left(
\boldsymbol{\Sigma}^{-1}\otimes\boldsymbol{\Sigma}^{-1}\right)  \boldsymbol{G}
_{p}\right]  ,
\end{aligned}
}\end{equation}

and

\begin{equation}\protect\phantomsection\label{eq-C}{
\boldsymbol{C}_{\boldsymbol{\theta}}=\boldsymbol{G}_{p}^{\top}\mathrm{vec}\left(
\boldsymbol{\Sigma}^{-1}\right)  \mathrm{vec}^{\top}\left(  \boldsymbol{\Sigma
}^{-1}\right)  \boldsymbol{G}_{p}.
}\end{equation}

\end{theorem}

\begin{theorem}[]\protect\hypertarget{thm-Th-Psi-BS1}{}\label{thm-Th-Psi-BS1}

The expression of \(\boldsymbol{\xi}_{\beta} (\boldsymbol{\theta})\) and
\(\boldsymbol{\xi}_{\beta}(\boldsymbol{\theta })\boldsymbol{\xi}_{\beta}^{\top}(\boldsymbol{\theta})\),where
\(\boldsymbol{\xi }_{\beta}(\boldsymbol{\theta})=E_{\boldsymbol{\theta}}[\boldsymbol{s} _{\boldsymbol{\theta}}(\boldsymbol{X})f_{\boldsymbol{\theta}}^{\beta }(\boldsymbol{X})]\),
are given as follows

\begin{equation}\protect\phantomsection\label{eq-psipsi}{
\begin{aligned}
\boldsymbol{\xi}_{\beta}(\boldsymbol{\theta})  &  =
\begin{pmatrix}
\boldsymbol{\xi}_{\beta}(\boldsymbol{\mu}) \\
\boldsymbol{\xi}_{\beta}(\mathrm{vech}(\boldsymbol{\Sigma}))
\end{pmatrix}
=
\begin{pmatrix}
\boldsymbol{0}_{p} \\
\boldsymbol{\xi}_{\beta}(\mathrm{vech}(\boldsymbol{\Sigma}))
\end{pmatrix}
, \\
\boldsymbol{\xi}_{\beta}(\boldsymbol{\theta})\boldsymbol{\xi}_{\beta}
^{\top}(\boldsymbol{\theta})  &  =
\begin{pmatrix}
\boldsymbol{0}_{p\times p} & \boldsymbol{0}_{p\times\frac{p(p+1)}{2}} \\
\boldsymbol{0}_{\frac{p(p+1)}{2}\times p} & \boldsymbol{\xi}_{\beta
}(\mathrm{vech}(\boldsymbol{\Sigma}))\boldsymbol{\xi}_{\beta}^{\top}
(\mathrm{vech}(\boldsymbol{\Sigma}))
\end{pmatrix}
\end{aligned}
}\end{equation}

where

\begin{equation}\protect\phantomsection\label{eq-psipsi2}{
\begin{aligned}
\boldsymbol{\xi}_{\beta}(\mathrm{vech}(\boldsymbol{\Sigma}))  &  =-\frac
{\beta}{2}\frac{(\beta+1)^{-(\frac{p}{2}+1)}}{(2\pi)^{\frac{\beta p}{2}
}\left\vert \boldsymbol{\Sigma}\right\vert ^{\frac{\beta}{2}}}\boldsymbol{G}
_{p}^{\top}\mathrm{vec}\left(  \boldsymbol{\Sigma}^{-1}\right)  , \\
\boldsymbol{\xi}_{\beta}(\mathrm{vech}(\boldsymbol{\Sigma}))\boldsymbol{\xi
}_{\beta}^{\top}(\mathrm{vech}(\boldsymbol{\Sigma}))  &  =\frac{\beta^{2}}
{4}\frac{(\beta+1)^{-(p+2)}}{(2\pi)^{\beta p}\left\vert \boldsymbol{\Sigma
}\right\vert ^{\beta}}\boldsymbol{C}_{\boldsymbol{\theta}}.
\end{aligned}
}\end{equation}

\end{theorem}

\begin{corollary}[]\protect\hypertarget{cor-Th-Psi-KS1}{}\label{cor-Th-Psi-KS1}

The expression of \(\boldsymbol{K}_{\beta}(\boldsymbol{\theta })\) is
given as follows

\[
\begin{aligned}
\boldsymbol{K}_{\beta}(\boldsymbol{\theta})=
\begin{pmatrix}
\boldsymbol{K}_{\beta}(\boldsymbol{\mu}) & \boldsymbol{0}_{p\times
\frac{p(p+1)}{2}}\\
\boldsymbol{0}_{\frac{p(p+1)}{2}\times p} & \boldsymbol{K}_{\beta
}(\mathrm{vech}(\boldsymbol{\Sigma}))
\end{pmatrix}
,
\end{aligned}
\]

where

\[
\boldsymbol{K}_{\beta}(\boldsymbol{\mu})=\boldsymbol{J}_{2\beta}
(\boldsymbol{\mu})=\frac{(2\beta+1)^{-\frac{p}{2}-1}}{(2\pi)^{\beta
p}\left\vert \boldsymbol{\Sigma}\right\vert ^{\beta}}\boldsymbol{\Sigma}
^{-1},
\]

\[
\boldsymbol{K}_{\beta}(\mathrm{vech}(\boldsymbol{\Sigma}))=\frac{1}
{4(2\pi)^{\beta p}\left\vert \boldsymbol{\Sigma}\right\vert ^{\beta}
}\boldsymbol{G}_{p}^{\top}\left[  \overline{\boldsymbol{J}}_{2\beta
}(\boldsymbol{\Sigma}^{-1})+\overline{\boldsymbol{\xi}}_{\beta}
(\boldsymbol{\Sigma}^{-1})\overline{\boldsymbol{\xi}}_{\beta}^{\top}
(\boldsymbol{\Sigma}^{-1})\right]  \boldsymbol{G}_{p},
\]

with

\begin{equation}\protect\phantomsection\label{eq-Kappa1}{
\begin{aligned}
\overline{\boldsymbol{J}}_{2\beta}(\boldsymbol{\Sigma}^{-1})  &  =\kappa
_{1}(p,\beta)\left(  \boldsymbol{\Sigma}^{-1}\otimes\boldsymbol{\Sigma}
^{-1}\right)  , \\
\kappa_{1}(p,\beta)  &  =2(2\beta+1)^{-\frac{p}{2}-2}, \\
\overline{\boldsymbol{\xi}}_{\beta}(\boldsymbol{\Sigma}^{-1})\overline
{\boldsymbol{\xi}}_{\beta}^{\top}(\boldsymbol{\Sigma}^{-1})  &  =\kappa
_{2}(p,\beta)\mathrm{vec}(\boldsymbol{\Sigma}^{-1})\mathrm{vec}^{\top}
(\boldsymbol{\Sigma}^{-1}) \\
\kappa_{2}(p,\beta)  &  =\beta^{2}\left[  4(2\beta+1)^{-\frac{p}{2}-2}
-(\beta+1)^{-(p+2)}\right]  .
\end{aligned}
}\end{equation}

\end{corollary}

\begin{proposition}[]\protect\hypertarget{prp-propK}{}\label{prp-propK}

The expression of
\(\boldsymbol{K}_{\beta}(\mathrm{vech} (\boldsymbol{\Sigma}))\) and its
inverse in terms of \(\boldsymbol{\Lambda}\) and \(\boldsymbol{R}_{0}\)
are given by

\[
\begin{aligned}
&  \boldsymbol{K}_{\beta}(\mathrm{vech}(\boldsymbol{\Lambda}^{\frac{1}{2}
}\boldsymbol{R}_{0}\boldsymbol{\Lambda}^{\frac{1}{2}}))=\frac{1}
{4(2\pi)^{\beta p}\left\vert \boldsymbol{\Lambda}\right\vert ^{\beta
}\left\vert \boldsymbol{R}_{0}\right\vert ^{\beta}} \\
&  \times\boldsymbol{G}_{p}^{\top}\left(  \boldsymbol{\Lambda}^{-\frac{1}{2}
}\otimes\boldsymbol{\Lambda}^{-\frac{1}{2}}\right)  \left(  \overline
{\boldsymbol{J}}_{2\beta}(\boldsymbol{R}_{0}^{-1})+\overline{\boldsymbol{\xi}
}_{\beta}(\boldsymbol{R}_{0}^{-1})\overline{\boldsymbol{\xi}}_{\beta}
^{\top}(\boldsymbol{R}_{0}^{-1})\right)  \left(  \boldsymbol{\Lambda}^{-\frac
{1}{2}}\otimes\boldsymbol{\Lambda}^{-\frac{1}{2}}\right)  \boldsymbol{G}_{p},
\end{aligned}
\]

with

\[
\begin{aligned}
\overline{\boldsymbol{J}}_{2\beta}(\boldsymbol{R}_{0}^{-1})  &  =\kappa
_{1}(p,\beta)\left(  \boldsymbol{R}_{0}^{-1}\otimes\boldsymbol{R}_{0}
^{-1}\right)  , \\
\overline{\boldsymbol{\xi}}_{\beta}(\boldsymbol{R}_{0}^{-1})\overline
{\boldsymbol{\xi}}_{\beta}^{\top}(\boldsymbol{R}_{0}^{-1})  &  =\kappa
_{2}(p,\beta)\mathrm{vec}(\boldsymbol{R}_{0}^{-1})\mathrm{vec}^{\top}
(\boldsymbol{R}_{0}^{-1}),
\end{aligned}
\]

and

\[
\begin{aligned}
&  \boldsymbol{K}_{\beta}^{-1}(\mathrm{vech}(\boldsymbol{\Lambda}^{\frac{1}
{2}}\boldsymbol{R}_{0}\boldsymbol{\Lambda}^{\frac{1}{2}}))=4(2\pi)^{\beta
p}\left\vert \boldsymbol{\Lambda}\right\vert ^{\beta}\left\vert \boldsymbol{R}
_{0}\right\vert ^{\beta} \\
&  \times\boldsymbol{G}_{p}^{+}\left(  \boldsymbol{\Lambda}^{\frac{1}{2}}
\otimes\boldsymbol{\Lambda}^{\frac{1}{2}}\right)  \left(  \overline
{\boldsymbol{J}}_{2\beta}(\boldsymbol{R}_{0}^{-1})+\overline{\boldsymbol{\xi}
}_{\beta}(\boldsymbol{R}_{0}^{-1})\overline{\boldsymbol{\xi}}_{\beta}
^{\top}(\boldsymbol{R}_{0}^{-1})\right)  ^{-1}\left(  \boldsymbol{\Lambda}
^{\frac{1}{2}}\otimes\boldsymbol{\Lambda}^{\frac{1}{2}}\right)  \boldsymbol{L}
_{p}^{\top},
\end{aligned}
\]

with

\begin{equation}\protect\phantomsection\label{eq-inv}{
\left(  \overline{\boldsymbol{J}}_{2\beta}(\boldsymbol{R}_{0}^{-1}
)+\overline{\boldsymbol{\xi}}_{\beta}(\boldsymbol{R}_{0}^{-1})\overline
{\boldsymbol{\xi}}_{\beta}^{\top}(\boldsymbol{R}_{0}^{-1})\right)  ^{-1}
=\kappa_{1}^{-1}(p,\beta)\left(  \left(  \boldsymbol{R}_{0}\otimes
\boldsymbol{R}_{0}\right)  -\frac{\kappa_{3}(p,\beta)\mathrm{vec}
(\boldsymbol{R}_{0})\mathrm{vec}^{\top}(\boldsymbol{R}_{0})}{1+p\kappa
_{3}(p,\beta)}\right)  ,
}\end{equation}

where

\begin{equation}\protect\phantomsection\label{eq-kappa3}{
\kappa_{3}(p,\beta)=\kappa_{1}^{-1}(p,\beta)\kappa_{2}(p,\beta).
}\end{equation}

\end{proposition}

\begin{theorem}[]\protect\hypertarget{thm-Th-U-BS1}{}\label{thm-Th-U-BS1}

The expression of the \(\beta\)-score statistic,

\[
\begin{aligned}
\boldsymbol{U}_{\beta,n}\left(  \boldsymbol{\theta}\right)   &  =\frac{1}
{n}\sum_{i=1}^{n}\boldsymbol{u}_{\beta}\left(  \boldsymbol{X}_{i}
,\boldsymbol{\theta}\right)  =(\boldsymbol{U}_{\beta,n}^{\top}\left(
\boldsymbol{\mu}\right)  ,\boldsymbol{U}_{\beta,n}^{\top}\left(  \mathrm{vech}(\boldsymbol{\Sigma
})\right)  )^{\top}, \\
\boldsymbol{u}_{\beta}\left(  \boldsymbol{x},\boldsymbol{\theta}\right)   &
=\boldsymbol{s}_{\boldsymbol{\theta}}(\boldsymbol{x})f_{\boldsymbol{\theta}}^{\beta
}(\boldsymbol{x})-\boldsymbol{\xi}_{\beta}(\boldsymbol{\theta}),
\end{aligned}
\]

is given by

\[
\boldsymbol{U}_{\beta,n}\left(  \boldsymbol{\mu}\right)  =\frac{1}
{(2\pi)^{\frac{\beta p}{2}}\left\vert \boldsymbol{\Sigma}\right\vert
^{\frac{\beta}{2}}}\boldsymbol{\Sigma}^{-1}\frac{1}{n}\sum_{i=1}^{n}
W_{i,\beta}(\boldsymbol{\theta})(\boldsymbol{X}_{i}-\boldsymbol{\mu}),
\]

with

\begin{equation}\protect\phantomsection\label{eq-w}{
W_{i,\beta}(\boldsymbol{\theta})=\exp\left\{  -\frac{\beta}{2}(\boldsymbol{X}
_{i}-\boldsymbol{\mu})^{\top}\boldsymbol{\Sigma}^{-1}(\boldsymbol{X}
_{i}-\boldsymbol{\mu})\right\}  ,
}\end{equation}

and

\[
\boldsymbol{U}_{\beta,n}\left(  \mathrm{vech}(\boldsymbol{\Sigma})\right)
=\boldsymbol{G}_{p}^{\top}\boldsymbol{V}_{\beta,n}\left(  \mathrm{vec}
(\boldsymbol{\Sigma})\right)
\]

with

\begin{equation}\protect\phantomsection\label{eq-vecU}{
\begin{aligned}
\boldsymbol{V}_{\beta,n}\left(  \mathrm{vec}(\boldsymbol{\Sigma})\right)
=&\frac{1}{2(2\pi)^{\frac{\beta p}{2}}\left\vert \boldsymbol{\Sigma
}\right\vert ^{\frac{\beta}{2}}}\Bigg[\frac{1}{n}\sum_{i=1}^{n}W_{i,\beta}(\boldsymbol{\theta})\left(
\boldsymbol{\Sigma}^{-1}(\boldsymbol{X}_{i}-\boldsymbol{\mu})\right)
\otimes\left(  \boldsymbol{\Sigma}^{-1}(\boldsymbol{X}_{i}-\boldsymbol{\mu
})\right) \\
&-\left(  \frac{1}{n}\sum_{i=1}^{n}W_{i,\beta}(\boldsymbol{\theta
})-\beta(\beta+1)^{-(\frac{p}{2}+1)}\right)  \mathrm{vec}\left(
\boldsymbol{\Sigma}^{-1}\right)\Bigg].
\end{aligned}
}\end{equation}

\end{theorem}

\begin{proposition}[]\protect\hypertarget{prp-PropV}{}\label{prp-PropV}

Let
\(\boldsymbol{V}_{\beta,n}(\widetilde{\boldsymbol{\Lambda} }_{\beta},\boldsymbol{R}_{0})\)
denote
\(\boldsymbol{V}_{\beta,n} (\widetilde{\boldsymbol{\theta}}_{2,\beta})\),
given in (\ref{eq-vecU}), a term of
\(\boldsymbol{U}_{\beta,n}(\widetilde{\boldsymbol{\theta}}_{2,\beta})\)
according to (\ref{eq-vecU0}), particularized to
\(\widetilde{\boldsymbol{\Sigma} }_{\beta}\boldsymbol{=}\widetilde{\boldsymbol{\Lambda}}_{\beta}^{1/2} \boldsymbol{R}_{0}\widetilde{\boldsymbol{\Lambda}}_{\beta}^{1/2}\),
then

\begin{equation}\protect\phantomsection\label{eq-V}{
\boldsymbol{V}_{\beta,n}(\widetilde{\boldsymbol{\Lambda}}_{\beta
},\boldsymbol{R}_{0})=\frac{\widetilde{\kappa}_{0}(p,\beta)}{2(2\pi
)^{\frac{\beta p}{2}}\left\vert \widetilde{\boldsymbol{\Sigma}}_{\beta
}\right\vert ^{\frac{\beta}{2}}}(\widetilde{\boldsymbol{\Lambda}}_{\beta
}^{-1/2}\otimes\widetilde{\boldsymbol{\Lambda}}_{\beta}^{-1/2})(\boldsymbol{R}
_{0}^{-1}\otimes\boldsymbol{R}_{0}^{-1})\left[  \mathrm{vec}\left(
\widetilde{\boldsymbol{R}}_{\boldsymbol{X},\beta}\right)  -\mathrm{vec}\left(
\boldsymbol{R}_{0}\right)  \right]  ,
}\end{equation}

where
\(\widetilde{\kappa}_{0}(p,\beta)=\frac{1}{n}\sum_{i=1}^{n}\widetilde{W}_{i,\beta}-\beta(\beta+1)^{-(\frac{p}{2}+1)}\),
which is the constant (\ref{eq-kappa0}) of Section \ref{sec-estim}, and

\[
\mathrm{vec}(\widetilde{\boldsymbol{R}}_{\boldsymbol{X},\beta})=\frac{\frac
{1}{n}\sum_{i=1}^{n}\widetilde{W}_{i,\beta}\left(
\widetilde{\boldsymbol{\Lambda}}_{\beta}^{-1/2}(\boldsymbol{X}_{i}
-\widetilde{\boldsymbol{\mu}}_{\beta})\right)  \otimes\left(
\widetilde{\boldsymbol{\Lambda}}_{\beta}^{-1/2}(\boldsymbol{X}_{i}
-\widetilde{\boldsymbol{\mu}}_{\beta})\right)  }{\widetilde{\kappa}
_{0}(p,\beta)}.
\]

\end{proposition}

\subsection{Validity of the short form of the
test-statistic}\label{sec-BR2}

For the hypotheses of this paper the parameter vector
\(\boldsymbol{\eta}\) is split into the nuisance parameters, estimated
under the null hypothesis, and the correlation parameters whose value is
fixed by the null hypothesis. Accordingly, let
\(\boldsymbol{J}_{\beta,ij}(\boldsymbol{\eta})\) and
\(\boldsymbol{K}_{\beta,ij}(\boldsymbol{\eta})\), \(i,j=1,2\), be the
blocks of \(\boldsymbol{J}_{\beta}(\boldsymbol{\eta})\) and
\(\boldsymbol{K}_{\beta}(\boldsymbol{\eta})\), and
\(\boldsymbol{U}_{\beta,n,1}(\boldsymbol{\eta})\) and
\(\boldsymbol{U}_{\beta,n,2}(\boldsymbol{\eta})\) the corresponding
components of the \(\beta\)-score function, for this partition. At the
restricted minimum DPD estimator
\(\widetilde{\boldsymbol{\eta}}_{\beta}\) the \(\beta\)-score function
satisfies
\(\boldsymbol{U}_{\beta,n,1}(\widetilde{\boldsymbol{\eta}}_{\beta})=\boldsymbol{0}\),
and the asymptotic variance-covariance matrix of
\(\sqrt{n}\boldsymbol{U}_{\beta,n,2}(\widetilde{\boldsymbol{\eta}}_{\beta})\)
is, in general,
\begin{equation}\protect\phantomsection\label{eq-Vgeneral}{
\begin{aligned}
\boldsymbol{V}_{\beta}(\boldsymbol{\eta})=&\boldsymbol{K}_{\beta,22}(\boldsymbol{\eta})
-\boldsymbol{J}_{\beta,21}(\boldsymbol{\eta})\boldsymbol{J}_{\beta,11}^{-1}(\boldsymbol{\eta})\boldsymbol{K}_{\beta,12}(\boldsymbol{\eta})
-\boldsymbol{K}_{\beta,21}(\boldsymbol{\eta})\boldsymbol{J}_{\beta,11}^{-1}(\boldsymbol{\eta})\boldsymbol{J}_{\beta,12}(\boldsymbol{\eta})\\
&+\boldsymbol{J}_{\beta,21}(\boldsymbol{\eta})\boldsymbol{J}_{\beta,11}^{-1}(\boldsymbol{\eta})\boldsymbol{K}_{\beta,11}(\boldsymbol{\eta})\boldsymbol{J}_{\beta,11}^{-1}(\boldsymbol{\eta})\boldsymbol{J}_{\beta,12}(\boldsymbol{\eta}),
\end{aligned}
}\end{equation} see Basu et al.~(2022b, Section 4). The short form
\(n\boldsymbol{U}_{\beta,n}^{\top}(\widetilde{\boldsymbol{\eta}}_{\beta})\boldsymbol{K}_{\beta}^{-1}(\widetilde{\boldsymbol{\eta}}_{\beta})\boldsymbol{U}_{\beta,n}(\widetilde{\boldsymbol{\eta}}_{\beta})\)
uses instead the inverse of
\(\boldsymbol{K}_{\beta,22}(\boldsymbol{\eta})-\boldsymbol{K}_{\beta,21}(\boldsymbol{\eta})\boldsymbol{K}_{\beta,11}^{-1}(\boldsymbol{\eta})\boldsymbol{K}_{\beta,12}(\boldsymbol{\eta})\),
which coincides with (\ref{eq-Vgeneral}) when
\(\boldsymbol{J}_{\beta,21}(\boldsymbol{\eta})\boldsymbol{J}_{\beta,11}^{-1}(\boldsymbol{\eta})=\boldsymbol{K}_{\beta,21}(\boldsymbol{\eta})\boldsymbol{K}_{\beta,11}^{-1}(\boldsymbol{\eta})\).
This always holds for \(\beta=0\), since
\(\boldsymbol{K}_{0}(\boldsymbol{\eta})=\boldsymbol{J}_{0}(\boldsymbol{\eta})\),
but not in general for \(\beta>0\). The next lemma shows that it holds
for the hypotheses considered in this paper.

\begin{lemma}[]\protect\hypertarget{lem-short}{}\label{lem-short}

Let the nuisance parameters be
\(\boldsymbol{\eta}_{1}=\boldsymbol{\mu}\), the variances
\(\boldsymbol{\eta}_{2,1}=\boldsymbol{\Lambda}\boldsymbol{1}_{p}\) and,
for the test of non-fixed equicorrelation, the common correlation
coefficient as well, and let the parameters fixed by the null hypothesis
be the remaining correlation parameters, a subvector of
\(\boldsymbol{\eta}_{2,2}=\mathrm{vecl}(\boldsymbol{R})\). Then, for
every \(\beta\geq0\), \[
\boldsymbol{J}_{\beta,21}(\boldsymbol{\eta})\boldsymbol{J}_{\beta,11}^{-1}(\boldsymbol{\eta})=\boldsymbol{K}_{\beta,21}(\boldsymbol{\eta})\boldsymbol{K}_{\beta,11}^{-1}(\boldsymbol{\eta}),
\] and the Rao's \(\beta\)-score test-statistic is \[
n\boldsymbol{U}_{\beta,n}^{\top}(\widetilde{\boldsymbol{\eta}}_{\beta})\boldsymbol{K}_{\beta}^{-1}(\widetilde{\boldsymbol{\eta}}_{\beta})\boldsymbol{U}_{\beta,n}(\widetilde{\boldsymbol{\eta}}_{\beta}).
\] Since the quadratic form is invariant under one-to-one differentiable
changes of the parameter, it is also equal to
\(n\boldsymbol{U}_{\beta,n}^{\top}(\widetilde{\boldsymbol{\theta}}_{\beta})\boldsymbol{K}_{\beta}^{-1}(\widetilde{\boldsymbol{\theta}}_{\beta})\boldsymbol{U}_{\beta,n}(\widetilde{\boldsymbol{\theta}}_{\beta})\),
with \(\widetilde{\boldsymbol{\theta}}_{\beta}\) the same estimator
expressed in the parameter \(\boldsymbol{\theta}\).

\end{lemma}

\section{Main results}\label{sec-Sec3}

\subsection{Restricted maximum likelihood estimators (MLEs) and minimum
density power divergences (MDPDs)}\label{sec-estim}

\begin{theorem}[]\protect\hypertarget{thm-RestrMDPDs}{}\label{thm-RestrMDPDs}

For known correlation matrix, \(\boldsymbol{R=R}_{0}\), the restricted
minimum DPD estimators of
\((\boldsymbol{\mu}^{\top},\boldsymbol{1} _{p}^{\top}\boldsymbol{\Lambda})^{\top}\),
are obtained as solution in
\((\widetilde{\boldsymbol{\mu}}_{\beta}^{\top},\boldsymbol{1}_{p}^{\top} \widetilde{\boldsymbol{\Lambda}}_{\beta})^{\top}\)
of

\begin{equation}\protect\phantomsection\label{eq-eqMuBeta}{
\begin{aligned}
\widetilde{\boldsymbol{\mu}}_{\beta}  &  =\frac{\sum\limits_{i=1}^{n}
W_{i}(\widetilde{\boldsymbol{\mu}}_{\beta},\widetilde{\boldsymbol{\Lambda}
}_{\beta})\boldsymbol{X}_{i}}{\sum\limits_{i=1}^{n}W_{i}
(\widetilde{\boldsymbol{\mu}}_{\beta},\widetilde{\boldsymbol{\Lambda}}_{\beta
})},
\end{aligned}
}\end{equation}

\begin{equation}\protect\phantomsection\label{eq-eqRRBeta}{
\begin{aligned}
\boldsymbol{1}_{p}  &  =\mathrm{diag}\{\boldsymbol{R}_{0}^{-1}\boldsymbol{R}
_{\boldsymbol{X},\beta}(\widetilde{\boldsymbol{\mu}}_{\beta}
,\widetilde{\boldsymbol{\Lambda}}_{\beta})\}\boldsymbol{1}_{p},
\end{aligned}
}\end{equation}

where

\[
W_{i}(\widetilde{\boldsymbol{\mu}}_{\beta},\widetilde{\boldsymbol{\Lambda}
}_{\beta})=\exp\left\{  -\frac{\beta}{2}(\boldsymbol{X}_{i}
-\widetilde{\boldsymbol{\mu}}_{\beta})^{\top}\widetilde{\boldsymbol{\Lambda}
}_{\beta}^{-1/2}\boldsymbol{R}_{0}^{-1}\widetilde{\boldsymbol{\Lambda}}
_{\beta}^{-1/2}(\boldsymbol{X}_{i}-\widetilde{\boldsymbol{\mu}}_{\beta
})\right\}
\]

and

\begin{equation}\protect\phantomsection\label{eq-RR}{
\begin{aligned}
\boldsymbol{R}_{\boldsymbol{X},\beta}(\widetilde{\boldsymbol{\mu}}_{\beta
},\widetilde{\boldsymbol{\Lambda}}_{\beta})  &
=\widetilde{\boldsymbol{\Lambda}}_{\beta}^{-1/2}\boldsymbol{S}_{\boldsymbol{X}
,\beta}(\widetilde{\boldsymbol{\mu}}_{\beta},\widetilde{\boldsymbol{\Lambda}
}_{\beta})\widetilde{\boldsymbol{\Lambda}}_{\beta}^{-1/2}, \\
\boldsymbol{S}_{\boldsymbol{X},\beta}(\widetilde{\boldsymbol{\mu}}_{\beta
},\widetilde{\boldsymbol{\Lambda}}_{\beta})  &  =\frac{\frac{1}{n}
\sum\limits_{i=1}^{n}W_{i}(\widetilde{\boldsymbol{\mu}}_{\beta}
,\widetilde{\boldsymbol{\Lambda}}_{\beta})(\boldsymbol{X}_{i}
-\widetilde{\boldsymbol{\mu}}_{\beta})(\boldsymbol{X}_{i}
-\widetilde{\boldsymbol{\mu}}_{\beta})^{\top}}{\widetilde{\kappa}_{0}(p,\beta
)}.
\end{aligned}
}\end{equation}

\end{theorem}

\begin{refremark}
For \(\beta=0\),
\(W_{i,\beta=0}(\widetilde{\boldsymbol{\mu}}_{\beta},\widetilde{\boldsymbol{\Lambda}}_{\beta})=1\),
\(i=1,\ldots,n\), so that
\(\widetilde{\boldsymbol{\mu}}_{\beta=0}=\overline{\boldsymbol{X}}\),
but the restricted MLEs of the variances have no explicit expression for
a general \(\boldsymbol{R}_{0}\), see Remark~\ref{rem-Remark1}. For
\(\beta>0\) all the restricted minimum DPD estimators require iterative
computations.

\label{rem-mle}

\end{refremark}

\begin{refremark}
The expression of
\(\boldsymbol{R}_{\boldsymbol{X},\beta }(\widetilde{\boldsymbol{\mu}}_{\beta},\widetilde{\boldsymbol{\Lambda}} _{\beta})\)
for \(\beta=0\)~does not match (in most cases) the sample Pearson
correlation matrix since the variance estimators fail to be the ordinary
sample variances. For example for \(p=2\),

\[
\begin{aligned}
\boldsymbol{R}_{0}=
\begin{pmatrix}
1 & \rho_{0}\\
\rho_{0} & 1
\end{pmatrix}
,\qquad\boldsymbol{R}_{\boldsymbol{X}}(\widetilde{\boldsymbol{\Lambda}})=
\begin{pmatrix}
\frac{S_{1}^{2}}{\widetilde{\sigma}_{1}^{2}} & \frac{S_{12}}{\widetilde{\sigma
}_{1}\widetilde{\sigma}_{2}}\\
\frac{S_{12}}{\widetilde{\sigma}_{1}\widetilde{\sigma}_{2}} & \frac{S_{2}^{2}
}{\widetilde{\sigma}_{2}^{2}}
\end{pmatrix}
,
\end{aligned}
\]

and the estimating equations for the variances are

\[
\begin{aligned}
\begin{pmatrix}
1\\
1
\end{pmatrix}
=\frac{1}{1-\rho_{0}^{2}}\mathrm{diag}\left\{
\begin{pmatrix}
1 & -\rho_{0}\\
-\rho_{0} & 1
\end{pmatrix}
\begin{pmatrix}
\frac{S_{1}^{2}}{\widetilde{\sigma}_{1}^{2}} & \frac{S_{12}}{\widetilde{\sigma
}_{1}\widetilde{\sigma}_{2}}\\
\frac{S_{12}}{\widetilde{\sigma}_{1}\widetilde{\sigma}_{2}} & \frac{S_{2}^{2}
}{\widetilde{\sigma}_{2}^{2}}
\end{pmatrix}
\right\}
\begin{pmatrix}
1\\
1
\end{pmatrix}
,
\end{aligned}
\]

i.e.

\[
\frac{1}{1-\rho_{0}^{2}}\left(  \frac{S_{j}^{2}}{\widetilde{\sigma}_{j}^{2}
}-\rho_{0}\frac{S_{12}}{\widetilde{\sigma}_{1}\widetilde{\sigma}_{2}}\right)
=1,\quad j=1,2,
\]

or equivalently

\[
\frac{S_{j}^{2}}{\widetilde{\sigma}_{j}^{2}}=1-\rho_{0}^{2}+\rho_{0}
\frac{S_{12}}{\widetilde{\sigma}_{1}\widetilde{\sigma}_{2}},\quad j=1,2,
\]

from which the ratio is deducted to be

\[
\frac{\widetilde{\sigma}_{1}}{\widetilde{\sigma}_{2}}=\frac{S_{1}}{S_{2}},
\]

but

\[
\begin{pmatrix}
\widetilde{\sigma}_{1}^{2}\\
\widetilde{\sigma}_{2}^{2}
\end{pmatrix}
=\mathrm{diag}\{\boldsymbol{R}_{0}^{-1}\boldsymbol{R}_{\boldsymbol{X}}\}
\begin{pmatrix}
S_{1}^{2}\\
S_{2}^{2}
\end{pmatrix}
,
\]

i.e.

\[
\begin{aligned}
\begin{pmatrix}
\widetilde{\sigma}_{1}^{2}\\
\widetilde{\sigma}_{2}^{2}
\end{pmatrix}
=\frac{1-\rho_{0}R_{12}}{1-\rho_{0}^{2}}
\begin{pmatrix}
1 & 0\\
0 & 1
\end{pmatrix}
\begin{pmatrix}
S_{1}^{2}\\
S_{2}^{2}
\end{pmatrix}
,
\end{aligned}
\]

or

\[
\widetilde{\sigma}_{j}^{2}=S_{j}^{2}\frac{1-\rho_{0}R_{12}}{1-\rho_{0}^{2}
},\quad j=1,2.
\]

In the foregoing subsections it is shown that both match under
uncorrelatedness, \(\rho_{0}=0\).

\label{rem-Remark1}

\end{refremark}

\begin{theorem}[]\protect\hypertarget{thm-RestrMDPDsEqui}{}\label{thm-RestrMDPDsEqui}

Under fixed equicorrelation,

\[
\boldsymbol{R}(\rho_{0})=(1-\rho_{0})\boldsymbol{I}_{p}+\rho_{0}
\boldsymbol{1}_{p}\boldsymbol{1}_{p}^{\top},
\]

the restricted minimum DPD estimators of
\((\boldsymbol{\mu}^{\top},\boldsymbol{1} _{p}^{\top}\boldsymbol{\Lambda})^{\top}\),
are obtained as solution in (\ref{eq-eqMuBeta}) and

\begin{equation}\protect\phantomsection\label{eq-estimEqui}{
\widetilde{R}_{jj,\beta}-(1-\rho_{0})=\frac{\rho_{0}}{1+(p-1)\rho_{0}
}\widetilde{R}_{\cdot j,\beta},
}\end{equation}

where
\(\widetilde{R}_{jj,\beta}=S_{j,\beta}^{2}/\widetilde{\sigma}_{j,\beta }^{2}\)
are the diagonal elements of (\ref{eq-RR}).

\end{theorem}

\begin{theorem}[]\protect\hypertarget{thm-ThEstimRho}{}\label{thm-ThEstimRho}

Under non-fixed equicorrelation,
\begin{equation}\protect\phantomsection\label{eq-R-Equi2}{
\boldsymbol{R}(\rho_{12})=(1-\rho_{12})\boldsymbol{I}_{p}+\rho_{12}
\boldsymbol{1}_{p}\boldsymbol{1}_{p}^{\top},
}\end{equation} being \(\rho_{12}\in(-\frac{1}{p-1},1)\) an unknown
parameter, the restricted minimum DPD estimators
\((\widetilde{\boldsymbol{\mu}}_{\beta}^{\top},\boldsymbol{1}_{p}^{\top}\widetilde{\boldsymbol{\Lambda}}_{\beta},\widetilde{\rho}_{12,\beta})^{\top}\)
are obtained as solution of (\ref{eq-eqMuBeta}),
\begin{equation}\protect\phantomsection\label{eq-rhoEstVar}{
\mathrm{diag}\{\boldsymbol{R}^{-1}(\widetilde{\rho}_{12,\beta})\widetilde{\boldsymbol{R}}_{\boldsymbol{X},\beta}\}\boldsymbol{1}_{p}=\boldsymbol{1}_{p},
}\end{equation} and
\begin{equation}\protect\phantomsection\label{eq-rhoEst}{
\boldsymbol{1}_{p}^{\top}\widetilde{\boldsymbol{M}}_{\beta}\boldsymbol{1}_{p}=\mathrm{trace}(\widetilde{\boldsymbol{M}}_{\beta}),
\qquad
\widetilde{\boldsymbol{M}}_{\beta}=\boldsymbol{R}^{-1}(\widetilde{\rho}_{12,\beta})\left(\widetilde{\boldsymbol{R}}_{\boldsymbol{X},\beta}-\boldsymbol{R}(\widetilde{\rho}_{12,\beta})\right)\boldsymbol{R}^{-1}(\widetilde{\rho}_{12,\beta}),
}\end{equation} where
\(\widetilde{\boldsymbol{R}}_{\boldsymbol{X},\beta}\) is (\ref{eq-RR})
with \(\boldsymbol{R}_{0}\) replaced by
\(\boldsymbol{R}(\widetilde{\rho}_{12,\beta})\) in the weights. The left
hand side of (\ref{eq-rhoEst}) is equal to \[
\{1+(p-1)\widetilde{\rho}_{12,\beta}\}^{-2}\{\boldsymbol{1}_{p}^{\top}\widetilde{\boldsymbol{R}}_{\boldsymbol{X},\beta}\boldsymbol{1}_{p}-p(1+(p-1)\widetilde{\rho}_{12,\beta})\}.
\]

\end{theorem}

\begin{refremark}
The estimating equations of Theorem~\ref{thm-ThEstimRho} have no
explicit solution, also for \(\beta=0\). In particular, the estimators
of the variances are not the weighted sample variances and
\(\widetilde{\rho}_{12,\beta}\) is not the mean of the off-diagonal
elements of the weighted correlation matrix, since, by
(\ref{eq-rhoEstVar}), this would require the \(p\) row means of these
elements to be equal, an event of probability zero for \(p\geq3\). They
can be computed by numerical minimization of the density power
divergence objective function under the restriction, with weights
\(\widetilde{W}_{i,\beta}\) given by (\ref{eq-w}) at the restricted
estimators, and the constant
\begin{equation}\protect\phantomsection\label{eq-kappa0}{
\widetilde{\kappa}_{0}(p,\beta)=\frac{1}{n}\sum_{i=1}^{n}\widetilde{W}
_{i,\beta}-\beta(\beta+1)^{-(\frac{p}{2}+1)}
}\end{equation} appears in all the test-statistics of this section. For
\(\beta=0\), \(\widetilde{W}_{i,\beta=0}=1\), \(i=1,\ldots,n\), and
\(\widetilde{\kappa}_{0}(p,\beta=0)=1\).

\label{rem-recursive}

\end{refremark}

\subsection{Testing specified values for correlation
matrix}\label{sec-test1}

For testing \(\boldsymbol{R}=\boldsymbol{R}_{0}\)~or more formally

\begin{equation}\protect\phantomsection\label{eq-nullVSalt}{
H_{0}\text{: }\boldsymbol{\eta}_{2,2}=\mathrm{vecl}(\boldsymbol{R}_{0})\quad
vs.\quad H_{1}\text{: }\boldsymbol{\eta}_{2,2}\neq\mathrm{vecl}(\boldsymbol{R}
_{0}),
}\end{equation}

let us consider a transformation of the original sample,
\(\widetilde{\boldsymbol{X}}_{i,\beta}\), \(i=1,...,n\), where

\[
\widetilde{\boldsymbol{X}}_{i,\beta}=\widetilde{\boldsymbol{\Lambda}}_{\beta
}^{-1/2}(\boldsymbol{X}_{i}-\widetilde{\boldsymbol{\mu}}_{\beta}),
\]

being \(\widetilde{\boldsymbol{\mu}}_{\beta}\) and
\(\widetilde{\boldsymbol{\Lambda}}_{\beta}\)~the~minimum DPD estimators
of \(\boldsymbol{\mu}\) and \(\boldsymbol{\Lambda}\) under
\(\boldsymbol{R} =\boldsymbol{R}_{0}\) (see Section \ref{sec-estim}).

\begin{theorem}[]\protect\hypertarget{thm-ThTestGeneral}{}\label{thm-ThTestGeneral}

The Rao's \(\beta\)-score test-statistic for (\ref{eq-nullVSalt}) is

\begin{equation}\protect\phantomsection\label{eq-RaoR0}{
\widetilde{R}_{\beta,n}=n\widetilde{\kappa}_{0}^{2}(p,\beta)\kappa_{1}
^{-1}(p,\beta)\mathrm{trace}\left(  \left(  \boldsymbol{R}_{0}^{-1}
\widetilde{\boldsymbol{R}}_{\boldsymbol{X},\beta}-\boldsymbol{I}_{p}\right)
^{2}\right)  ,
}\end{equation}

with
\(\widetilde{\boldsymbol{R}}_{\boldsymbol{X},\beta}=\boldsymbol{R} _{\boldsymbol{X},\beta}(\widetilde{\boldsymbol{\mu}}_{\beta} ,\widetilde{\boldsymbol{\Lambda}}_{\beta})\),
given previously in (\ref{eq-RR}), i.e.

\[
\widetilde{\boldsymbol{R}}_{\boldsymbol{X},\beta}=\frac{\frac{1}{n}
\sum\limits_{i=1}^{n}\widetilde{W}_{i,\beta}\widetilde{\boldsymbol{X}
}_{i,\beta}\widetilde{\boldsymbol{X}}_{i,\beta}^{\top}}{\widetilde{\kappa}
_{0}(p,\beta)},
\]

\(\widetilde{W}_{i,\beta}=W_{i,\beta}(\widetilde{\boldsymbol{\theta}})\),
given previously in (\ref{eq-w}), i.e.

\[
\widetilde{W}_{i,\beta}=\exp\left\{  -\frac{\beta}{2}\widetilde{\boldsymbol{X}
}_{i,\beta}^{\top}\boldsymbol{R}_{0}^{-1}\widetilde{\boldsymbol{X}}_{i,\beta
}\right\}  ,
\]

\(\widetilde{\kappa}_{0}(p,\beta)\) was given in (\ref{eq-kappa0}) and
\(\kappa _{1}(p,\beta)\) in (\ref{eq-Kappa1}). The asymptotic
distribution of (\ref{eq-RaoR0}) is \(\chi^{2}\) with
\(\frac{p(p-1)}{2}\) degrees of freedom.

\end{theorem}

\begin{refremark}
Write
\(\widetilde{\boldsymbol{D}}_{\beta}=\boldsymbol{R}_{0}^{-1}\widetilde{\boldsymbol{R}}_{\boldsymbol{X},\beta}-\boldsymbol{I}_{p}\).
By the estimating equations of the variances its diagonal elements are
zero, and therefore
\begin{equation}\protect\phantomsection\label{eq-trace-vecl}{
\mathrm{trace}\left(\widetilde{\boldsymbol{D}}_{\beta}^{2}\right)=2\sum_{i<j}\widetilde{D}_{ij,\beta}\widetilde{D}_{ji,\beta}.
}\end{equation} The matrix \(\widetilde{\boldsymbol{D}}_{\beta}\) is not
symmetric in general, because \(\boldsymbol{R}_{0}^{-1}\) and
\(\widetilde{\boldsymbol{R}}_{\boldsymbol{X},\beta}\) do not commute, so
the product \(\widetilde{D}_{ij,\beta}\widetilde{D}_{ji,\beta}\) cannot
be replaced by \(\widetilde{D}_{ij,\beta}^{2}\), except when
\(\boldsymbol{R}_{0}=\boldsymbol{I}_{p}\).

\label{rem-trace}

\end{refremark}

\begin{refremark}
For \(\beta=0\), i.e.~for all the classical Rao's score tests, it hold

\[
\widetilde{\kappa}_{0}^{2}(p,\beta=0)\kappa_{1}^{-1}(p,\beta=0)=\frac{1}{2}.
\]

\label{rem-beta0}

\end{refremark}

\begin{refremark}
Particularization of the Rao's \(\beta\)-score test-statistic for
testing specified values of correlation in the bidimensional case
(\(p=2\)):

\begin{equation}\protect\phantomsection\label{eq-RaoR0d}{
\widetilde{R}_{\beta,n}=2n\frac{\widetilde{\kappa}_{0}^{2}(p=2,\beta)}
{\kappa_{1}(p=2,\beta)}\left(  \frac{R_{12,\beta}-\rho_{0}}{1-\rho
_{0}R_{12,\beta}}\right)  ^{2},
}\end{equation}

where the weights to calculate the sample correlation is given by

\[
\widetilde{W}_{i,\beta}=\exp\left\{  -\frac{\beta}{2(1-\rho_{0})}\left(
\widetilde{\boldsymbol{X}}_{i}^{\top}\widetilde{\boldsymbol{X}}_{i}-\frac
{\rho_{0}}{1+(p-1)\rho_{0}}(\widetilde{\boldsymbol{X}}_{i}^{\top}\boldsymbol{1}
_{2})^{2}\right)  \right\}  .
\]

The likelihood ratio test statistic can be found in Chapter 4 of
Anderson (2003), but to our knowledge, the explicit expression of the
classical score test statistic (\(\beta=0\)),

\[
\widetilde{R}_{n}=n\left(  \frac{R_{12}-\rho_{0}}{1-\rho_{0}R_{12}}\right)
^{2},
\]

had not been published yet. It suits results such as (5.4.34)-(5.4.40)
of Lehman (1999, page 316) for any distribution, but it is not an
equivalent expression. In fact, the most well-known test-statistic for
the correlation coefficient of the bivariate normal distribution is the
one for the Fisher's transform of the sample correlation coefficient
(see page (5.4.41)-(5.4.42) in Lehmann (1999)).

\label{rem-p2}

\end{refremark}

\subsection{Testing fixed equicorrelation}\label{sec-test2}

The equicorrelation structure establishes homogeneity for all
off-diagonal elements of the correlation matrix, i.e.

\[
H_{0}\text{: }\rho_{ij}=\rho_{0}\text{, }\forall i\neq j\quad vs.\quad
H_{1}\text{: }\exists i\neq j\text{ s.t. }\rho_{ij}\neq\rho_{0},
\]

with \(\rho_{0}\)~being the so called intraclass correlation satisfying

\[
-\frac{1}{p-1}<\rho_{0}<1,
\]

to guarantee positive definiteness. In matrix form, the equicorrelation
structure can be expressed as

\[
H_{0}\text{: }\boldsymbol{\eta}_{2,2}=\mathrm{vecl}(\boldsymbol{R}(\rho
_{0}))\quad vs.\quad H_{1}\text{: }\boldsymbol{\eta}_{2,2}\neq\mathrm{vecl}
(\boldsymbol{R}(\rho_{0})),
\]

where

\[
\boldsymbol{R}(\rho_{0})=(1-\rho_{0})\boldsymbol{I}_{p}+\rho_{0}
\boldsymbol{1}_{p}\boldsymbol{1}_{p}^{\top}.
\]

\begin{corollary}[]\protect\hypertarget{cor-ThTestEquiF}{}\label{cor-ThTestEquiF}

The Rao's \(\beta\)-score test-statistic for testing equicorrelation, as
particular case of (\ref{eq-RaoR0}), is given by
\begin{equation}\protect\phantomsection\label{eq-RaoR0f}{
\widetilde{R}_{\beta,n}=\frac{2n\widetilde{\kappa}_{0}^{2}(p,\beta)\kappa
_{1}^{-1}(p,\beta)}{(1-\rho_{0})^{2}}\sum_{i<j}\left(\widetilde{R}_{ij,\beta}-\widetilde{R}_{ii,\beta}+1-\rho_{0}\right)\left(\widetilde{R}_{ij,\beta}-\widetilde{R}_{jj,\beta}+1-\rho_{0}\right),
}\end{equation} where
\(\widetilde{R}_{ij,\beta}=S_{ij,\beta}/(\widetilde{\sigma}_{i,\beta}\widetilde{\sigma}_{j,\beta})\)
and
\(\widetilde{R}_{jj,\beta}=S_{j,\beta}^{2}/\widetilde{\sigma}_{j,\beta}^{2}\)
are the elements of (\ref{eq-RR}), with \[
\begin{aligned}
\widetilde{W}_{i,\beta}  &  =\exp\left\{  -\tfrac{\beta}{2(1-\rho_{0})}\left[
\sum_{j=1}^{p}\widetilde{X}_{ij,\beta}^{2}-\tfrac{\rho_{0}}{1+(p-1)\rho_{0}
}\left(  \sum_{j=1}^{p}\widetilde{X}_{ij,\beta}\right)  ^{2}\right]  \right\}
, \\
\widetilde{X}_{ij,\beta}  &  =\frac{X_{ij}-\widetilde{\mu}_{j,\beta}
}{\widetilde{\sigma}_{j,\beta}},\quad i=1,\ldots,n,\quad j=1,\ldots,p.
\end{aligned}
\] Its asymptotic distribution under the null hypothesis is \(\chi^{2}\)
with \(\frac{p(p-1)}{2}\) degrees of freedom.

\end{corollary}

\subsection{Testing complete uncorrelatedness or
independence}\label{sec-test3}

\begin{corollary}[]\protect\hypertarget{cor-ThTestInd}{}\label{cor-ThTestInd}

The Rao's \(\beta\)-score test-statistic for testing uncorrelatedness or
independence (\(\boldsymbol{R}(\rho_{0}=0)\)), as particular case of
(\ref{eq-RaoR0}), is given by

\begin{equation}\protect\phantomsection\label{eq-RaoR0bb}{
\widetilde{R}_{\beta,n}=2n\widetilde{\kappa}_{0}^{2}(p,\beta)\kappa_{1}
^{-1}(p,\beta)\sum_{i<j}^{p}R_{ij,\beta}^{2},
}\end{equation}

where \(R_{ij,\beta}\)~are extradiagonal elements of (\ref{eq-RR}) when
\(\widetilde{\sigma}_{j,\beta}^{2}=S_{j,\beta}^{2}\), \(j=1,...,p\). Its
asymptotic distribution is \(\chi^{2}\) with \(\frac{p(p-1)}{2}\)
degrees of freedom.

\end{corollary}

For the specific case of \(p=2\),
\(\widetilde{R}_{\beta,n}=2n\widetilde{\kappa }_{0}^{2}(p=2,\beta)\kappa_{1}^{-1}(p=2,\beta)R_{12,\beta}^{2}\),
with \(2\kappa_{1}^{-1}(p,\beta)=(2\beta+1)^{3}\) (see Example 4 in Basu
et al.~(2021)).

\begin{refremark}
Notice that (\ref{eq-RaoR0bb}) has an explicit expression for
\(\beta=0\), which is based on the classical Pearson sample correlations
(without constraints),
\(\widetilde{R}_{n}=\widetilde{R}_{\beta=0,n}=n\sum_{i<j}^{p}R_{ij,\beta}^{2}\),
while for \(\beta>0\), the correlations are calculated as

\[
\begin{aligned}
\boldsymbol{R}_{\boldsymbol{X},\beta}  &  =\mathrm{diag}^{-1}\{S_{j,\beta
}\}_{j=1}^{p}\boldsymbol{S}_{\boldsymbol{X},\beta}\mathrm{diag}^{-1}
\{S_{j,\beta}\}_{j=1}^{p}, \\
\boldsymbol{S}_{\boldsymbol{X},\beta}  &  =\frac{\frac{1}{n}\sum
\limits_{i=1}^{n}\widetilde{W}_{i,\beta}(\boldsymbol{X}_{i}
-\widetilde{\boldsymbol{\mu}}_{\beta})(\boldsymbol{X}_{i}
-\widetilde{\boldsymbol{\mu}}_{\beta})^{\top}}{\widetilde{\kappa}_{0}(p,\beta
)}, \\
S_{i,\beta}^{2}  &  =\mathrm{diag}(\boldsymbol{S}_{\boldsymbol{X},\beta}), \\
\widetilde{\boldsymbol{\mu}}_{\beta}  &  =\frac{\sum\limits_{i=1}^{n}\widetilde{W}_{i,\beta}\boldsymbol{X}_{i}}{\sum\limits_{i=1}^{n}\widetilde{W}_{i,\beta}}, \\
\widetilde{W}_{i,\beta}  &  =\exp\left\{  -\tfrac{\beta}{2}(\boldsymbol{X}
_{i}-\widetilde{\boldsymbol{\mu}}_{\beta})^{\top}\mathrm{diag}^{-1}\{S_{j,\beta
}^{2}\}_{j=1}^{p}(\boldsymbol{X}_{i}-\widetilde{\boldsymbol{\mu}}_{\beta
})\right\}  ,
\end{aligned}
\]

\(\widetilde{\kappa}_{0}(p,\beta)\) was given in (\ref{eq-kappa0}).

\label{rem-indep}

\end{refremark}

\begin{refremark}
In the traditional multivariate analysis, when \(p\) is small relative
to \(n\), Bartlett (1954) established the likelihood ratio test for the
complete independence as

\[
-\left(  n-1-\frac{2p+5}{6}\right)  \log\left\vert \boldsymbol{R}
_{\boldsymbol{X},\beta=0}\right\vert \underset{n\rightarrow\infty
}{\overset{\mathcal{L}}{\longrightarrow}}\chi_{\frac{p(p-1)}{2}}^{2},
\]

which is not longer valid for \(p>n\), since the lack of positive
definiteness makes the determinant to be null. Such a problem of the
likelihood ratio test does not exist for the Rao's score test.

\label{rem-bartlett}

\end{refremark}

\begin{refremark}
In Kallenberg et al.~(1997) it is mentioned
\(\widetilde{R}_{n}=nR_{12}^{2}\) ((\ref{eq-RaoR0bb}), with \(p=2\) and
\(\beta=0\)) to be the classical Rao test statistic for testing
independence for bivariate normal random variables and its asymptotic
standard normality is very well-known from different sources such as the
example given in page 293 of Lehmann (1999).

\label{rem-kallenberg}

\end{refremark}

\subsection{Testing non-fixed equicorrelation}\label{sec-test4}

The non-fixed equicorrelation structure establishes homogeneity for all
off-diagonal elements of the correlation matrix, i.e.

\[
H_{0}\text{: }\rho_{ij}=\rho_{12}\text{, }\forall i\neq j\quad vs.\quad
H_{1}\text{: }\exists i\neq j\text{ s.t. }\rho_{ij}\neq\rho_{12},
\]

for unknown value of \(\rho_{12}\). In addition, it is assumed

\[
\frac{-1}{p-1}<\rho_{12}<1,
\]

for imposing positive definiteness. In matrix form, the equicorrelation
can be expressed as

\[
\boldsymbol{R}(\rho_{12})=(1-\rho_{12})\boldsymbol{I}_{p}+\rho_{12}
\boldsymbol{1}_{p}\boldsymbol{1}_{p}^{\top}.
\]

\begin{corollary}[]\protect\hypertarget{cor-ThTestEquiNF}{}\label{cor-ThTestEquiNF}

The Rao's \(\beta\)-score test-statistic for testing equicorrelation
with an unknown common correlation is
\begin{equation}\protect\phantomsection\label{eq-RaoR0g}{
\widetilde{R}_{\beta,n}=n\widetilde{\kappa}_{0}^{2}(p,\beta)\kappa_{1}^{-1}(p,\beta)\mathrm{trace}\left(\left(\boldsymbol{R}^{-1}(\widetilde{\rho}_{12,\beta})\widetilde{\boldsymbol{R}}_{\boldsymbol{X},\beta}-\boldsymbol{I}_{p}\right)^{2}\right),
}\end{equation} with the restricted estimators of
Theorem~\ref{thm-ThEstimRho}. Its asymptotic distribution under the null
hypothesis is \(\chi^{2}\) with
\(\frac{p(p-1)}{2}-1=\frac{(p+1)(p-2)}{2}\) degrees of freedom.

\end{corollary}

\begin{refremark}
For \(\beta=0\), let
\(\overline{r}_{k}=\frac{1}{p-1}\sum_{i\neq k}R_{ik}\) and
\(\overline{r}=\frac{2}{p(p-1)}\sum_{i<j}R_{ij}\) be computed from the
sample Pearson correlation matrix, \(\widehat{\lambda}=1-\overline{r}\)
and
\(\widehat{\mu}=(p-1)^{2}(1-\widehat{\lambda}^{2})/\{p-(p-2)\widehat{\lambda}^{2}\}\).
The test-statistic of Lawley (1963), \[
L_{n}=\frac{n}{\widehat{\lambda}^{2}}\left(\sum_{i<j}(R_{ij}-\overline{r})^{2}-\widehat{\mu}\sum_{k=1}^{p}(\overline{r}_{k}-\overline{r})^{2}\right),
\] has the same asymptotic distribution under the null hypothesis, and
Lawley (1963) shows that it is asymptotically equal to \(-2\) times the
logarithm of the likelihood ratio test-statistic. Therefore
(\ref{eq-RaoR0g}) with \(\beta=0\) and \(L_{n}\) are asymptotically
equivalent, and the second term of \(L_{n}\) accounts for the fact that
the sample variances are not the restricted estimators of the variances.
The test-statistics with \(\beta>0\) have no classical counterpart.

\label{rem-lawley}

\end{refremark}

\section{Asymptotic properties and robustness}\label{sec-Asym}

Throughout this section
\(\boldsymbol{\theta}_{0}=(\boldsymbol{\mu}_{0}^{\top},\mathrm{vech}^{\top}(\boldsymbol{\Sigma}_{0}))^{\top}\),
with
\(\boldsymbol{\Sigma}_{0}=\boldsymbol{\Lambda}_{0}^{1/2}\boldsymbol{R}_{0}\boldsymbol{\Lambda}_{0}^{1/2}\),
denotes the true value of the parameter under the null hypothesis, where
\(\boldsymbol{R}_{0}=\boldsymbol{R}(\rho_{12,0})\) for the test of
Section \ref{sec-test4}, \(\boldsymbol{\eta}_{0}\) is the value of
\(\boldsymbol{\eta}\) which corresponds to \(\boldsymbol{\theta}_{0}\),
and \(r\) denotes the number of degrees of freedom of the test,
\(r=\frac{p(p-1)}{2}\) in Sections \ref{sec-test1}, \ref{sec-test2} and
\ref{sec-test3} and \(r=\frac{(p+1)(p-2)}{2}\), with \(p\geq3\), in
Section \ref{sec-test4}. For \(\boldsymbol{x}\in\mathbb{R}^{p}\) let \[
\boldsymbol{z}(\boldsymbol{x})=\boldsymbol{\Lambda}_{0}^{-1/2}(\boldsymbol{x}-\boldsymbol{\mu}_{0}),\qquad
w_{\beta}(\boldsymbol{x})=\exp\left\{-\frac{\beta}{2}\boldsymbol{z}^{\top}(\boldsymbol{x})\boldsymbol{R}_{0}^{-1}\boldsymbol{z}(\boldsymbol{x})\right\},
\] so that \(w_{\beta}\) is the weight (\ref{eq-w}) as a function of the
point at which it is evaluated, and
\(W_{i,\beta}(\boldsymbol{\theta}_{0})=w_{\beta}(\boldsymbol{X}_{i})\).
Let also \begin{equation}\protect\phantomsection\label{eq-Psi}{
\boldsymbol{\Psi}_{\beta}(\boldsymbol{x})=w_{\beta}(\boldsymbol{x})\boldsymbol{z}(\boldsymbol{x})\boldsymbol{z}^{\top}(\boldsymbol{x})
-\left(w_{\beta}(\boldsymbol{x})-\beta(\beta+1)^{-(\frac{p}{2}+1)}\right)\boldsymbol{R}_{0},
}\end{equation} which is the contribution of the observation
\(\boldsymbol{x}\) to
\(\widetilde{\kappa}_{0}(p,\beta)(\widetilde{\boldsymbol{R}}_{\boldsymbol{X},\beta}-\boldsymbol{R}_{0})\)
when the restricted estimators are replaced by the true values, and
verifies
\(\mathrm{E}_{\boldsymbol{\theta}_{0}}[\boldsymbol{\Psi}_{\beta}(\boldsymbol{X})]=\boldsymbol{0}_{p\times p}\).

The test-statistics of Section \ref{sec-Sec3} are, by
Theorem~\ref{thm-ThTestGeneral},
\(\mathrm{trace}((\boldsymbol{R}_{0}^{-1}\widetilde{\boldsymbol{R}}_{\boldsymbol{X},\beta}-\boldsymbol{I}_{p})^{2})\)
up to a constant, and the diagonal elements of
\(\boldsymbol{R}_{0}^{-1}\widetilde{\boldsymbol{R}}_{\boldsymbol{X},\beta}-\boldsymbol{I}_{p}\)
are zero because of the estimating equations of the nuisance parameters,
see Remark~\ref{rem-trace}. The results below have the same structure,
with
\(\widetilde{\boldsymbol{R}}_{\boldsymbol{X},\beta}-\boldsymbol{R}_{0}\)
replaced by a symmetric matrix \(\boldsymbol{M}\) which has to be
adjusted in the same way. Accordingly, for a symmetric matrix
\(\boldsymbol{M}\) of order \(p\) we denote by
\begin{equation}\protect\phantomsection\label{eq-adjust}{
\boldsymbol{M}^{\ast}=\boldsymbol{M}-\boldsymbol{D}(\boldsymbol{M})\boldsymbol{R}_{0}-\boldsymbol{R}_{0}\boldsymbol{D}(\boldsymbol{M}),
\qquad \boldsymbol{D}(\boldsymbol{M})=\mathrm{diag}\{\boldsymbol{d}(\boldsymbol{M})\},
}\end{equation} the unique matrix of this form,
\(\boldsymbol{D}(\boldsymbol{M})\) being diagonal, such that
\begin{equation}\protect\phantomsection\label{eq-adjustvar}{
\mathrm{diag}\{\boldsymbol{R}_{0}^{-1}\boldsymbol{M}^{\ast}\}\boldsymbol{1}_{p}=\boldsymbol{0}_{p},
}\end{equation} for the tests of Sections \ref{sec-test1},
\ref{sec-test2} and \ref{sec-test3}. For the test of Section
\ref{sec-test4} a multiple of
\(\boldsymbol{1}_{p}\boldsymbol{1}_{p}^{\top}-\boldsymbol{I}_{p}\) is
subtracted as well, and \(\boldsymbol{M}^{\ast}\) is the unique matrix
\(\boldsymbol{M}-\boldsymbol{D}(\boldsymbol{M})\boldsymbol{R}_{0}-\boldsymbol{R}_{0}\boldsymbol{D}(\boldsymbol{M})-\lambda(\boldsymbol{M})(\boldsymbol{1}_{p}\boldsymbol{1}_{p}^{\top}-\boldsymbol{I}_{p})\)
which verifies (\ref{eq-adjustvar}) and
\begin{equation}\protect\phantomsection\label{eq-adjustrho}{
\boldsymbol{1}_{p}^{\top}\boldsymbol{R}_{0}^{-1}\boldsymbol{M}^{\ast}\boldsymbol{R}_{0}^{-1}\boldsymbol{1}_{p}=\mathrm{trace}(\boldsymbol{R}_{0}^{-1}\boldsymbol{M}^{\ast}\boldsymbol{R}_{0}^{-1}).
}\end{equation} Conditions (\ref{eq-adjustvar}) and (\ref{eq-adjustrho})
are (\ref{eq-estimEqui}) and (\ref{eq-rhoEst}), the estimating equations
of the variances and of the common correlation coefficient, written for
\(\boldsymbol{M}^{\ast}\), so that \(\boldsymbol{M}^{\ast}\) is the part
of \(\boldsymbol{M}\) which is not absorbed by the nuisance parameters.
In particular \(\boldsymbol{M}^{\ast}=\boldsymbol{0}_{p\times p}\) if
and only if \(\boldsymbol{M}\) is a direction in which the null
hypothesis remains true,
\(\mathrm{trace}((\boldsymbol{R}_{0}^{-1}\boldsymbol{M}^{\ast})^{2})\geq0\),
and \(\boldsymbol{M}^{\ast}=\boldsymbol{M}\) for the test of complete
uncorrelatedness whenever the diagonal elements of \(\boldsymbol{M}\)
are zero.

\subsection{Regularity conditions and choice of the
root}\label{sec-root}

The minimum DPD estimators are associated with the objective function
\begin{equation}\protect\phantomsection\label{eq-Hn}{
\begin{aligned}
H_{\beta,n}(\boldsymbol{\theta})&=\int_{\mathbb{R}^{p}}f_{\boldsymbol{\theta}}^{1+\beta}(\boldsymbol{x})d\boldsymbol{x}-\left(1+\frac{1}{\beta}\right)\frac{1}{n}\sum_{i=1}^{n}f_{\boldsymbol{\theta}}^{\beta}(\boldsymbol{X}_{i})\\
&=\frac{1}{(2\pi)^{\frac{\beta p}{2}}\left\vert \boldsymbol{\Sigma}\right\vert ^{\frac{\beta}{2}}}\left[(\beta+1)^{-\frac{p}{2}}-\left(1+\frac{1}{\beta}\right)\frac{1}{n}\sum_{i=1}^{n}W_{i,\beta}(\boldsymbol{\theta})\right],
\end{aligned}
}\end{equation} for \(\beta>0\), and with minus the mean log-likelihood
for \(\beta=0\), and
\(\frac{\partial}{\partial\boldsymbol{\theta}}H_{\beta,n}(\boldsymbol{\theta})=-(\beta+1)\boldsymbol{U}_{\beta,n}(\boldsymbol{\theta})\).
Hence the estimating equations of Section \ref{sec-estim} characterize
the stationary points of (\ref{eq-Hn}) restricted to the null
hypothesis, which need not be unique for \(\beta>0\). Moreover, the
restricted minimum DPD estimator cannot be defined, for every sample
size, as the global minimizer of (\ref{eq-Hn}), as the following remark
shows.

\begin{refremark}
Let \(\beta>0\) and \(n<(\beta+1)^{\frac{p}{2}+1}/\beta\). With
\(\boldsymbol{\mu}=\boldsymbol{X}_{1}\) and
\(\boldsymbol{\Lambda}=t\boldsymbol{I}_{p}\),
\(W_{1,\beta}(\boldsymbol{\theta})=1\) and, with probability one,
\(W_{i,\beta}(\boldsymbol{\theta})\) tends to zero as \(t\) tends to
zero, \(i=2,\ldots,n\), so that the term in brackets of (\ref{eq-Hn})
tends to \((\beta+1)^{-\frac{p}{2}}-\frac{\beta+1}{\beta n}<0\), while
\(\left\vert \boldsymbol{\Sigma}\right\vert ^{-\frac{\beta}{2}}=t^{-\frac{\beta p}{2}}\left\vert \boldsymbol{R}_{0}\right\vert ^{-\frac{\beta}{2}}\)
tends to infinity. Therefore (\ref{eq-Hn}) is not bounded from below
under the null hypothesis. For example, for \(\beta=0.05\) and \(p=5\)
this happens for every \(n\leq23\).

\label{rem-unbounded}

\end{refremark}

For this reason the restricted minimum DPD estimator is defined as the
solution of the estimating equations with the smallest value of
(\ref{eq-Hn}) in a compact set of the parameter space, and the next
theorem establishes that this is a consistent sequence of roots, its
asymptotic representation and how the root is identified. In this
section
\(\boldsymbol{\eta}_{(1)}=(\boldsymbol{\eta}_{1}^{\top},\boldsymbol{\eta}_{2,1}^{\top})^{\top}\),
formed by \(\boldsymbol{\mu}\) and the variances
\(\boldsymbol{\Lambda}\boldsymbol{1}_{p}\), and by the common
correlation coefficient \(\rho_{12}\) as well for the test of Section
\ref{sec-test4}, denotes the vector of nuisance parameters, with true
value \(\boldsymbol{\eta}_{(1),0}\) and restricted minimum DPD estimator
\(\widetilde{\boldsymbol{\eta}}_{(1),\beta}\), and the subscripts \(1\)
and \(2\) in the blocks of
\(\boldsymbol{J}_{\beta}(\boldsymbol{\eta})\),
\(\boldsymbol{K}_{\beta}(\boldsymbol{\eta})\) and
\(\boldsymbol{U}_{\beta,n}(\boldsymbol{\eta})\) refer to
\(\boldsymbol{\eta}_{(1)}\) and to the correlation parameters whose
value is fixed by the null hypothesis, as in Lemma~\ref{lem-short}.

\begin{theorem}[]\protect\hypertarget{thm-root}{}\label{thm-root}

Let \(\boldsymbol{X}_{1},\ldots,\boldsymbol{X}_{n}\) be independent and
identically distributed random vectors with distribution
\(\mathcal{N}_{p}(\boldsymbol{\mu}_{0},\boldsymbol{\Sigma}_{0})\), and
let \(\beta\geq0\). Then the following statements hold.

\begin{enumerate}
\def\labelenumi{(\roman{enumi})}
\item
  With probability tending to one as \(n\rightarrow\infty\), the
  restricted estimating equations of Section \ref{sec-estim} have a root
  \(\widetilde{\boldsymbol{\eta}}_{(1),\beta}\) such that
  \(\widetilde{\boldsymbol{\eta}}_{(1),\beta}\underset{n\rightarrow\infty}{\overset{P}{\longrightarrow}}\boldsymbol{\eta}_{(1),0}\).
\item
  For any such sequence of roots,
  \begin{equation}\protect\phantomsection\label{eq-linear}{
  \sqrt{n}(\widetilde{\boldsymbol{\eta}}_{(1),\beta}-\boldsymbol{\eta}_{(1),0})=\boldsymbol{J}_{\beta,11}^{-1}(\boldsymbol{\eta}_{0})\sqrt{n}\boldsymbol{U}_{\beta,n,1}(\boldsymbol{\eta}_{0})+\boldsymbol{r}_{n},\qquad
  \boldsymbol{r}_{n}\underset{n\rightarrow\infty}{\overset{P}{\longrightarrow}}\boldsymbol{0}.
  }\end{equation}
\item
  There exists \(r_{0}>0\) such that, with probability tending to one as
  \(n\rightarrow\infty\), the function
  \(\boldsymbol{\eta}_{(1)}\mapsto H_{\beta,n}(\boldsymbol{\eta}_{(1)},\mathrm{vecl}(\boldsymbol{R}_{0}))\)
  is strictly convex on the ball
  \(\{\boldsymbol{\eta}_{(1)}:\left\Vert \boldsymbol{\eta}_{(1)}-\boldsymbol{\eta}_{(1),0}\right\Vert \leq r_{0}\}\),
  the root of (i) is its unique stationary point in this ball, and it is
  the minimizer of (\ref{eq-Hn}) on it.
\item
  For any compact set \(C\) of values of \(\boldsymbol{\eta}_{(1)}\)
  which contains \(\boldsymbol{\eta}_{(1),0}\) in its interior, with
  probability tending to one as \(n\rightarrow\infty\), the root of (i)
  is the solution of the restricted estimating equations in \(C\) with
  the smallest value of (\ref{eq-Hn}).
\end{enumerate}

The asymptotic distributions of Theorem~\ref{thm-ThTestGeneral} and of
Corollary~\ref{cor-ThTestEquiF}, Corollary~\ref{cor-ThTestInd} and
Corollary~\ref{cor-ThTestEquiNF} hold for the test-statistics evaluated
at the sequence of roots of (i).

\end{theorem}

\begin{refremark}
By (iii), the root is identified by any descent algorithm applied to
(\ref{eq-Hn}), under the null hypothesis, started at a consistent
estimator of \(\boldsymbol{\eta}_{(1),0}\), since such a starting value
belongs to the ball with probability tending to one. Under contamination
the starting value must also be robust, and the coordinatewise medians,
the squared normalized median absolute deviations and, for the test of
Section \ref{sec-test4}, the mean of \(\sin(\frac{\pi}{2}\tau_{ij})\),
with \(\tau_{ij}\) the sample Kendall correlation coefficients, are
consistent under the multivariate normal distribution. When the
algorithm, started at several values, reaches different solutions of the
estimating equations, the one with the smallest value of (\ref{eq-Hn})
is taken, which by (iv) is the consistent root as long as the solutions
remain in a compact set, that is, away from the degenerate solutions of
Remark~\ref{rem-unbounded}.

\label{rem-rootpractice}

\end{refremark}

\subsection{Asymptotic power under contiguous
alternatives}\label{sec-contig}

Let \(\boldsymbol{\Delta}\) be a fixed symmetric matrix of order \(p\)
with zero diagonal elements and consider the sequence of alternative
hypotheses \begin{equation}\protect\phantomsection\label{eq-H1n}{
H_{1,n}:\boldsymbol{R}=\boldsymbol{R}_{n}=\boldsymbol{R}_{0}+n^{-1/2}\boldsymbol{\Delta},
}\end{equation} with \(\boldsymbol{\mu}=\boldsymbol{\mu}_{0}\) and
\(\boldsymbol{\Lambda}=\boldsymbol{\Lambda}_{0}\), for \(n\) large
enough for \(\boldsymbol{R}_{n}\) to be positive definite. We denote by
\(\chi_{r}^{2}(\delta)\) the noncentral chi-square distribution with
\(r\) degrees of freedom and noncentrality parameter \(\delta\).

\begin{theorem}[]\protect\hypertarget{thm-contig}{}\label{thm-contig}

Under (\ref{eq-H1n}), for every \(\beta\geq0\), the Rao's
\(\beta\)-score test-statistics (\ref{eq-RaoR0}) and (\ref{eq-RaoR0g})
verify \begin{equation}\protect\phantomsection\label{eq-noncentral}{
\widetilde{R}_{\beta,n}\underset{n\rightarrow\infty}{\overset{\mathcal{L}}{\longrightarrow}}\chi_{r}^{2}\left(\delta_{\beta}(\boldsymbol{\Delta})\right),
}\end{equation} with
\begin{equation}\protect\phantomsection\label{eq-ARE}{
\delta_{\beta}(\boldsymbol{\Delta})=\frac{1}{2}\mathrm{ARE}(\beta,p)\,\mathrm{trace}\left(\left(\boldsymbol{R}_{0}^{-1}\boldsymbol{\Delta}^{\ast}\right)^{2}\right),
\qquad
\mathrm{ARE}(\beta,p)=\frac{(1+2\beta)^{\frac{p}{2}+2}}{(1+\beta)^{p+4}},
}\end{equation} where \(\boldsymbol{\Delta}^{\ast}\) is given by
(\ref{eq-adjust}). In particular, for testing complete uncorrelatedness,
\(\boldsymbol{\Delta}^{\ast}=\boldsymbol{\Delta}\) and
\(\delta_{\beta}(\boldsymbol{\Delta})=\mathrm{ARE}(\beta,p)\sum_{i<j}\Delta_{ij}^{2}\).

\end{theorem}

\begin{corollary}[]\protect\hypertarget{cor-are}{}\label{cor-are}

The asymptotic relative efficiency, in the sense of Pitman, of the Rao's
\(\beta\)-score test with respect to the classical Rao's score test
(\(\beta=0\)) is \(\mathrm{ARE}(\beta,p)\), for the four tests of
Section \ref{sec-Sec3}. It depends neither on \(\boldsymbol{R}_{0}\),
nor on \(\boldsymbol{\Delta}\), nor on the nuisance parameters, it is
equal to one only for \(\beta=0\), it is strictly decreasing in
\(\beta\) and, for \(\beta>0\), strictly decreasing in \(p\).

\end{corollary}

\begin{figure}[!htb]

\centering{

\pandocbounded{\includegraphics[keepaspectratio]{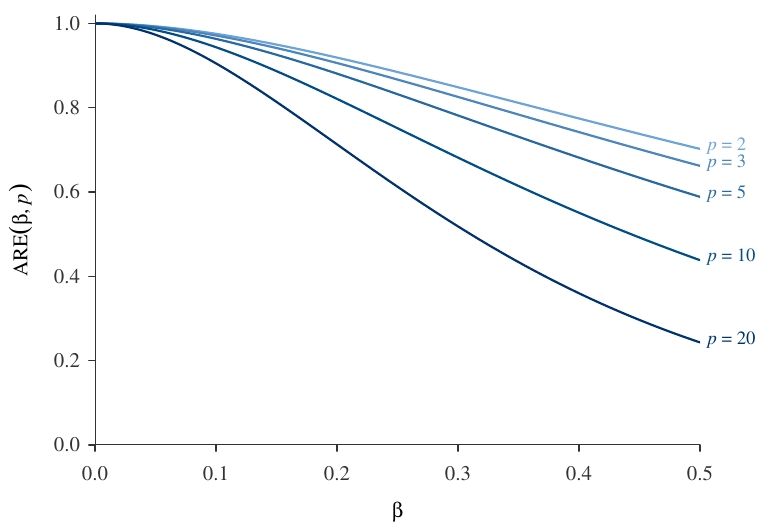}}

}

\caption{\label{fig-are}Asymptotic relative efficiency
\(\mathrm{ARE}(\beta,p)\) of the Rao's \(\beta\)-score tests with
respect to the classical Rao's score tests, as a function of \(\beta\),
for several values of the dimension \(p\).}

\end{figure}

Table~\ref{tbl-are} and Figure~\ref{fig-are} show
\(\mathrm{ARE}(\beta,p)\) for some values of \(\beta\) and \(p\),
together with the largest value of \(\beta\) for which the efficiency is
at least \(0.95\) and \(0.90\). The loss of efficiency for a fixed
\(\beta\) increases with the dimension, so the tuning parameter should
be chosen smaller as \(p\) increases.

\begin{longtable}[]{@{}
  >{\raggedright\arraybackslash}p{(\linewidth - 10\tabcolsep) * \real{0.1667}}
  >{\centering\arraybackslash}p{(\linewidth - 10\tabcolsep) * \real{0.1667}}
  >{\centering\arraybackslash}p{(\linewidth - 10\tabcolsep) * \real{0.1667}}
  >{\centering\arraybackslash}p{(\linewidth - 10\tabcolsep) * \real{0.1667}}
  >{\centering\arraybackslash}p{(\linewidth - 10\tabcolsep) * \real{0.1667}}
  >{\centering\arraybackslash}p{(\linewidth - 10\tabcolsep) * \real{0.1667}}@{}}
\caption{Asymptotic relative efficiency \(\mathrm{ARE}(\beta,p)\) of the
Rao's \(\beta\)-score tests with respect to the classical Rao's score
tests.}\label{tbl-are}\tabularnewline
\toprule\noalign{}
\begin{minipage}[b]{\linewidth}\raggedright
\end{minipage} & \begin{minipage}[b]{\linewidth}\centering
\(p=2\)
\end{minipage} & \begin{minipage}[b]{\linewidth}\centering
\(p=3\)
\end{minipage} & \begin{minipage}[b]{\linewidth}\centering
\(p=5\)
\end{minipage} & \begin{minipage}[b]{\linewidth}\centering
\(p=10\)
\end{minipage} & \begin{minipage}[b]{\linewidth}\centering
\(p=20\)
\end{minipage} \\
\midrule\noalign{}
\endfirsthead
\toprule\noalign{}
\begin{minipage}[b]{\linewidth}\raggedright
\end{minipage} & \begin{minipage}[b]{\linewidth}\centering
\(p=2\)
\end{minipage} & \begin{minipage}[b]{\linewidth}\centering
\(p=3\)
\end{minipage} & \begin{minipage}[b]{\linewidth}\centering
\(p=5\)
\end{minipage} & \begin{minipage}[b]{\linewidth}\centering
\(p=10\)
\end{minipage} & \begin{minipage}[b]{\linewidth}\centering
\(p=20\)
\end{minipage} \\
\midrule\noalign{}
\endhead
\bottomrule\noalign{}
\endlastfoot
\(\beta=0.05\) & 0.993 & 0.992 & 0.990 & 0.984 & 0.973 \\
\(\beta=0.1\) & 0.975 & 0.971 & 0.963 & 0.944 & 0.905 \\
\(\beta=0.2\) & 0.919 & 0.906 & 0.881 & 0.821 & 0.713 \\
\(\beta=0.3\) & 0.849 & 0.826 & 0.782 & 0.682 & 0.519 \\
\(\beta=0.5\) & 0.702 & 0.662 & 0.589 & 0.438 & 0.243 \\
largest \(\beta\) with \(\mathrm{ARE}(\beta,p)\geq0.95\) & 0.150 & 0.137
& 0.119 & 0.093 & 0.070 \\
largest \(\beta\) with \(\mathrm{ARE}(\beta,p)\geq0.90\) & 0.228 & 0.208
& 0.179 & 0.139 & 0.103 \\
\end{longtable}

\begin{refremark}
Let \(\boldsymbol{X}_{1},\ldots,\boldsymbol{X}_{n}\) be independent and
identically distributed with an elliptical distribution with location
\(\boldsymbol{\mu}_{0}\), scatter matrix
\(\boldsymbol{\Sigma}_{0}=\boldsymbol{\Lambda}_{0}^{1/2}\boldsymbol{R}_{0}\boldsymbol{\Lambda}_{0}^{1/2}\)
and radial density \(g\), not necessarily normal. If \(\boldsymbol{A}\)
is nonsingular and \(\boldsymbol{b}\in\mathbb{R}^{p}\), the objective
function (\ref{eq-Hn}) at the distribution of
\(\boldsymbol{A}\boldsymbol{X}_{i}+\boldsymbol{b}\) and at the parameter
\((\boldsymbol{A}\boldsymbol{\mu}+\boldsymbol{b},\boldsymbol{A}\boldsymbol{\Sigma}\boldsymbol{A}^{\top})\)
is \(\left\vert \boldsymbol{A}\right\vert ^{-\beta}\) times its value at
the distribution of \(\boldsymbol{X}_{i}\) and at
\((\boldsymbol{\mu},\boldsymbol{\Sigma})\), so the minimum DPD
functional is affine equivariant. Taking
\(\boldsymbol{A}=\boldsymbol{\Sigma}_{0}^{-1/2}\) the distribution
becomes spherical and the functional, being equivariant under orthogonal
transformations, is \((\boldsymbol{0}_{p},c\boldsymbol{I}_{p})\) for
some \(c>0\). Therefore the minimum DPD functional of
\((\boldsymbol{\mu},\boldsymbol{\Sigma})\) is
\((\boldsymbol{\mu}_{0},c_{\beta}(g)\boldsymbol{\Sigma}_{0})\), the
correlation matrix is estimated without asymptotic bias over the whole
elliptical family, and the tests of Section \ref{sec-Sec3} remain
consistent outside the multivariate normal model.

The hypotheses of Section \ref{sec-Sec3} are moreover invariant under a
common scale change of the scatter matrix, and the restricted minimum
DPD estimators are affine equivariant M-estimators of scatter, so by
Theorem 1 of Tyler (1983) their asymptotic variance-covariance matrix
enters the test-statistics only through a scalar
\(\sigma_{1,\beta}(g)\), and \[
\widetilde{R}_{\beta,n}\underset{n\rightarrow\infty}{\overset{\mathcal{L}}{\longrightarrow}}\frac{\sigma_{1,\beta}(g)}{\sigma_{1,\beta}(\phi)}\chi_{r}^{2},
\] under the null hypothesis, \(\phi\) denoting the radial density of
the multivariate normal distribution. In particular
\(\mathrm{ARE}(\beta,p)=\sigma_{1,0}(\phi)/\sigma_{1,\beta}(\phi)\),
which explains why it depends neither on \(\boldsymbol{R}_{0}\), nor on
\(\boldsymbol{\Delta}\), nor on the nuisance parameters. The same scalar
governs the limiting distribution of the likelihood ratio test outside
the multivariate normal model, where it reduces to \(1+\kappa\), with
\(\kappa\) the kurtosis parameter (Muirhead and Waternaux, 1980; Tyler,
1983).

\label{rem-elliptical}

\end{refremark}

\subsection{Influence functions}\label{sec-IF}

Let \(G\) be a distribution on \(\mathbb{R}^{p}\). The functional
associated with the restricted minimum DPD estimator of the nuisance
parameters, \(\widetilde{\boldsymbol{T}}_{\beta}(G)\), is the solution
in \(\boldsymbol{\eta}_{(1)}\) of
\(\mathrm{E}_{G}[\boldsymbol{u}_{\beta,1}(\boldsymbol{X},(\boldsymbol{\eta}_{(1)},\mathrm{vecl}(\boldsymbol{R}_{0})))]=\boldsymbol{0}\)
which is equal to \(\boldsymbol{\eta}_{(1),0}\) at
\(G=F_{\boldsymbol{\theta}_{0}}\), where \(\boldsymbol{u}_{\beta,1}\) is
formed by the components of the \(\beta\)-score function which
correspond to \(\boldsymbol{\eta}_{(1)}\), and the functional associated
with the test-statistic is \[
\widetilde{R}_{\beta}(G)=\overline{\boldsymbol{U}}_{\beta}^{\top}(G)\boldsymbol{K}_{\beta}^{-1}(\widetilde{\boldsymbol{\eta}}_{\beta}(G))\overline{\boldsymbol{U}}_{\beta}(G),\qquad
\overline{\boldsymbol{U}}_{\beta}(G)=\mathrm{E}_{G}[\boldsymbol{u}_{\beta}(\boldsymbol{X},\widetilde{\boldsymbol{\eta}}_{\beta}(G))],
\] with
\(\widetilde{\boldsymbol{\eta}}_{\beta}(G)=(\widetilde{\boldsymbol{T}}_{\beta}^{\top}(G),\mathrm{vecl}^{\top}(\boldsymbol{R}_{0}))^{\top}\),
so that \(\widetilde{R}_{\beta,n}=n\widetilde{R}_{\beta}(G_{n})\) for
the empirical distribution \(G_{n}\). We denote by
\(\wedge_{\boldsymbol{x}}\) the distribution degenerate at
\(\boldsymbol{x}\) and
\(F_{\epsilon,\boldsymbol{x}}=(1-\epsilon)F_{\boldsymbol{\theta}_{0}}+\epsilon\wedge_{\boldsymbol{x}}\),
and the influence function of a functional \(T\) at
\(F_{\boldsymbol{\theta}_{0}}\) is
\(\mathrm{IF}(\boldsymbol{x},T,F_{\boldsymbol{\theta}_{0}})=\left.\frac{\partial}{\partial\epsilon}T(F_{\epsilon,\boldsymbol{x}})\right\vert _{\epsilon=0}\).

\begin{theorem}[]\protect\hypertarget{thm-IFest}{}\label{thm-IFest}

The influence function of the restricted minimum DPD estimator of the
nuisance parameters is
\begin{equation}\protect\phantomsection\label{eq-IFgeneral}{
\mathrm{IF}(\boldsymbol{x},\widetilde{\boldsymbol{T}}_{\beta},F_{\boldsymbol{\theta}_{0}})=\boldsymbol{J}_{\beta,11}^{-1}(\boldsymbol{\eta}_{0})\boldsymbol{u}_{\beta,1}(\boldsymbol{x},\boldsymbol{\eta}_{0}).
}\end{equation} In particular, for a fixed correlation matrix
\(\boldsymbol{R}_{0}\) (Sections \ref{sec-test1}, \ref{sec-test2} and
\ref{sec-test3}), the influence functions of the estimators of
\(\boldsymbol{\mu}\) and of the vector of variances
\(\boldsymbol{\Lambda}\boldsymbol{1}_{p}\) are
\begin{equation}\protect\phantomsection\label{eq-IFexplicit}{
\begin{aligned}
\mathrm{IF}(\boldsymbol{x},\widetilde{\boldsymbol{\mu}}_{\beta},F_{\boldsymbol{\theta}_{0}})  &  =(\beta+1)^{\frac{p}{2}+1}w_{\beta}(\boldsymbol{x})(\boldsymbol{x}-\boldsymbol{\mu}_{0}),\\
\mathrm{IF}(\boldsymbol{x},\widetilde{\boldsymbol{\Lambda}}_{\beta}\boldsymbol{1}_{p},F_{\boldsymbol{\theta}_{0}})  &  =2(\beta+1)^{\frac{p}{2}+2}\boldsymbol{\Lambda}_{0}\left[\boldsymbol{d}(\boldsymbol{\Psi}_{\beta}(\boldsymbol{x}))-\frac{\beta^{2}\,\mathrm{trace}(\boldsymbol{R}_{0}^{-1}\boldsymbol{\Psi}_{\beta}(\boldsymbol{x}))}{2(2+\beta^{2}p)}\boldsymbol{1}_{p}\right],
\end{aligned}
}\end{equation} with \(\boldsymbol{\Psi}_{\beta}(\boldsymbol{x})\) given
by (\ref{eq-Psi}) and \(\boldsymbol{d}(\cdot)\) by (\ref{eq-adjust}).
The influence function (\ref{eq-IFgeneral}) is bounded in
\(\boldsymbol{x}\) if and only if \(\beta>0\), and in this case it tends
to a constant vector, zero for the components of \(\boldsymbol{\mu}\),
as \(\left\Vert \boldsymbol{x}\right\Vert \rightarrow\infty\).

\end{theorem}

For \(\beta=0\) and \(\boldsymbol{R}_{0}=\boldsymbol{I}_{p}\),
(\ref{eq-IFexplicit}) reduces to the influence functions of the sample
mean and of the sample variances,
\(\boldsymbol{x}-\boldsymbol{\mu}_{0}\) and
\(((x_{j}-\mu_{j,0})^{2}-\sigma_{j,0}^{2})_{j=1,\ldots,p}\).

\begin{theorem}[]\protect\hypertarget{thm-IF2}{}\label{thm-IF2}

The first-order influence function of the test functional vanishes under
the null hypothesis,
\(\mathrm{IF}(\boldsymbol{x},\widetilde{R}_{\beta},F_{\boldsymbol{\theta}_{0}})=0\),
and its second-order influence function, or self-standardized influence
function, is \begin{equation}\protect\phantomsection\label{eq-IF2}{
\mathrm{IF}^{(2)}(\boldsymbol{x},\widetilde{R}_{\beta},F_{\boldsymbol{\theta}_{0}})=\frac{1}{2}\left.\frac{\partial^{2}}{\partial\epsilon^{2}}\widetilde{R}_{\beta}(F_{\epsilon,\boldsymbol{x}})\right\vert _{\epsilon=0}
=\frac{1}{2}(1+2\beta)^{\frac{p}{2}+2}\mathrm{trace}\left(\left(\boldsymbol{R}_{0}^{-1}\boldsymbol{\Psi}_{\beta}^{\ast}(\boldsymbol{x})\right)^{2}\right).
}\end{equation} It is bounded in \(\boldsymbol{x}\) if and only if
\(\beta>0\). In particular, for testing complete uncorrelatedness,
\begin{equation}\protect\phantomsection\label{eq-IF2ind}{
\begin{aligned}
\mathrm{IF}^{(2)}(\boldsymbol{x},\widetilde{R}_{\beta},F_{\boldsymbol{\theta}_{0}})&=(1+2\beta)^{\frac{p}{2}+2}\exp\left\{-\beta\boldsymbol{z}^{\top}(\boldsymbol{x})\boldsymbol{z}(\boldsymbol{x})\right\}\sum_{i<j}z_{i}^{2}(\boldsymbol{x})z_{j}^{2}(\boldsymbol{x}),\\
\sup_{\boldsymbol{x}\in\mathbb{R}^{p}}\mathrm{IF}^{(2)}(\boldsymbol{x},\widetilde{R}_{\beta},F_{\boldsymbol{\theta}_{0}})&=\frac{2(p-1)(1+2\beta)^{\frac{p}{2}+2}}{pe^{2}\beta^{2}},
\end{aligned}
}\end{equation} for \(\beta>0\), the supremum being the gross-error
sensitivity of the test, while for \(\beta=0\) it is the unbounded
function
\(\sum_{i<j}z_{i}^{2}(\boldsymbol{x})z_{j}^{2}(\boldsymbol{x})\).

\end{theorem}

\begin{figure}[!htb]

\centering{

\pandocbounded{\includegraphics[keepaspectratio]{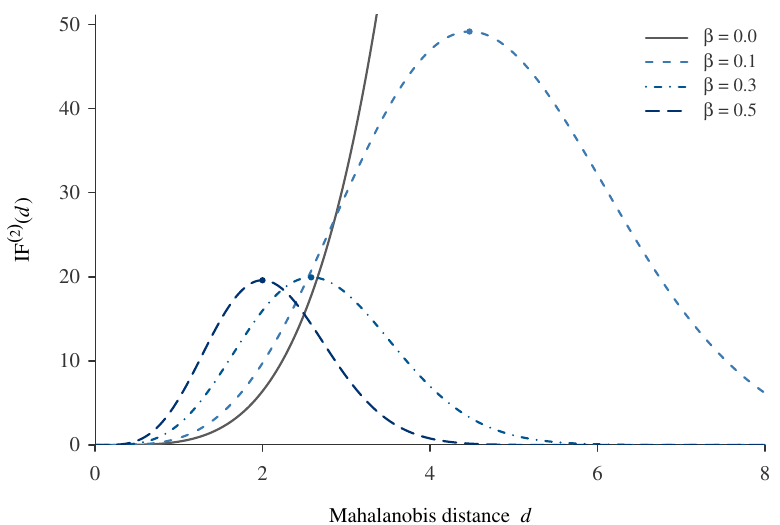}}

}

\caption{\label{fig-if2}Upper envelope of the second-order influence
function (\ref{eq-IF2ind}) of the test of complete uncorrelatedness for
\(p=5\), as a function of the Mahalanobis distance \(d\). The dots mark
the maxima, attained at \(d=\sqrt{2/\beta}\).}

\end{figure}

Figure~\ref{fig-if2} displays (\ref{eq-IF2ind}) for \(p=5\) along the
direction of equal components of \(\boldsymbol{z}(\boldsymbol{x})\), on
which \(\sum_{i<j}z_{i}^{2}(\boldsymbol{x})z_{j}^{2}(\boldsymbol{x})\)
attains its maximum for a given value of the Mahalanobis distance
\(d=\left\Vert \boldsymbol{z}(\boldsymbol{x})\right\Vert\), so that the
curves are the upper envelope of the second-order influence function.
For \(\beta=0\) the function increases without bound, while for
\(\beta>0\) it attains its maximum at \(d=\sqrt{2/\beta}\) and decreases
to zero afterwards, so that a sufficiently remote observation has no
effect on the test-statistic. The same shape is obtained for other
values of \(p\), with the vertical scale multiplied by
\((p-1)(1+2\beta)^{\frac{p}{2}+2}/p\).

The level and power influence functions (Hampel et al., 1986; Heritier
and Ronchetti, 1994) describe the effect of an infinitesimal
contamination on the asymptotic level and power of the test. Let
\(\boldsymbol{\theta}_{n}\) be the value of the parameter under
(\ref{eq-H1n}) and consider the contaminated distributions \[
F_{n,\epsilon,\boldsymbol{x}}^{P}=\left(1-\frac{\epsilon}{\sqrt{n}}\right)F_{\boldsymbol{\theta}_{n}}+\frac{\epsilon}{\sqrt{n}}\wedge_{\boldsymbol{x}},
\qquad
F_{n,\epsilon,\boldsymbol{x}}^{L}=\left(1-\frac{\epsilon}{\sqrt{n}}\right)F_{\boldsymbol{\theta}_{0}}+\frac{\epsilon}{\sqrt{n}}\wedge_{\boldsymbol{x}},
\] the second one being the first one with
\(\boldsymbol{\Delta}=\boldsymbol{0}_{p\times p}\). For a test of
asymptotic size \(\alpha\), with critical value \(\chi_{r,\alpha}^{2}\),
the quantile of order \(1-\alpha\) of the \(\chi_{r}^{2}\) distribution,
let
\(\pi_{\epsilon,\beta}(\boldsymbol{x})=\lim_{n\rightarrow\infty}P_{F_{n,\epsilon,\boldsymbol{x}}^{P}}(\widetilde{R}_{\beta,n}>\chi_{r,\alpha}^{2})\)
and
\(\alpha_{\epsilon,\beta}(\boldsymbol{x})=\lim_{n\rightarrow\infty}P_{F_{n,\epsilon,\boldsymbol{x}}^{L}}(\widetilde{R}_{\beta,n}>\chi_{r,\alpha}^{2})\)
be the asymptotic power and level under contamination. The power and
level influence functions are
\(\mathrm{PIF}_{\beta}(\boldsymbol{x})=\left.\frac{\partial}{\partial\epsilon}\pi_{\epsilon,\beta}(\boldsymbol{x})\right\vert _{\epsilon=0}\)
and
\(\mathrm{LIF}_{\beta}(\boldsymbol{x})=\left.\frac{\partial}{\partial\epsilon}\alpha_{\epsilon,\beta}(\boldsymbol{x})\right\vert _{\epsilon=0}\).

\begin{theorem}[]\protect\hypertarget{thm-PIF}{}\label{thm-PIF}

If the observations are independent and identically distributed with
distribution \(F_{n,\epsilon,\boldsymbol{x}}^{P}\), then, for every
\(\beta\geq0\),
\begin{equation}\protect\phantomsection\label{eq-deltastar}{
\widetilde{R}_{\beta,n}\underset{n\rightarrow\infty}{\overset{\mathcal{L}}{\longrightarrow}}\chi_{r}^{2}\left(\delta_{\beta}^{\ast}(\epsilon,\boldsymbol{x})\right),
}\end{equation} with \[
\delta_{\beta}^{\ast}(\epsilon,\boldsymbol{x})=\frac{1}{2}\mathrm{ARE}(\beta,p)\,\mathrm{trace}\left(\left(\boldsymbol{R}_{0}^{-1}\boldsymbol{\Delta}_{\beta}^{\ast}(\epsilon,\boldsymbol{x})\right)^{2}\right),
\qquad
\boldsymbol{\Delta}_{\beta}(\epsilon,\boldsymbol{x})=\boldsymbol{\Delta}+\epsilon(1+\beta)^{\frac{p}{2}+2}\boldsymbol{\Psi}_{\beta}(\boldsymbol{x}),
\] or equivalently,
\begin{equation}\protect\phantomsection\label{eq-deltaexp}{
\delta_{\beta}^{\ast}(\epsilon,\boldsymbol{x})=\delta_{\beta}(\boldsymbol{\Delta})+\epsilon\left(\frac{1+2\beta}{1+\beta}\right)^{\frac{p}{2}+2}\mathrm{trace}\left(\boldsymbol{R}_{0}^{-1}\boldsymbol{\Delta}^{\ast}\boldsymbol{R}_{0}^{-1}\boldsymbol{\Psi}_{\beta}^{\ast}(\boldsymbol{x})\right)+\epsilon^{2}\mathrm{IF}^{(2)}(\boldsymbol{x},\widetilde{R}_{\beta},F_{\boldsymbol{\theta}_{0}}).
}\end{equation} Consequently,
\begin{equation}\protect\phantomsection\label{eq-PIF}{
\mathrm{PIF}_{\beta}(\boldsymbol{x})=\left(\frac{1+2\beta}{1+\beta}\right)^{\frac{p}{2}+2}f_{\chi_{r+2}^{2}(\delta_{\beta}(\boldsymbol{\Delta}))}\left(\chi_{r,\alpha}^{2}\right)\mathrm{trace}\left(\boldsymbol{R}_{0}^{-1}\boldsymbol{\Delta}^{\ast}\boldsymbol{R}_{0}^{-1}\boldsymbol{\Psi}_{\beta}^{\ast}(\boldsymbol{x})\right),
}\end{equation} where \(f_{\chi_{r+2}^{2}(\delta)}\) is the density
function of the \(\chi_{r+2}^{2}(\delta)\) distribution,
\(\mathrm{LIF}_{\beta}(\boldsymbol{x})=0\), and the asymptotic level
under contamination is
\begin{equation}\protect\phantomsection\label{eq-levelcont}{
\alpha_{\epsilon,\beta}(\boldsymbol{x})=P\left(\chi_{r}^{2}\left(\epsilon^{2}\mathrm{IF}^{(2)}(\boldsymbol{x},\widetilde{R}_{\beta},F_{\boldsymbol{\theta}_{0}})\right)>\chi_{r,\alpha}^{2}\right).
}\end{equation} If
\(\boldsymbol{\Delta}^{\ast}\neq\boldsymbol{0}_{p\times p}\), the power
influence function is bounded in \(\boldsymbol{x}\) if and only if
\(\beta>0\). For \(\beta>0\),
\(\sup_{\boldsymbol{x}\in\mathbb{R}^{p}}\alpha_{\epsilon,\beta}(\boldsymbol{x})<1\)
for every \(\epsilon\), while for \(\beta=0\) this supremum is equal to
one for every \(\epsilon\neq0\).

\end{theorem}

\section{Simulation study}\label{sec-Sim}

The finite-sample behaviour of the test-statistics is examined by Monte
Carlo simulation, with \(B=1000\) replications in every scenario and
nominal level \(\alpha=0.05\), so that the Monte Carlo standard error of
a proportion of rejections close to the nominal level is \(0.007\). The
data are generated from a \(p\)-variate distribution with unequal
variances, \(\sigma_{j}^{2}=\exp\{-1+2(j-1)/(p-1)\}\), \(j=1,\ldots,p\),
so that the correlation matrix and the variances are not confounded, and
the specified correlation matrix is the autoregressive one,
\(R_{0,ij}=0.5^{\left\vert i-j\right\vert }\). The tuning parameter
takes the values \(\beta=0\), \(0.1\), \(0.2\), \(0.3\) and \(0.5\), the
first of them being the classical Rao's score test. For each hypothesis
the classical competitors of Section \ref{sec-Discuss} are included,
namely the likelihood ratio test in the four cases, the test of Jennrich
(1970) for a specified correlation matrix, the tests of Bartlett (1954)
and the one based on the Spearman rank correlation coefficients for
complete uncorrelatedness, and the test of Lawley (1963) for an
equicorrelation with an unknown common coefficient. The restricted
estimators are computed as described in Remark~\ref{rem-rootpractice},
from the classical and from the robust starting values, and in every
replication the two resulting values of (\ref{eq-Hn}) are compared in
order to check that the same root is reached.

\subsection{Size}\label{sec-SimSize}

Table~\ref{tbl-size} reports the simulated size. With \(n=400\) every
test-statistic keeps the nominal level, the differences with respect to
\(0.05\) being of the order of the Monte Carlo standard error. The
differences appear as the sample size decreases, and they do so in the
direction that favours the test-statistics of this paper. With \(n=50\)
and \(p=5\) the likelihood ratio test is liberal under the four
hypotheses, \(0.070\), \(0.065\), \(0.085\) and \(0.075\), and the test
of Jennrich (1970) reaches \(0.075\), whereas the Rao's \(\beta\)-score
test-statistics keep the nominal level under three of the four. The
corrections of Bartlett (1954) and Lawley (1963), which were devised
precisely for this difficulty, do most of their work but not all of it,
\(0.062\) and \(0.058\) against the \(0.085\) and \(0.075\) of the
uncorrected likelihood ratio test. The test based on the Spearman rank
correlation coefficients, which is not reported in Table~\ref{tbl-size},
keeps the nominal level in the six scenarios of complete
uncorrelatedness, between \(0.039\) and \(0.058\).

Two qualifications are in order. Under complete uncorrelatedness with
\(n=50\) and \(p=5\) the Rao's \(\beta\)-score test-statistics are
liberal as well, between \(0.049\) and \(0.067\), so that the advantage
there is only with respect to the uncorrected likelihood ratio test. And
with a specified correlation matrix, \(p=3\) and \(n=50\), the value
\(0.071\) attained at \(\beta=0.5\) is three Monte Carlo standard errors
above the nominal level. The behaviour is therefore not uniform, and
what the simulation supports is the weaker statement that, in the
scenarios considered, no test-statistic of the family is more liberal
than its classical counterpart.

It is worth noting that this holds already at \(\beta=0\), so that it is
not an effect of robustness but of the form of the test-statistic, which
involves only the estimators restricted to the null hypothesis. That is
the argument of Rao (1948) himself, and with \(n=50\) and \(2p+1=11\)
parameters it stops being a computational convenience and becomes a
statistical advantage.

\begin{longtable}[]{@{}
  >{\raggedright\arraybackslash}p{(\linewidth - 18\tabcolsep) * \real{0.2292}}
  >{\centering\arraybackslash}p{(\linewidth - 18\tabcolsep) * \real{0.0833}}
  >{\centering\arraybackslash}p{(\linewidth - 18\tabcolsep) * \real{0.1042}}
  >{\centering\arraybackslash}p{(\linewidth - 18\tabcolsep) * \real{0.0833}}
  >{\centering\arraybackslash}p{(\linewidth - 18\tabcolsep) * \real{0.0833}}
  >{\centering\arraybackslash}p{(\linewidth - 18\tabcolsep) * \real{0.0833}}
  >{\centering\arraybackslash}p{(\linewidth - 18\tabcolsep) * \real{0.0833}}
  >{\centering\arraybackslash}p{(\linewidth - 18\tabcolsep) * \real{0.0833}}
  >{\centering\arraybackslash}p{(\linewidth - 18\tabcolsep) * \real{0.0833}}
  >{\centering\arraybackslash}p{(\linewidth - 18\tabcolsep) * \real{0.0833}}@{}}
\caption{Simulated size of the tests at the nominal level \(0.05\), with
\(B=1000\) replications. The last column is the test of Jennrich (1970)
for a specified correlation matrix, the one of Bartlett (1954) for
complete uncorrelatedness and the one of Lawley (1963) for an
equicorrelation with an unknown common coefficient. The hypotheses are
those of Sections \ref{sec-test1} to \ref{sec-test4}, in this
order.}\label{tbl-size}\tabularnewline
\toprule\noalign{}
\begin{minipage}[b]{\linewidth}\raggedright
hypothesis
\end{minipage} & \begin{minipage}[b]{\linewidth}\centering
\(p\)
\end{minipage} & \begin{minipage}[b]{\linewidth}\centering
\(n\)
\end{minipage} & \begin{minipage}[b]{\linewidth}\centering
\(\beta=0\)
\end{minipage} & \begin{minipage}[b]{\linewidth}\centering
\(0.1\)
\end{minipage} & \begin{minipage}[b]{\linewidth}\centering
\(0.2\)
\end{minipage} & \begin{minipage}[b]{\linewidth}\centering
\(0.3\)
\end{minipage} & \begin{minipage}[b]{\linewidth}\centering
\(0.5\)
\end{minipage} & \begin{minipage}[b]{\linewidth}\centering
LRT
\end{minipage} & \begin{minipage}[b]{\linewidth}\centering
other
\end{minipage} \\
\midrule\noalign{}
\endfirsthead
\toprule\noalign{}
\begin{minipage}[b]{\linewidth}\raggedright
hypothesis
\end{minipage} & \begin{minipage}[b]{\linewidth}\centering
\(p\)
\end{minipage} & \begin{minipage}[b]{\linewidth}\centering
\(n\)
\end{minipage} & \begin{minipage}[b]{\linewidth}\centering
\(\beta=0\)
\end{minipage} & \begin{minipage}[b]{\linewidth}\centering
\(0.1\)
\end{minipage} & \begin{minipage}[b]{\linewidth}\centering
\(0.2\)
\end{minipage} & \begin{minipage}[b]{\linewidth}\centering
\(0.3\)
\end{minipage} & \begin{minipage}[b]{\linewidth}\centering
\(0.5\)
\end{minipage} & \begin{minipage}[b]{\linewidth}\centering
LRT
\end{minipage} & \begin{minipage}[b]{\linewidth}\centering
other
\end{minipage} \\
\midrule\noalign{}
\endhead
\bottomrule\noalign{}
\endlastfoot
\(\boldsymbol{R}=\boldsymbol{R}_{0}\) & 3 & 50 & 0.053 & 0.057 & 0.063 &
0.067 & 0.071 & 0.068 & 0.066 \\
\(\boldsymbol{R}=\boldsymbol{R}_{0}\) & 3 & 100 & 0.051 & 0.060 & 0.064
& 0.059 & 0.058 & 0.057 & 0.058 \\
\(\boldsymbol{R}=\boldsymbol{R}_{0}\) & 3 & 400 & 0.039 & 0.041 & 0.049
& 0.049 & 0.049 & 0.043 & 0.043 \\
\(\boldsymbol{R}=\boldsymbol{R}_{0}\) & 5 & 50 & 0.053 & 0.054 & 0.046 &
0.043 & 0.046 & 0.070 & 0.075 \\
\(\boldsymbol{R}=\boldsymbol{R}_{0}\) & 5 & 100 & 0.049 & 0.046 & 0.046
& 0.051 & 0.044 & 0.060 & 0.070 \\
\(\boldsymbol{R}=\boldsymbol{R}_{0}\) & 5 & 400 & 0.052 & 0.052 & 0.049
& 0.049 & 0.048 & 0.055 & 0.060 \\
\(\rho_{ij}=\rho_{0}\) & 3 & 50 & 0.050 & 0.048 & 0.046 & 0.050 & 0.051
& 0.057 & --- \\
\(\rho_{ij}=\rho_{0}\) & 3 & 100 & 0.053 & 0.054 & 0.060 & 0.052 & 0.053
& 0.061 & --- \\
\(\rho_{ij}=\rho_{0}\) & 3 & 400 & 0.045 & 0.048 & 0.044 & 0.047 & 0.041
& 0.047 & --- \\
\(\rho_{ij}=\rho_{0}\) & 5 & 50 & 0.050 & 0.046 & 0.037 & 0.042 & 0.041
& 0.065 & --- \\
\(\rho_{ij}=\rho_{0}\) & 5 & 100 & 0.056 & 0.061 & 0.064 & 0.053 & 0.051
& 0.065 & --- \\
\(\rho_{ij}=\rho_{0}\) & 5 & 400 & 0.043 & 0.045 & 0.054 & 0.057 & 0.048
& 0.045 & --- \\
\(\boldsymbol{R}=\boldsymbol{I}_{p}\) & 3 & 50 & 0.048 & 0.051 & 0.041 &
0.042 & 0.040 & 0.054 & 0.043 \\
\(\boldsymbol{R}=\boldsymbol{I}_{p}\) & 3 & 100 & 0.038 & 0.039 & 0.042
& 0.041 & 0.040 & 0.046 & 0.040 \\
\(\boldsymbol{R}=\boldsymbol{I}_{p}\) & 3 & 400 & 0.048 & 0.052 & 0.057
& 0.058 & 0.059 & 0.051 & 0.049 \\
\(\boldsymbol{R}=\boldsymbol{I}_{p}\) & 5 & 50 & 0.067 & 0.063 & 0.060 &
0.061 & 0.049 & 0.085 & 0.062 \\
\(\boldsymbol{R}=\boldsymbol{I}_{p}\) & 5 & 100 & 0.039 & 0.040 & 0.034
& 0.038 & 0.036 & 0.050 & 0.038 \\
\(\boldsymbol{R}=\boldsymbol{I}_{p}\) & 5 & 400 & 0.042 & 0.045 & 0.046
& 0.047 & 0.048 & 0.045 & 0.043 \\
\(\rho_{ij}=\rho_{12}\) & 3 & 50 & 0.057 & 0.053 & 0.047 & 0.042 & 0.046
& 0.060 & 0.051 \\
\(\rho_{ij}=\rho_{12}\) & 3 & 100 & 0.048 & 0.046 & 0.045 & 0.046 &
0.043 & 0.054 & 0.046 \\
\(\rho_{ij}=\rho_{12}\) & 3 & 400 & 0.056 & 0.057 & 0.056 & 0.058 &
0.046 & 0.057 & 0.056 \\
\(\rho_{ij}=\rho_{12}\) & 5 & 50 & 0.055 & 0.050 & 0.050 & 0.050 & 0.042
& 0.075 & 0.058 \\
\(\rho_{ij}=\rho_{12}\) & 5 & 100 & 0.037 & 0.045 & 0.050 & 0.051 &
0.046 & 0.051 & 0.037 \\
\(\rho_{ij}=\rho_{12}\) & 5 & 400 & 0.053 & 0.057 & 0.059 & 0.055 &
0.039 & 0.054 & 0.052 \\
\end{longtable}

\subsection{Power}\label{sec-SimPower}

Table~\ref{tbl-power} reports the simulated power under the fixed
alternative obtained by adding \(\delta\) to the elements of the first
subdiagonal of the correlation matrix under the null hypothesis, for
\(p=3\). The power decreases with \(\beta\), as the asymptotic relative
efficiency of Corollary~\ref{cor-are} predicts, and the loss is moderate
for the smaller values of the tuning parameter; with \(\delta=0.10\) and
\(n=200\) it goes from \(0.871\) at \(\beta=0\) to \(0.866\) at
\(\beta=0.1\) and \(0.705\) at \(\beta=0.5\), against
\(\mathrm{ARE}(0.1,3)=0.971\) and \(\mathrm{ARE}(0.5,3)=0.662\). For
\(\delta\geq0.2\) all the test-statistics reject with probability close
to one and the comparison is uninformative.

The test of Jennrich (1970) is markedly less powerful than the Rao's
score test-statistic against these alternatives, \(0.071\) against
\(0.139\) with \(n=100\) and \(\delta=0.05\), and \(0.368\) against
\(0.544\) with \(\delta=0.10\), while the two have the same size in the
same scenario. The two test-statistics are asymptotically equivalent
under the null hypothesis, which says nothing about their power in
finite samples, and the difference is plausibly explained by the
standardization: the test-statistic of Jennrich (1970) starts from the
sample Pearson correlation matrix and corrects afterwards for the
estimation of the variances, whereas (\ref{eq-RaoR0}) starts from the
restricted estimators, in which the restriction is already built into
the standardization.

\begin{longtable}[]{@{}
  >{\raggedright\arraybackslash}p{(\linewidth - 18\tabcolsep) * \real{0.2340}}
  >{\centering\arraybackslash}p{(\linewidth - 18\tabcolsep) * \real{0.0851}}
  >{\centering\arraybackslash}p{(\linewidth - 18\tabcolsep) * \real{0.0851}}
  >{\centering\arraybackslash}p{(\linewidth - 18\tabcolsep) * \real{0.0851}}
  >{\centering\arraybackslash}p{(\linewidth - 18\tabcolsep) * \real{0.0851}}
  >{\centering\arraybackslash}p{(\linewidth - 18\tabcolsep) * \real{0.0851}}
  >{\centering\arraybackslash}p{(\linewidth - 18\tabcolsep) * \real{0.0851}}
  >{\centering\arraybackslash}p{(\linewidth - 18\tabcolsep) * \real{0.0851}}
  >{\centering\arraybackslash}p{(\linewidth - 18\tabcolsep) * \real{0.0851}}
  >{\centering\arraybackslash}p{(\linewidth - 18\tabcolsep) * \real{0.0851}}@{}}
\caption{Simulated power of the tests at the nominal level \(0.05\),
with \(p=3\) and \(B=1000\) replications, under the alternative obtained
by adding \(\delta\) to the elements of the first subdiagonal. The last
column is as in Table~\ref{tbl-size}. For \(\delta\geq0.2\) all the
proportions of rejections are close to
one.}\label{tbl-power}\tabularnewline
\toprule\noalign{}
\begin{minipage}[b]{\linewidth}\raggedright
hypothesis
\end{minipage} & \begin{minipage}[b]{\linewidth}\centering
\(n\)
\end{minipage} & \begin{minipage}[b]{\linewidth}\centering
\(\delta\)
\end{minipage} & \begin{minipage}[b]{\linewidth}\centering
\(\beta=0\)
\end{minipage} & \begin{minipage}[b]{\linewidth}\centering
\(0.1\)
\end{minipage} & \begin{minipage}[b]{\linewidth}\centering
\(0.2\)
\end{minipage} & \begin{minipage}[b]{\linewidth}\centering
\(0.3\)
\end{minipage} & \begin{minipage}[b]{\linewidth}\centering
\(0.5\)
\end{minipage} & \begin{minipage}[b]{\linewidth}\centering
LRT
\end{minipage} & \begin{minipage}[b]{\linewidth}\centering
other
\end{minipage} \\
\midrule\noalign{}
\endfirsthead
\toprule\noalign{}
\begin{minipage}[b]{\linewidth}\raggedright
hypothesis
\end{minipage} & \begin{minipage}[b]{\linewidth}\centering
\(n\)
\end{minipage} & \begin{minipage}[b]{\linewidth}\centering
\(\delta\)
\end{minipage} & \begin{minipage}[b]{\linewidth}\centering
\(\beta=0\)
\end{minipage} & \begin{minipage}[b]{\linewidth}\centering
\(0.1\)
\end{minipage} & \begin{minipage}[b]{\linewidth}\centering
\(0.2\)
\end{minipage} & \begin{minipage}[b]{\linewidth}\centering
\(0.3\)
\end{minipage} & \begin{minipage}[b]{\linewidth}\centering
\(0.5\)
\end{minipage} & \begin{minipage}[b]{\linewidth}\centering
LRT
\end{minipage} & \begin{minipage}[b]{\linewidth}\centering
other
\end{minipage} \\
\midrule\noalign{}
\endhead
\bottomrule\noalign{}
\endlastfoot
\(\boldsymbol{R}=\boldsymbol{R}_{0}\) & 100 & 0.05 & 0.139 & 0.136 &
0.127 & 0.124 & 0.101 & 0.149 & 0.071 \\
\(\boldsymbol{R}=\boldsymbol{R}_{0}\) & 100 & 0.10 & 0.544 & 0.535 &
0.498 & 0.449 & 0.388 & 0.584 & 0.368 \\
\(\boldsymbol{R}=\boldsymbol{R}_{0}\) & 100 & 0.15 & 0.979 & 0.975 &
0.965 & 0.945 & 0.876 & 0.987 & 0.928 \\
\(\boldsymbol{R}=\boldsymbol{R}_{0}\) & 200 & 0.05 & 0.259 & 0.264 &
0.245 & 0.221 & 0.180 & 0.275 & 0.169 \\
\(\boldsymbol{R}=\boldsymbol{R}_{0}\) & 200 & 0.10 & 0.871 & 0.866 &
0.844 & 0.816 & 0.705 & 0.884 & 0.797 \\
\(\boldsymbol{R}=\boldsymbol{R}_{0}\) & 200 & 0.15 & 1.000 & 1.000 &
1.000 & 1.000 & 0.994 & 1.000 & 1.000 \\
\(\boldsymbol{R}=\boldsymbol{I}_{p}\) & 100 & 0.05 & 0.084 & 0.077 &
0.075 & 0.076 & 0.065 & 0.085 & 0.081 \\
\(\boldsymbol{R}=\boldsymbol{I}_{p}\) & 100 & 0.10 & 0.201 & 0.195 &
0.180 & 0.155 & 0.136 & 0.215 & 0.207 \\
\(\boldsymbol{R}=\boldsymbol{I}_{p}\) & 100 & 0.15 & 0.421 & 0.417 &
0.396 & 0.347 & 0.280 & 0.439 & 0.428 \\
\(\boldsymbol{R}=\boldsymbol{I}_{p}\) & 200 & 0.05 & 0.119 & 0.114 &
0.103 & 0.098 & 0.086 & 0.122 & 0.116 \\
\(\boldsymbol{R}=\boldsymbol{I}_{p}\) & 200 & 0.10 & 0.356 & 0.344 &
0.327 & 0.301 & 0.256 & 0.366 & 0.360 \\
\(\boldsymbol{R}=\boldsymbol{I}_{p}\) & 200 & 0.15 & 0.720 & 0.715 &
0.686 & 0.654 & 0.537 & 0.728 & 0.722 \\
\(\rho_{ij}=\rho_{12}\) & 100 & 0.05 & 0.096 & 0.099 & 0.097 & 0.083 &
0.071 & 0.100 & 0.098 \\
\(\rho_{ij}=\rho_{12}\) & 100 & 0.10 & 0.305 & 0.299 & 0.267 & 0.242 &
0.204 & 0.320 & 0.307 \\
\(\rho_{ij}=\rho_{12}\) & 100 & 0.15 & 0.653 & 0.653 & 0.619 & 0.580 &
0.481 & 0.681 & 0.659 \\
\(\rho_{ij}=\rho_{12}\) & 200 & 0.05 & 0.143 & 0.136 & 0.132 & 0.129 &
0.107 & 0.153 & 0.145 \\
\(\rho_{ij}=\rho_{12}\) & 200 & 0.10 & 0.502 & 0.501 & 0.460 & 0.433 &
0.341 & 0.515 & 0.504 \\
\(\rho_{ij}=\rho_{12}\) & 200 & 0.15 & 0.937 & 0.930 & 0.906 & 0.877 &
0.778 & 0.944 & 0.937 \\
\end{longtable}

\subsection{Contamination}\label{sec-SimCont}

\begin{figure}[!htb]

\centering{

\pandocbounded{\includegraphics[keepaspectratio]{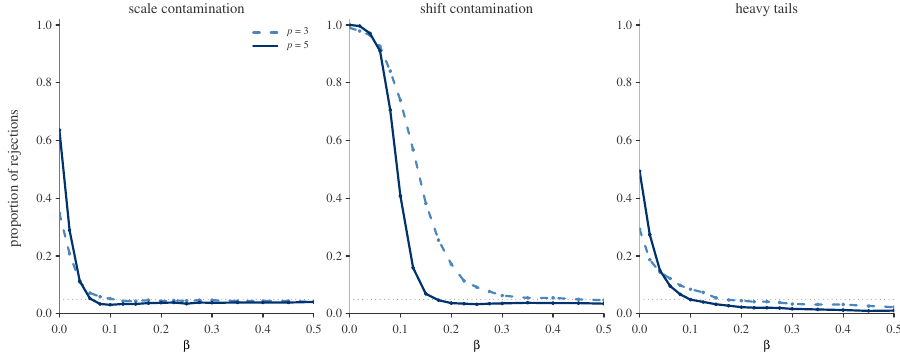}}

}

\caption{\label{fig-cont}Simulated proportion of rejections under
contamination for the test of complete uncorrelatedness, with
\(\epsilon=0.05\), \(n=200\) and \(B=1000\) replications, on a grid of
eleven values of \(\beta\). The dotted line is the nominal level
\(0.05\).}

\end{figure}

The third study is the one of interest, and is reported in
Table~\ref{tbl-cont} and displayed in Figure~\ref{fig-cont}. The null
hypothesis is true for the uncontaminated part of the data, so that the
proportions of rejections are sizes and measure the effect of the
outliers alone. Three schemes are considered, a fraction \(\epsilon\) of
the observations generated with the variances multiplied by a constant,
the same fraction shifted in the mean, and a multivariate Student
distribution with five degrees of freedom.

The classical test-statistic breaks down in the three of them, and the
more so as the dimension grows: with contamination in scale,
\(\epsilon=0.05\), \(n=200\) and a fixed equicorrelation, the proportion
of rejections is \(0.637\) at \(\beta=0\) and \(0.039\) at
\(\beta=0.1\), while the likelihood ratio test gives \(0.646\). The
breakdown of the likelihood ratio test and of the test of Bartlett
(1954) is indistinguishable from that of the classical Rao's score test,
as was to be expected, since the three are built on the sample Pearson
correlation matrix. The test based on the Spearman rank correlation
coefficients, which is distribution-free, resists the contamination in
scale, \(0.062\) against \(0.334\) for the test of Bartlett (1954) with
\(p=3\), but not the contamination in the mean, where it reaches
\(0.910\) with \(p=5\).

Figure~\ref{fig-cont}, computed on a grid of eighteen values of
\(\beta\) which is finer where the transition takes place, shows that
the value of the tuning parameter that is needed depends on the
contamination scheme, and by a factor of four. With contamination in
scale the nominal level is recovered at \(\beta\simeq0.08\) when
\(p=3\), the proportion of rejections falling from \(0.348\) to
\(0.208\), \(0.110\), \(0.072\) and \(0.059\) at
\(\beta=0.02,0.04,0.06\) and \(0.08\). With contamination in the mean it
is not recovered until \(\beta\simeq0.28\), the curve passing through
\(0.990\), \(0.965\), \(0.840\), \(0.568\), \(0.255\) and \(0.114\) at
\(\beta=0,0.04,0.08,0.125,0.175\) and \(0.225\). Read together with
Table~\ref{tbl-are}, this is the practical rule for the choice of the
tuning parameter: protection against contamination in scale is almost
free, \(\mathrm{ARE}(0.08,3)=0.981\), whereas protection against
contamination in the mean costs about one sixth of the asymptotic
efficiency, \(\mathrm{ARE}(0.28,3)=0.842\).

The effect of the dimension is twofold and the two parts must not be
confused. At \(\beta=0\) the breakdown worsens as \(p\) grows, \(0.348\)
against \(0.636\) with contamination in scale and \(0.293\) against
\(0.494\) under the Student distribution, which agrees with the mean
weight of an observation being \((1+\beta)^{-p/2}\). The value of the
tuning parameter needed to recover the nominal level, on the contrary,
decreases with the dimension under the three schemes,
\(\beta\simeq0.06\) against \(\beta\simeq0.08\) in scale,
\(\beta\simeq0.16\) against \(\beta\simeq0.28\) in the mean and
\(\beta\simeq0.10\) against \(\beta\simeq0.16\) under the Student
distribution. The reason is that a given contamination falls farther in
Mahalanobis distance the larger the dimension is, so that the weight
penalizes it sooner, and the two effects compensate each other to some
extent.

Under the Student distribution the proportion of rejections does not
settle at the nominal level but keeps decreasing, down to \(0.023\) with
\(p=3\) and to \(0.011\) with \(p=5\) at \(\beta=0.5\), a fifth of the
nominal level. This is not a failure of the test-statistics but a
consequence of Remark~\ref{rem-elliptical}: the constants
\(\kappa_{1}(p,\beta)\) and \(\kappa_{2}(p,\beta)\) are computed under
the multivariate normal distribution, so that over an elliptical family
the limiting distribution is a scalar multiple of the chi-square one,
and the scalar departs from one as \(\beta\) and as \(p\) grow, in the
same way as the asymptotic relative efficiency does. The test-statistics
remain valid in level, since the departure makes them conservative and
not liberal, at the price of some power.

Finally, in none of the scenarios of this section, and in none of the
forty combinations of hypothesis, dimension, sample size and tuning
parameter of Table~\ref{tbl-size}, did the two starting values of
Remark~\ref{rem-rootpractice} lead to different solutions of the
restricted estimating equations. This agrees with
Remark~\ref{rem-unbounded}, since the objective function is bounded from
below as soon as \(n\geq(\beta+1)^{\frac{p}{2}+1}/\beta\), which for
\(p=5\) and \(\beta=0.5\) amounts to \(8.3\) observations. The
diagnostic is therefore inactive in the designs considered here, and
becomes informative only for the much smaller sample sizes relative to
the dimension of Section \ref{sec-Appl}.

\begin{longtable}[]{@{}
  >{\raggedright\arraybackslash}p{(\linewidth - 18\tabcolsep) * \real{0.2245}}
  >{\raggedright\arraybackslash}p{(\linewidth - 18\tabcolsep) * \real{0.1224}}
  >{\centering\arraybackslash}p{(\linewidth - 18\tabcolsep) * \real{0.0816}}
  >{\centering\arraybackslash}p{(\linewidth - 18\tabcolsep) * \real{0.0816}}
  >{\centering\arraybackslash}p{(\linewidth - 18\tabcolsep) * \real{0.0816}}
  >{\centering\arraybackslash}p{(\linewidth - 18\tabcolsep) * \real{0.0816}}
  >{\centering\arraybackslash}p{(\linewidth - 18\tabcolsep) * \real{0.0816}}
  >{\centering\arraybackslash}p{(\linewidth - 18\tabcolsep) * \real{0.0816}}
  >{\centering\arraybackslash}p{(\linewidth - 18\tabcolsep) * \real{0.0816}}
  >{\centering\arraybackslash}p{(\linewidth - 18\tabcolsep) * \real{0.0816}}@{}}
\caption{Simulated proportion of rejections under contamination at the
nominal level \(0.05\), with \(n=200\) and \(B=1000\) replications. The
null hypothesis is true for the uncontaminated part of the data, so that
the entries are sizes. The tests of Bartlett (1954) and the one based on
the Spearman rank correlation coefficients are discussed in the
text.}\label{tbl-cont}\tabularnewline
\toprule\noalign{}
\begin{minipage}[b]{\linewidth}\raggedright
hypothesis
\end{minipage} & \begin{minipage}[b]{\linewidth}\raggedright
scheme
\end{minipage} & \begin{minipage}[b]{\linewidth}\centering
\(\epsilon\)
\end{minipage} & \begin{minipage}[b]{\linewidth}\centering
\(p\)
\end{minipage} & \begin{minipage}[b]{\linewidth}\centering
\(\beta=0\)
\end{minipage} & \begin{minipage}[b]{\linewidth}\centering
\(0.1\)
\end{minipage} & \begin{minipage}[b]{\linewidth}\centering
\(0.2\)
\end{minipage} & \begin{minipage}[b]{\linewidth}\centering
\(0.3\)
\end{minipage} & \begin{minipage}[b]{\linewidth}\centering
\(0.5\)
\end{minipage} & \begin{minipage}[b]{\linewidth}\centering
LRT
\end{minipage} \\
\midrule\noalign{}
\endfirsthead
\toprule\noalign{}
\begin{minipage}[b]{\linewidth}\raggedright
hypothesis
\end{minipage} & \begin{minipage}[b]{\linewidth}\raggedright
scheme
\end{minipage} & \begin{minipage}[b]{\linewidth}\centering
\(\epsilon\)
\end{minipage} & \begin{minipage}[b]{\linewidth}\centering
\(p\)
\end{minipage} & \begin{minipage}[b]{\linewidth}\centering
\(\beta=0\)
\end{minipage} & \begin{minipage}[b]{\linewidth}\centering
\(0.1\)
\end{minipage} & \begin{minipage}[b]{\linewidth}\centering
\(0.2\)
\end{minipage} & \begin{minipage}[b]{\linewidth}\centering
\(0.3\)
\end{minipage} & \begin{minipage}[b]{\linewidth}\centering
\(0.5\)
\end{minipage} & \begin{minipage}[b]{\linewidth}\centering
LRT
\end{minipage} \\
\midrule\noalign{}
\endhead
\bottomrule\noalign{}
\endlastfoot
\(\rho_{ij}=\rho_{0}\) & scale & 0.05 & 3 & 0.347 & 0.057 & 0.041 &
0.041 & 0.042 & 0.351 \\
\(\rho_{ij}=\rho_{0}\) & scale & 0.05 & 5 & 0.637 & 0.039 & 0.039 &
0.036 & 0.044 & 0.646 \\
\(\rho_{ij}=\rho_{0}\) & scale & 0.10 & 3 & 0.402 & 0.069 & 0.038 &
0.036 & 0.036 & 0.412 \\
\(\rho_{ij}=\rho_{0}\) & scale & 0.10 & 5 & 0.743 & 0.037 & 0.026 &
0.027 & 0.029 & 0.755 \\
\(\rho_{ij}=\rho_{0}\) & shift & 0.05 & 3 & 0.919 & 0.692 & 0.285 &
0.120 & 0.053 & 0.916 \\
\(\rho_{ij}=\rho_{0}\) & shift & 0.05 & 5 & 0.939 & 0.608 & 0.203 &
0.071 & 0.044 & 0.938 \\
\(\rho_{ij}=\rho_{0}\) & shift & 0.10 & 3 & 1.000 & 0.989 & 0.846 &
0.460 & 0.089 & 1.000 \\
\(\rho_{ij}=\rho_{0}\) & shift & 0.10 & 5 & 1.000 & 0.976 & 0.673 &
0.227 & 0.047 & 1.000 \\
\(\rho_{ij}=\rho_{0}\) & \(t_{5}\) & 0.05 & 3 & 0.290 & 0.081 & 0.049 &
0.034 & 0.028 & 0.294 \\
\(\rho_{ij}=\rho_{0}\) & \(t_{5}\) & 0.05 & 5 & 0.535 & 0.066 & 0.030 &
0.019 & 0.011 & 0.540 \\
\(\boldsymbol{R}=\boldsymbol{I}_{p}\) & scale & 0.05 & 3 & 0.331 & 0.054
& 0.040 & 0.038 & 0.038 & 0.337 \\
\(\boldsymbol{R}=\boldsymbol{I}_{p}\) & scale & 0.05 & 5 & 0.669 & 0.043
& 0.035 & 0.034 & 0.037 & 0.672 \\
\(\boldsymbol{R}=\boldsymbol{I}_{p}\) & scale & 0.10 & 3 & 0.434 & 0.068
& 0.039 & 0.035 & 0.029 & 0.439 \\
\(\boldsymbol{R}=\boldsymbol{I}_{p}\) & scale & 0.10 & 5 & 0.746 & 0.036
& 0.027 & 0.027 & 0.033 & 0.759 \\
\(\boldsymbol{R}=\boldsymbol{I}_{p}\) & shift & 0.05 & 3 & 0.994 & 0.752
& 0.161 & 0.045 & 0.027 & 0.994 \\
\(\boldsymbol{R}=\boldsymbol{I}_{p}\) & shift & 0.05 & 5 & 0.999 & 0.393
& 0.038 & 0.028 & 0.029 & 0.999 \\
\(\boldsymbol{R}=\boldsymbol{I}_{p}\) & shift & 0.10 & 3 & 1.000 & 0.994
& 0.600 & 0.115 & 0.035 & 1.000 \\
\(\boldsymbol{R}=\boldsymbol{I}_{p}\) & shift & 0.10 & 5 & 1.000 & 0.950
& 0.045 & 0.026 & 0.035 & 1.000 \\
\(\boldsymbol{R}=\boldsymbol{I}_{p}\) & \(t_{5}\) & 0.05 & 3 & 0.288 &
0.074 & 0.041 & 0.030 & 0.027 & 0.293 \\
\(\boldsymbol{R}=\boldsymbol{I}_{p}\) & \(t_{5}\) & 0.05 & 5 & 0.510 &
0.068 & 0.036 & 0.024 & 0.016 & 0.523 \\
\end{longtable}

\subsection{Power under contiguous alternatives}\label{sec-SimContig}

The last study verifies numerically Theorem \ref{thm-contig}. The
observations are generated from the correlation matrix
\(\boldsymbol{R}_{n}\) of (\ref{eq-H1n}), with \(\boldsymbol{\Delta}\)
having ones on the first subdiagonal and zeros elsewhere, the same
direction of departure from the null hypothesis as in the study of the
power, so that the alternative hypotheses approach the null one at the
rate \(n^{-1/2}\). The asymptotic power is then the probability that a
\(\chi_{r}^{2}(\delta_{\beta}(\boldsymbol{\Delta}))\) random variable
exceeds \(\chi_{r,\alpha}^{2}\), with
\(\delta_{\beta}(\boldsymbol{\Delta})\) given by (\ref{eq-ARE}), and it
does not depend on the sample size. Table~\ref{tbl-contig} reports it
against the simulated power for \(p=3\) and \(B=500\) replications, for
the test of a specified correlation matrix, where \(r=3\) and
\(\mathrm{trace}((\boldsymbol{R}_{0}^{-1}\boldsymbol{\Delta}^{\ast})^{2})=9.877\),
and for the test of complete uncorrelatedness, where \(r=3\),
\(\boldsymbol{\Delta}^{\ast}=\boldsymbol{\Delta}\) and the trace is
\(4\), and for the test of an equicorrelation with an unknown common
coefficient, where \(r=2\) and the trace is \(3.879\).

The simulated power agrees with the asymptotic one and the agreement
improves with the sample size. For the specified correlation matrix and
\(\beta=0\) the simulated power exceeds the asymptotic one by \(0.054\)
with \(n=500\) and by \(0.024\) with \(n=2000\), and the ordering in
\(\beta\) predicted by Corollary \ref{cor-are} is reproduced in every
scenario. Inverting the power function at each value of \(\beta\) and
dividing the resulting noncentrality parameter by the one obtained for
\(\beta=0\) estimates \(\mathrm{ARE}(\beta,p)\) without using
(\ref{eq-ARE}). For the specified correlation matrix with \(n=500\)
these ratios are \(0.991\), \(0.887\), \(0.839\) and \(0.671\), and for
complete uncorrelatedness with \(n=2000\) they are \(0.966\), \(0.920\),
\(0.828\) and \(0.674\), against the values \(0.971\), \(0.906\),
\(0.826\) and \(0.662\) of Table~\ref{tbl-are}. Each ratio carries a
considerable simulation error, because it inverts the power function at
two estimated proportions, and with \(B=500\) replications every
confidence interval at the \(95\%\) level obtained in this way contains
the corresponding value of \(\mathrm{ARE}(\beta,p)\).

The convergence is slowest for the equicorrelation with an unknown
common coefficient, which is the hypothesis with one more nuisance
parameter to be estimated. Its noncentrality parameter can be estimated
more precisely from the mean of the test-statistic, which is
\(r+\delta_{\beta}(\boldsymbol{\Delta})\), than from the proportion of
rejections, and with \(n=8000\) the two estimates of
\(\mathrm{trace}((\boldsymbol{R}_{0}^{-1}\boldsymbol{\Delta}^{\ast})^{2})\)
so obtained are \(4.09\) and \(4.26\), the latter with a standard error
of \(0.31\), against the value \(3.879\) given by (\ref{eq-adjust}).

\begin{longtable}[]{@{}
  >{\raggedright\arraybackslash}p{(\linewidth - 14\tabcolsep) * \real{0.2391}}
  >{\centering\arraybackslash}p{(\linewidth - 14\tabcolsep) * \real{0.0870}}
  >{\raggedright\arraybackslash}p{(\linewidth - 14\tabcolsep) * \real{0.2391}}
  >{\centering\arraybackslash}p{(\linewidth - 14\tabcolsep) * \real{0.0870}}
  >{\centering\arraybackslash}p{(\linewidth - 14\tabcolsep) * \real{0.0870}}
  >{\centering\arraybackslash}p{(\linewidth - 14\tabcolsep) * \real{0.0870}}
  >{\centering\arraybackslash}p{(\linewidth - 14\tabcolsep) * \real{0.0870}}
  >{\centering\arraybackslash}p{(\linewidth - 14\tabcolsep) * \real{0.0870}}@{}}
\caption{Simulated power of the Rao's \(\beta\)-score tests under the
contiguous alternatives (\ref{eq-H1n}) at the nominal level \(0.05\),
with \(p=3\) and \(B=500\) replications, and asymptotic power given by
Theorem \ref{thm-contig}. The matrix \(\boldsymbol{\Delta}\) has ones on
the first subdiagonal and zeros elsewhere. The asymptotic power does not
depend on the sample size.}\label{tbl-contig}\tabularnewline
\toprule\noalign{}
\begin{minipage}[b]{\linewidth}\raggedright
hypothesis
\end{minipage} & \begin{minipage}[b]{\linewidth}\centering
\(n\)
\end{minipage} & \begin{minipage}[b]{\linewidth}\raggedright
power
\end{minipage} & \begin{minipage}[b]{\linewidth}\centering
\(\beta=0\)
\end{minipage} & \begin{minipage}[b]{\linewidth}\centering
\(0.1\)
\end{minipage} & \begin{minipage}[b]{\linewidth}\centering
\(0.2\)
\end{minipage} & \begin{minipage}[b]{\linewidth}\centering
\(0.3\)
\end{minipage} & \begin{minipage}[b]{\linewidth}\centering
\(0.5\)
\end{minipage} \\
\midrule\noalign{}
\endfirsthead
\toprule\noalign{}
\begin{minipage}[b]{\linewidth}\raggedright
hypothesis
\end{minipage} & \begin{minipage}[b]{\linewidth}\centering
\(n\)
\end{minipage} & \begin{minipage}[b]{\linewidth}\raggedright
power
\end{minipage} & \begin{minipage}[b]{\linewidth}\centering
\(\beta=0\)
\end{minipage} & \begin{minipage}[b]{\linewidth}\centering
\(0.1\)
\end{minipage} & \begin{minipage}[b]{\linewidth}\centering
\(0.2\)
\end{minipage} & \begin{minipage}[b]{\linewidth}\centering
\(0.3\)
\end{minipage} & \begin{minipage}[b]{\linewidth}\centering
\(0.5\)
\end{minipage} \\
\midrule\noalign{}
\endhead
\bottomrule\noalign{}
\endlastfoot
\(\boldsymbol{R}=\boldsymbol{R}_{0}\) & 500 & simulated & 0.490 & 0.486
& 0.440 & 0.418 & 0.340 \\
& & asymptotic & 0.436 & 0.424 & 0.398 & 0.365 & 0.297 \\
\(\boldsymbol{R}=\boldsymbol{R}_{0}\) & 2000 & simulated & 0.460 & 0.434
& 0.400 & 0.354 & 0.284 \\
& & asymptotic & 0.436 & 0.424 & 0.398 & 0.365 & 0.297 \\
\(\boldsymbol{R}=\boldsymbol{I}_{p}\) & 500 & simulated & 0.180 & 0.180
& 0.164 & 0.144 & 0.128 \\
& & asymptotic & 0.192 & 0.188 & 0.177 & 0.165 & 0.140 \\
\(\boldsymbol{R}=\boldsymbol{I}_{p}\) & 2000 & simulated & 0.206 & 0.200
& 0.192 & 0.176 & 0.150 \\
& & asymptotic & 0.192 & 0.188 & 0.177 & 0.165 & 0.140 \\
\(\rho_{ij}=\rho_{12}\) & 500 & simulated & 0.242 & 0.226 & 0.212 &
0.206 & 0.188 \\
& & asymptotic & 0.220 & 0.214 & 0.202 & 0.188 & 0.158 \\
\(\rho_{ij}=\rho_{12}\) & 2000 & simulated & 0.256 & 0.266 & 0.252 &
0.250 & 0.210 \\
& & asymptotic & 0.220 & 0.214 & 0.202 & 0.188 & 0.158 \\
\end{longtable}

\section{A real data example}\label{sec-Appl}

The heptathlon data of the 1988 Seoul Olympic Games, available in the
\texttt{HSAUR2} package of \texttt{R}, record the results of \(n=25\)
athletes in the \(p=7\) events of the competition. The three running
events are measured in seconds and the remaining four in metres or
points, so that the former are replaced by their negatives in order that
a large value is a good result in every event. The hypothesis of
interest is the equicorrelation of Section \ref{sec-test4}, under which
a single factor of general athletic ability loads equally on the seven
events, with the common correlation coefficient \(\rho_{12}\) left
unspecified. The degrees of freedom of the asymptotic distribution under
the null hypothesis are \(\frac{(p+1)(p-2)}{2}=20\) and
\(\chi_{20,0.95}^{2}=31.41\).

The off-diagonal elements of the sample Pearson correlation matrix have
mean \(0.518\) and standard deviation \(0.301\), while the corresponding
matrix with entries \(\sin(\frac{\pi}{2}\tau_{ij})\), built from the
sample Kendall correlation coefficients and used as a starting value in
Remark~\ref{rem-rootpractice}, has mean \(0.521\) and standard deviation
\(0.247\). The two means agree and the two standard deviations do not,
which already suggests that a small number of observations is
responsible for part of the dispersion of the sample correlation
coefficients. The squared Mahalanobis distances computed from the
coordinatewise medians, the median absolute deviations and the latter
matrix isolate a single athlete, the last one in the final
classification, whose value is \(60.45\) against \(11.38\) for the
second largest, and the correlation coefficients of the javelin throw
with the remaining events are the ones that this observation distorts
the most, from values between \(-0.02\) and \(0.33\) for the Pearson
coefficients to values between \(0.19\) and \(0.38\) for the robust
ones.

\begin{longtable}[]{@{}
  >{\raggedright\arraybackslash}p{(\linewidth - 10\tabcolsep) * \real{0.1667}}
  >{\centering\arraybackslash}p{(\linewidth - 10\tabcolsep) * \real{0.1667}}
  >{\centering\arraybackslash}p{(\linewidth - 10\tabcolsep) * \real{0.1667}}
  >{\centering\arraybackslash}p{(\linewidth - 10\tabcolsep) * \real{0.1667}}
  >{\centering\arraybackslash}p{(\linewidth - 10\tabcolsep) * \real{0.1667}}
  >{\centering\arraybackslash}p{(\linewidth - 10\tabcolsep) * \real{0.1667}}@{}}
\caption{Rao's \(\beta\)-score test-statistics (\ref{eq-RaoR0g}) for the
equicorrelation of the heptathlon data, with all the athletes and with
the most outlying one deleted, and the relative change of the
test-statistic produced by the deletion.}\label{tbl-hept}\tabularnewline
\toprule\noalign{}
\begin{minipage}[b]{\linewidth}\raggedright
\(\beta\)
\end{minipage} & \begin{minipage}[b]{\linewidth}\centering
\(\widetilde{R}_{\beta,n}\), all
\end{minipage} & \begin{minipage}[b]{\linewidth}\centering
\(p\)-value
\end{minipage} & \begin{minipage}[b]{\linewidth}\centering
\(\widetilde{R}_{\beta,n}\), deleted
\end{minipage} & \begin{minipage}[b]{\linewidth}\centering
\(p\)-value
\end{minipage} & \begin{minipage}[b]{\linewidth}\centering
change
\end{minipage} \\
\midrule\noalign{}
\endfirsthead
\toprule\noalign{}
\begin{minipage}[b]{\linewidth}\raggedright
\(\beta\)
\end{minipage} & \begin{minipage}[b]{\linewidth}\centering
\(\widetilde{R}_{\beta,n}\), all
\end{minipage} & \begin{minipage}[b]{\linewidth}\centering
\(p\)-value
\end{minipage} & \begin{minipage}[b]{\linewidth}\centering
\(\widetilde{R}_{\beta,n}\), deleted
\end{minipage} & \begin{minipage}[b]{\linewidth}\centering
\(p\)-value
\end{minipage} & \begin{minipage}[b]{\linewidth}\centering
change
\end{minipage} \\
\midrule\noalign{}
\endhead
\bottomrule\noalign{}
\endlastfoot
\(0\) & 82.40 & \(1.5\times10^{-9}\) & 47.51 & \(5.0\times10^{-4}\) &
\(-42.3\%\) \\
\(0.05\) & 49.36 & \(2.7\times10^{-4}\) & 49.22 & \(2.9\times10^{-4}\) &
\(-0.3\%\) \\
\(0.10\) & 47.72 & \(4.7\times10^{-4}\) & 49.76 & \(2.4\times10^{-4}\) &
\(4.3\%\) \\
\(0.20\) & 45.66 & \(9.0\times10^{-4}\) & 47.67 & \(4.7\times10^{-4}\) &
\(4.4\%\) \\
\(0.30\) & 39.86 & \(5.2\times10^{-3}\) & 41.14 & \(3.6\times10^{-3}\) &
\(3.2\%\) \\
\(0.40\) & 27.37 & 0.125 & 27.07 & 0.133 & \(-1.1\%\) \\
\(0.50\) & 11.40 & 0.935 & 10.78 & 0.952 & \(-5.5\%\) \\
\end{longtable}

\begin{figure}[!htb]

\centering{

\pandocbounded{\includegraphics[keepaspectratio]{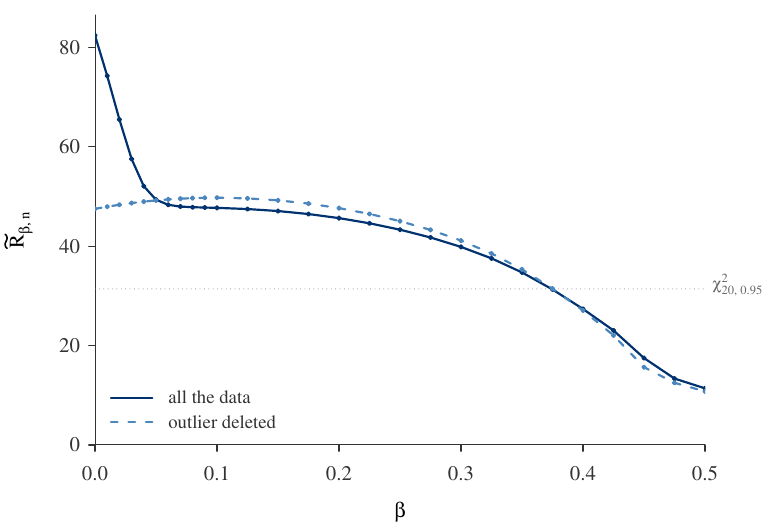}}

}

\caption{\label{fig-hept}Rao's \(\beta\)-score test-statistics for the
equicorrelation of the heptathlon data, with all the athletes and with
the most outlying one deleted. The dotted line is
\(\chi_{20,0.95}^{2}\).}

\end{figure}

The test-statistics are reported in Table~\ref{tbl-hept}, for a
selection of values of the tuning parameter, and displayed in
Figure~\ref{fig-hept} over a grid of thirty-one values. The
equicorrelation is rejected at the five per cent level by every
test-statistic with \(\beta\leq0.35\), with and without the outlying
athlete, so that the deletion does not reverse the decision. What the
deletion does reverse is the strength of the evidence reported by the
classical test-statistic, which falls from \(82.40\) to \(47.51\), a
reduction of \(42.3\%\), and whose \(p\)-value increases by five orders
of magnitude. A single athlete out of twenty five is therefore
responsible for most of the apparent evidence against the
equicorrelation when \(\beta=0\). From \(\beta=0.05\) onwards the two
test-statistics differ by less than five per cent of their value, the
two curves of Figure~\ref{fig-hept} are indistinguishable at the scale
of the figure, and the conclusion no longer depends on whether the
athlete is kept or removed. The finer grid of Figure~\ref{fig-hept}
shows that the descent is gradual rather than abrupt, the test-statistic
taking the values \(74.27\), \(65.51\), \(57.55\) and \(52.07\) at
\(\beta=0.01,0.02,0.03\) and \(0.04\), so that the weight discounts the
outlying athlete continuously. Both test-statistics cross
\(\chi_{20,0.95}^{2}\) at \(\beta=0.375\), where they take the values
\(31.30\) and \(31.46\). The price of this stability is small, since
\(\mathrm{ARE}(0.1,7)=0.955\) by Corollary~\ref{cor-are}, that is, less
than five per cent of asymptotic efficiency under the multivariate
normal distribution.

The behaviour of the test-statistics for the larger values of \(\beta\)
illustrates Remark~\ref{rem-rootpractice}. For \(\beta\leq0.5\) the
descent algorithm reaches the same solution of the restricted estimating
equations from the classical and from the robust starting values, the
two values of (\ref{eq-Hn}) differing by less than \(10^{-8}\), and the
estimates \(\widetilde{\rho}_{12,\beta}\) remain between \(0.52\) and
\(0.56\). The diagnostic fails first at \(\beta=0.525\) with the
outlying athlete deleted and at \(\beta=0.55\) with all the data, the
two values of (\ref{eq-Hn}) differing then by \(1.4\times10^{2}\) and by
\(5.4\times10^{4}\), and by more than \(10^{9}\) for the larger values
of the tuning parameter, while the estimate of the common correlation
coefficient drops to \(-0.16\), so that the solutions no longer remain
in a compact set away from the degenerate solutions of
Remark~\ref{rem-unbounded} and the root is not identified. The
diagnostic of Remark~\ref{rem-rootpractice} thus detects, in a sample as
small as this one, the failure that Theorem~\ref{thm-root} only excludes
asymptotically, and the corresponding test-statistics are discarded.
This is the reason why Figure~\ref{fig-hept} is drawn for
\(\beta\leq0.5\).

\section{Discussion}\label{sec-Discuss}

The classical Rao's score test-statistics derived in this paper are
closely related to the classical tests on correlation matrices. For a
specified correlation matrix, the test-statistic of Jennrich (1970) is a
quadratic form in the difference between the sample and the specified
correlation matrices, whose second term corrects for the estimation of
the variances, as the restricted estimators do here. Since it is
algebraically equal to a Wald test-statistic (Neudecker and Satorra,
1996), it is asymptotically equivalent under the null hypothesis to
(\ref{eq-RaoR0}) with \(\beta=0\). For an unknown common correlation,
(\ref{eq-RaoR0g}) with \(\beta=0\) is asymptotically equivalent to the
test-statistic of Lawley (1963), see Remark~\ref{rem-lawley}. The Rao's
\(\beta\)-score test-statistics with \(\beta>0\) have no classical
counterpart.

By Theorem~\ref{thm-contig}, the noncentrality parameter of the
classical Rao's score test for a specified correlation matrix is
\(\frac{1}{2}\mathrm{trace}((\boldsymbol{R}_{0}^{-1}\boldsymbol{\Delta}^{\ast})^{2})\),
a quadratic form with the same structure as the test-statistic of
Jennrich (1970), and the price paid for robustness is the factor
\(\mathrm{ARE}(\beta,p)\), which does not depend on the hypothesis
tested. For a fixed \(\beta\) this efficiency decreases with the
dimension. The reason is that, under the model,
\((\boldsymbol{X}_{i}-\boldsymbol{\mu})^{\top}\boldsymbol{\Sigma}^{-1}(\boldsymbol{X}_{i}-\boldsymbol{\mu})\sim\chi_{p}^{2}\),
so the mean weight of an observation is \[
E_{\boldsymbol{\theta}}\left[W_{i,\beta}(\boldsymbol{\theta})\right]=(1+\beta)^{-p/2},
\] and smaller values of \(\beta\) should be used as \(p\) increases,
see Table~\ref{tbl-are}. On the other hand, by Theorem~\ref{thm-IF2} and
Theorem~\ref{thm-PIF}, any \(\beta>0\) gives bounded influence
functions, but the supremum in (\ref{eq-IF2ind}) is of order
\(\beta^{-2}\) as \(\beta\) tends to zero, so very small values of
\(\beta\) give little protection, and the choice of \(\beta\) is a
compromise between both properties. By Remark~\ref{rem-elliptical} the
tests are also consistent under elliptical distributions, and a single
scalar separates their asymptotic distribution from the chi-square one,
so a distribution-free version over the elliptical family is available
at the price of estimating that scalar, an alternative being the tests
based on the distribution-free M-estimator of scatter of Tyler (1987).
The behaviour of the Rao's score tests when \(p\) increases with \(n\)
requires a different asymptotic framework and is not addressed here.

\section*{References}\label{references}
\addcontentsline{toc}{section}{References}

Aitchison, J. and Silvey, S. D. (1958). Maximum-likelihood estimation of
parameters subject to restraints. \emph{Annals of Mathematical
Statistics}, \textbf{29}, 813--828.

Anderson, T. W. (2003). \emph{An Introduction to Multivariate
Statistical Analysis}, 3rd ed.~Hoboken, NJ: John Wiley \& Sons.

Bartlett, M. S. (1954). A note on the multiplying factors for various
\(\chi^{2}\) approximations. \emph{Journal of the Royal Statistical
Society, Series B}, \textbf{16}, 296--298.

Basu, A., Chakraborty, S., Ghosh, A. and Pardo, L. (2022a). Robust
density power divergence based tests in multivariate analysis, a
comparative overview of different approaches. \emph{Journal of
Multivariate Analysis}, \textbf{188}, 104846.

Basu, A., Ghosh, A., Martín, N. and Pardo, L. (2022b). A robust
generalization of the Rao test. \emph{Journal of Business \& Economic
Statistics}, \textbf{40}, 868--879.

Basu, A., Harris, I. R., Hjort, N. L. and Jones, M. C. (1998). Robust
and efficient estimation by minimising a density power divergence.
\emph{Biometrika}, \textbf{85}, 549--559.

Basu, A., Shioya, H. and Park, C. (2011). \emph{Statistical Inference,
the Minimum Distance Approach}. Boca Raton, FL: Chapman \& Hall/CRC.

Browne, M. W. (1974). Generalized least squares estimators in the
analysis of covariance structures. \emph{South African Statistical
Journal}, \textbf{8}, 1--24.

Ferrari, D. and Yang, Y. (2010). Maximum Lq-likelihood estimation.
\emph{Annals of Statistics}, \textbf{38}, 753--783.

Ghosh, A., Mandal, A., Martín, N. and Pardo, L. (2016). Influence
analysis of robust Wald-type tests. \emph{Journal of Multivariate
Analysis}, \textbf{147}, 102--126.

Hampel, F. R., Ronchetti, E. M., Rousseeuw, P. J. and Stahel, W. A.
(1986). \emph{Robust Statistics, the Approach Based on Influence
Functions}. New York: John Wiley \& Sons.

Henderson, H. V. and Searle, S. R. (1979). Vec and vech operators for
matrices, with some uses in Jacobians and multivariate statistics.
\emph{Canadian Journal of Statistics}, \textbf{7}, 65--81.

Heritier, S. and Ronchetti, E. (1994). Robust bounded-influence tests in
general parametric models. \emph{Journal of the American Statistical
Association}, \textbf{89}, 897--904.

Jennrich, R. I. (1970). An asymptotic \(\chi^{2}\) test for the equality
of two correlation matrices. \emph{Journal of the American Statistical
Association}, \textbf{65}, 904--912.

Kallenberg, W. C. M., Ledwina, T. and Rafaj\l{}owicz, E. (1997). Testing
bivariate independence and normality. \emph{Sankhy=\{a\}, Series A},
\textbf{59}, 42--59.

Lawley, D. N. (1963). On testing a set of correlation coefficients for
equality. \emph{Annals of Mathematical Statistics}, \textbf{34},
149--151.

Lehmann, E. L. (1999). \emph{Elements of Large-Sample Theory}. New York:
Springer.

Leung, D. and Drton, M. (2018). Testing independence in high dimensions
with sums of rank correlations. \emph{Annals of Statistics},
\textbf{46}, 280--307.

Magnus, J. R. and Neudecker, H. (1980). The elimination matrix, some
lemmas and applications. \emph{SIAM Journal on Algebraic and Discrete
Methods}, \textbf{1}, 422--449.

Martín, N. (2020). Rao's score tests on correlation matrices.
arXiv:2012.14238, first version.

McCulloch, C. E. (1982). Symmetric matrix derivatives with applications.
\emph{Journal of the American Statistical Association}, \textbf{77},
679--682.

Muirhead, R. J. and Waternaux, C. M. (1980). Asymptotic distributions in
canonical correlation analysis and other multivariate procedures for
nonnormal populations. \emph{Biometrika}, \textbf{67}, 31--43.

Neudecker, H. and Satorra, A. (1996). The algebraic equality of two
asymptotic tests for the hypothesis that a normal distribution has a
specified correlation matrix. \emph{Statistics \& Probability Letters},
\textbf{30}, 99--103.

Rao, C. R. (1948). Large sample tests of statistical hypotheses
concerning several parameters with applications to problems of
estimation. \emph{Mathematical Proceedings of the Cambridge
Philosophical Society}, \textbf{44}, 50--57.

Silvey, S. D. (1959). The Lagrangian multiplier test. \emph{Annals of
Mathematical Statistics}, \textbf{30}, 389--407.

Tyler, D. E. (1983). Robustness and efficiency properties of scatter
matrices. \emph{Biometrika}, \textbf{70}, 411--420.

Tyler, D. E. (1987). A distribution-free M-estimator of multivariate
scatter. \emph{Annals of Statistics}, \textbf{15}, 234--251.

van der Vaart, A. W. (1998). \emph{Asymptotic Statistics}. Cambridge:
Cambridge University Press.

Wald, A. (1943). Tests of statistical hypotheses concerning several
parameters when the number of observations is large. \emph{Transactions
of the American Mathematical Society}, \textbf{54}, 426--482.

\appendix

\section{Proofs}\label{sec-Proofs}

\subsection{\texorpdfstring{Proof of Proposition
\ref{prp-Props}}{Proof of Proposition }}\label{sec-ProofProps}

The detailed proof of the expressions of the score functions~was given
in McCulloch (1982), and in particular the expression of the second
partition one needs careful derivations,

\[
\begin{aligned}
\boldsymbol{s}_{\boldsymbol{\Sigma}}(\boldsymbol{X})  &  =-\frac{1}{2}
\frac{\partial}{\partial\boldsymbol{\Sigma}}\left(  \log\left\vert
\boldsymbol{\Sigma}\right\vert +(\boldsymbol{X}-\boldsymbol{\mu}
)^{\top}\boldsymbol{\Sigma}^{-1}(\boldsymbol{X}-\boldsymbol{\mu})\right) \\
&  =-\frac{1}{2}\boldsymbol{\Sigma}^{-1}+\frac{\partial}{\partial
\boldsymbol{\Sigma}}\left[  (\boldsymbol{X}-\boldsymbol{\mu})^{\top}
\boldsymbol{\Sigma}^{-1}(\boldsymbol{X}-\boldsymbol{\mu})\right] \\
&  =-\frac{1}{2}\boldsymbol{\Sigma}^{-1}-\boldsymbol{\Sigma}^{-1}
(\boldsymbol{X}-\boldsymbol{\mu})(\boldsymbol{X}-\boldsymbol{\mu}
)^{\top}\boldsymbol{\Sigma}^{-1}
\end{aligned}
\]

and hence,

\[
\begin{aligned}
\boldsymbol{s}_{\mathrm{vech}(\boldsymbol{\Sigma})}(\boldsymbol{X})  &
=-\frac{1}{2}\mathrm{vech}\left(  \boldsymbol{\Sigma}^{-1}-\boldsymbol{\Sigma
}^{-1}(\boldsymbol{X}-\boldsymbol{\mu})(\boldsymbol{X}-\boldsymbol{\mu}
)^{\top}\boldsymbol{\Sigma}^{-1}\right) \\
&  =-\frac{1}{2}\boldsymbol{G}_{p}^{\top}\mathrm{vec}(\boldsymbol{\Sigma}
^{-1})+\frac{1}{2}\boldsymbol{G}_{p}^{\top}\mathrm{vec}\left(  \boldsymbol{\Sigma
}^{-1}(\boldsymbol{X}-\boldsymbol{\mu})(\boldsymbol{X}-\boldsymbol{\mu}
)^{\top}\boldsymbol{\Sigma}^{-1}\right) \\
&  =-\frac{1}{2}\boldsymbol{G}_{p}^{\top}\mathrm{vec}(\boldsymbol{\Sigma}
^{-1})+\frac{1}{2}\boldsymbol{G}_{p}^{\top}\left(  \boldsymbol{\Sigma}
^{-1}\otimes\boldsymbol{\Sigma}^{-1}\right)  \mathrm{vec}\left(
(\boldsymbol{X}-\boldsymbol{\mu})(\boldsymbol{X}-\boldsymbol{\mu})^{\top}\right) \\
&  =-\frac{1}{2}\boldsymbol{G}_{p}^{\top}\mathrm{vec}(\boldsymbol{\Sigma}
^{-1})+\frac{1}{2}\boldsymbol{G}_{p}^{\top}\left(  \boldsymbol{\Sigma}
^{-1}\otimes\boldsymbol{\Sigma}^{-1}\right)  \left(  (\boldsymbol{X}
-\boldsymbol{\mu})\otimes(\boldsymbol{X}-\boldsymbol{\mu})\right) \\
&  =-\frac{1}{2}\boldsymbol{G}_{p}^{\top}\mathrm{vec}(\boldsymbol{\Sigma}
^{-1})+\frac{1}{2}\boldsymbol{G}_{p}^{\top}\left[  \left(  \boldsymbol{\Sigma
}^{-1}(\boldsymbol{X}-\boldsymbol{\mu})\right)  \otimes\left(  (\boldsymbol{X}
-\boldsymbol{\mu})^{\top}\boldsymbol{\Sigma}^{-1}\right)  \right]  .
\end{aligned}
\]

\subsection{\texorpdfstring{Proof of Theorem
\ref{thm-Th-J-BS1}}{Proof of Theorem }}\label{sec-Proof-Th-J-BS1}

Let us consider

\[
\begin{aligned}
f_{\theta}^{\beta+1}(\boldsymbol{x})  &  =\frac{1}{(2\pi)^{\frac{(\beta
+1)p}{2}}\left\vert \boldsymbol{\Sigma}\right\vert ^{\frac{\beta+1}{2}}}
\exp\left\{  -\frac{1}{2}(\boldsymbol{x}-\boldsymbol{\mu})^{\top}\left(
\tfrac{1}{\beta+1}\boldsymbol{\Sigma}\right)  ^{-1}(\boldsymbol{x}
-\boldsymbol{\mu})\right\} \\
&  =\frac{(2\pi)^{\frac{p}{2}}\left\vert \tfrac{1}{\beta+1}\boldsymbol{\Sigma
}\right\vert ^{\frac{1}{2}}}{(2\pi)^{\frac{p}{2}}\left\vert \tfrac{1}{\beta
+1}\boldsymbol{\Sigma}\right\vert ^{\frac{1}{2}}}\frac{1}{(2\pi)^{\frac
{(\beta+1)p}{2}}\left\vert \boldsymbol{\Sigma}\right\vert ^{\frac{\beta+1}{2}
}}\exp\left\{  -\frac{1}{2}(\boldsymbol{x}-\boldsymbol{\mu})^{\top}\left(
\tfrac{1}{\beta+1}\boldsymbol{\Sigma}\right)  ^{-1}(\boldsymbol{x}
-\boldsymbol{\mu})\right\} \\
&  =\frac{(2\pi)^{\frac{p}{2}}\left\vert \tfrac{1}{\beta+1}\boldsymbol{\Sigma
}\right\vert ^{\frac{1}{2}}}{(2\pi)^{\frac{(\beta+1)p}{2}}\left\vert
\boldsymbol{\Sigma}\right\vert ^{\frac{\beta+1}{2}}}\frac{1}{(2\pi)^{\frac
{p}{2}}\left\vert \tfrac{1}{\beta+1}\boldsymbol{\Sigma}\right\vert ^{\frac
{1}{2}}}\exp\left\{  -\frac{1}{2}(\boldsymbol{x}-\boldsymbol{\mu})^{\top}\left(
\tfrac{1}{\beta+1}\boldsymbol{\Sigma}\right)  ^{-1}(\boldsymbol{x}
-\boldsymbol{\mu})\right\} \\
&  =\frac{(\beta+1)^{-\frac{p}{2}}}{(2\pi)^{\frac{\beta p}{2}}\left\vert
\boldsymbol{\Sigma}\right\vert ^{\frac{\beta}{2}}}f_{\theta^{\ast}
}(\boldsymbol{x}),
\end{aligned}
\]

where the parameter vector is given by
\(\boldsymbol{\theta}^{\ast }=(\boldsymbol{\mu}^{\ast T},\mathrm{vech}^{\top}(\boldsymbol{\Sigma}^{\ast }))^{\top}\),
\(\boldsymbol{\mu}^{\ast}=\boldsymbol{\mu}\) and
\(\boldsymbol{\Sigma }^{\ast}=\tfrac{1}{\beta+1}\boldsymbol{\Sigma}\).
Notice that if we call

\begin{equation}\protect\phantomsection\label{eq-paramEta}{
\boldsymbol{\vartheta}=(\boldsymbol{\mu}^{\top},(\beta+1)\mathrm{vech}
^{\top}(\boldsymbol{\Sigma}))^{\top}
}\end{equation}

then

\[
\begin{aligned}
f_{\boldsymbol{\vartheta}}^{\beta+1}(\boldsymbol{x})  &  =\frac{(\beta
+1)^{-\frac{p}{2}}}{(2\pi)^{\frac{\beta p}{2}}\left\vert (\beta
+1)\boldsymbol{\Sigma}\right\vert ^{\frac{\beta}{2}}}f_{\boldsymbol{\theta}
}(\boldsymbol{x}) \\
&  =\frac{(\beta+1)^{-\frac{p}{2}(\beta+1)}}{(2\pi)^{\frac{\beta p}{2}
}\left\vert \boldsymbol{\Sigma}\right\vert ^{\frac{\beta}{2}}}
f_{\boldsymbol{\theta}}(\boldsymbol{x}),
\end{aligned}
\]

\[
\begin{aligned}
\boldsymbol{J}_{\beta}^{\ast}(\boldsymbol{\vartheta})  &
=E_{\boldsymbol{\vartheta}}[\boldsymbol{\varsigma}_{\boldsymbol{\vartheta}
}(\boldsymbol{X})\boldsymbol{\varsigma}_{\boldsymbol{\vartheta}}
^{\top}(\boldsymbol{X})f_{\boldsymbol{\vartheta}}^{\beta}] \\
&  =\int_{\mathbb{R}^{p}}\boldsymbol{\varsigma}_{\boldsymbol{\vartheta}}(\boldsymbol{x}
)\boldsymbol{\varsigma}_{\boldsymbol{\vartheta}}^{\top}(\boldsymbol{x}
)f_{\boldsymbol{\vartheta}}^{\beta+1}(\boldsymbol{x})d\boldsymbol{x} \\
&  =\frac{(\beta+1)^{-\frac{p}{2}(\beta+1)}}{(2\pi)^{\frac{\beta p}{2}
}\left\vert \boldsymbol{\Sigma}\right\vert ^{\frac{\beta}{2}}}\int_{\mathbb{R}^{p}}\boldsymbol{\varsigma}_{\boldsymbol{\vartheta}}(\boldsymbol{x}
)\boldsymbol{\varsigma}_{\boldsymbol{\vartheta}}^{\top}(\boldsymbol{x}
)f_{\boldsymbol{\theta}}(\boldsymbol{x})d\boldsymbol{x} \\
&  =\frac{(\beta+1)^{-\frac{p}{2}(\beta+1)}}{(2\pi)^{\frac{\beta p}{2}
}\left\vert \boldsymbol{\Sigma}\right\vert ^{\frac{\beta}{2}}}
E_{\boldsymbol{\theta}}[\boldsymbol{\varsigma}_{\boldsymbol{\vartheta}
}(\boldsymbol{X})\boldsymbol{\varsigma}_{\boldsymbol{\vartheta}}
^{\top}(\boldsymbol{X})],
\end{aligned}
\]

where

\[
\begin{aligned}
\boldsymbol{\varsigma}_{\boldsymbol{\vartheta}}(\boldsymbol{x})  &
=\frac{\partial}{\partial\boldsymbol{\vartheta}}\log f_{\boldsymbol{\vartheta
}}(\boldsymbol{x})=\binom{\frac{\partial}{\partial\boldsymbol{\mu}}\log
f_{\boldsymbol{\vartheta}}(\boldsymbol{x})}{\frac{\partial}{\partial
(\beta+1)\mathrm{vech}(\boldsymbol{\Sigma})}\log f_{\boldsymbol{\vartheta}
}(\boldsymbol{x})} \\
&  =\binom{\frac{1}{\beta+1}\frac{\partial}{\partial\boldsymbol{\mu}}\log
f_{\boldsymbol{\theta}}(\boldsymbol{x})}{\frac{1}{(\beta+1)^{2}}\left(
-\frac{\beta}{2}\frac{\partial}{\partial\mathrm{vech}(\boldsymbol{\Sigma}
)}\log\left\vert \boldsymbol{\Sigma}\right\vert +\frac{\partial}
{\partial\mathrm{vech}(\boldsymbol{\Sigma})}\log f_{\boldsymbol{\theta}
}(\boldsymbol{x})\right)  },
\end{aligned}
\]

\[
\begin{aligned}
\frac{\partial}{\partial\boldsymbol{\mu}}\log f_{\boldsymbol{\vartheta}
}(\boldsymbol{x})  &  =\frac{\partial}{\partial\boldsymbol{\mu}}\left(
-\frac{p}{2}\log(2\pi)-\frac{1}{2}\log\left\vert (\beta+1)\boldsymbol{\Sigma
}\right\vert -\frac{1}{2(\beta+1)}(\boldsymbol{x}-\boldsymbol{\mu}
)^{\top}\boldsymbol{\Sigma}^{-1}(\boldsymbol{x}-\boldsymbol{\mu})\right) \\
&  =-\frac{\partial}{\partial\boldsymbol{\mu}}\left(  \frac{1}{2(\beta
+1)}(\boldsymbol{x}-\boldsymbol{\mu})^{\top}\boldsymbol{\Sigma}^{-1}
(\boldsymbol{x}-\boldsymbol{\mu})\right) \\
&  =-\frac{1}{\beta+1}\frac{\partial}{\partial\boldsymbol{\mu}}\left(
\frac{1}{2}(\boldsymbol{x}-\boldsymbol{\mu})^{\top}\boldsymbol{\Sigma}
^{-1}(\boldsymbol{x}-\boldsymbol{\mu})\right) \\
&  =\frac{1}{\beta+1}\frac{\partial}{\partial\boldsymbol{\mu}}\log
f_{\boldsymbol{\theta}}(\boldsymbol{x}),
\end{aligned}
\]

\[
\begin{aligned}
&  \frac{\partial}{\partial(\beta+1)\mathrm{vech}(\boldsymbol{\Sigma})}\log
f_{\boldsymbol{\vartheta}}(\boldsymbol{x}) \\
&  =\frac{\partial}{\partial(\beta+1)\mathrm{vech}(\boldsymbol{\Sigma}
)}\left(  -\frac{p}{2}\log(2\pi)-\frac{1}{2}\log\left\vert (\beta
+1)\boldsymbol{\Sigma}\right\vert -\frac{1}{2}(\boldsymbol{x}-\boldsymbol{\mu
})^{\top}(\beta+1)^{-1}\boldsymbol{\Sigma}^{-1}(\boldsymbol{x}-\boldsymbol{\mu
})\right) \\
&  =\frac{\partial}{\partial(\beta+1)\mathrm{vech}(\boldsymbol{\Sigma}
)}\left(  -\frac{1}{2}\log\left\vert (\beta+1)\boldsymbol{\Sigma}\right\vert
-\frac{1}{2}(\boldsymbol{x}-\boldsymbol{\mu})^{\top}(\beta+1)^{-1}
\boldsymbol{\Sigma}^{-1}(\boldsymbol{x}-\boldsymbol{\mu})\right) \\
&  =\frac{\partial\mathrm{vech}^{\top}(\boldsymbol{\Sigma})}{\partial
(\beta+1)\mathrm{vech}(\boldsymbol{\Sigma})}\frac{\partial}{\partial
\mathrm{vech}(\boldsymbol{\Sigma})}\left(  -\frac{1}{2}\log\left\vert
\boldsymbol{\Sigma}\right\vert -\frac{1}{2(\beta+1)}(\boldsymbol{x}
-\boldsymbol{\mu})^{\top}\boldsymbol{\Sigma}^{-1}(\boldsymbol{x}-\boldsymbol{\mu
})\right) \\
&  =\frac{1}{\beta+1}\frac{\partial}{\partial\mathrm{vech}(\boldsymbol{\Sigma
})}\left(  -\frac{\beta+1}{2(\beta+1)}\log\left\vert \boldsymbol{\Sigma
}\right\vert -\frac{1}{2(\beta+1)}(\boldsymbol{x}-\boldsymbol{\mu}
)^{\top}\boldsymbol{\Sigma}^{-1}(\boldsymbol{x}-\boldsymbol{\mu})\right) \\
&  =\frac{1}{(\beta+1)^{2}}\frac{\partial}{\partial\mathrm{vech}
(\boldsymbol{\Sigma})}\left(  -\frac{\beta+1}{2}\log\left\vert
\boldsymbol{\Sigma}\right\vert -\frac{1}{2}(\boldsymbol{x}-\boldsymbol{\mu
})^{\top}\boldsymbol{\Sigma}^{-1}(\boldsymbol{x}-\boldsymbol{\mu})\right) \\
&  =\frac{1}{(\beta+1)^{2}}\left[  -\frac{\beta}{2}\frac{\partial}
{\partial\mathrm{vech}(\boldsymbol{\Sigma})}\log\left\vert \boldsymbol{\Sigma
}\right\vert +\frac{\partial}{\partial\mathrm{vech}(\boldsymbol{\Sigma}
)}\left(  -\frac{1}{2}\log\left\vert \boldsymbol{\Sigma}\right\vert -\frac
{1}{2}(\boldsymbol{x}-\boldsymbol{\mu})^{\top}\boldsymbol{\Sigma}^{-1}
(\boldsymbol{x}-\boldsymbol{\mu})\right)  \right] \\
&  =\frac{1}{(\beta+1)^{2}}\left(  -\frac{\beta}{2}\frac{\partial}
{\partial\mathrm{vech}(\boldsymbol{\Sigma})}\log\left\vert \boldsymbol{\Sigma
}\right\vert +\frac{\partial}{\partial\mathrm{vech}(\boldsymbol{\Sigma})}\log
f_{\boldsymbol{\theta}}(\boldsymbol{x})\right)  ,
\end{aligned}
\]

\[
\begin{aligned}
\boldsymbol{\varsigma}_{\boldsymbol{\vartheta}}(\boldsymbol{X})  &
\boldsymbol{=}
\begin{pmatrix}
\boldsymbol{\varsigma}_{\boldsymbol{\mu}}(\boldsymbol{X}) \\
\boldsymbol{\varsigma}_{(\beta+1)\mathrm{vech}(\boldsymbol{\Sigma}
)}(\boldsymbol{X})
\end{pmatrix}
,\qquad\boldsymbol{\varsigma}_{\boldsymbol{\vartheta}}(\boldsymbol{X}
)\boldsymbol{=}\frac{\partial}{\partial\boldsymbol{\vartheta}}\log
f_{\boldsymbol{\vartheta}}(\boldsymbol{X}), \\
\boldsymbol{\varsigma}_{\boldsymbol{\mu}}(\boldsymbol{X})  &  =\frac{\partial
}{\partial\boldsymbol{\mu}}\log f_{\boldsymbol{\vartheta}}(\boldsymbol{x}
),\qquad\boldsymbol{\varsigma}_{(\beta+1)\mathrm{vech}(\boldsymbol{\Sigma}
)}(\boldsymbol{X})=\frac{\partial}{\partial(\beta+1)\mathrm{vech}
(\boldsymbol{\Sigma})}\log f_{\boldsymbol{\vartheta}}(\boldsymbol{X}), \\
\boldsymbol{s}_{\boldsymbol{\theta}}(\boldsymbol{X})  &  \boldsymbol{=}
\begin{pmatrix}
\boldsymbol{s}_{\boldsymbol{\mu}}(\boldsymbol{X}) \\
\boldsymbol{s}_{\mathrm{vech}(\boldsymbol{\Sigma})}(\boldsymbol{X})
\end{pmatrix}
,\qquad\boldsymbol{s}_{\boldsymbol{\theta}}(\boldsymbol{X})\boldsymbol{=}
\frac{\partial}{\partial\boldsymbol{\theta}}\log f_{\boldsymbol{\theta}
}(\boldsymbol{X}), \\
\boldsymbol{s}_{\boldsymbol{\mu}}(\boldsymbol{X})  &  =\frac{\partial
}{\partial\boldsymbol{\mu}}\log f_{\boldsymbol{\theta}}(\boldsymbol{X}
),\qquad\boldsymbol{s}_{\mathrm{vech}(\boldsymbol{\Sigma})}(\boldsymbol{X}
)=\frac{\partial}{\partial\mathrm{vech}(\boldsymbol{\Sigma})}\log
f_{\boldsymbol{\theta}}(\boldsymbol{X}),
\end{aligned}
\]

\[
\begin{aligned}
\boldsymbol{\varsigma}_{\boldsymbol{\vartheta}}(\boldsymbol{X}
)\boldsymbol{\varsigma}_{\boldsymbol{\vartheta}}^{\top}(\boldsymbol{X})=
\begin{pmatrix}
\boldsymbol{\varsigma}_{\boldsymbol{\mu}}(\boldsymbol{X})\boldsymbol{\varsigma
}_{\boldsymbol{\mu}}^{\top}(\boldsymbol{X}) & \boldsymbol{\varsigma
}_{\boldsymbol{\mu}}(\boldsymbol{X})\boldsymbol{\varsigma}_{(\beta
+1)\mathrm{vech}(\boldsymbol{\Sigma})}^{\top}(\boldsymbol{X})\\
\boldsymbol{\varsigma}_{(\beta+1)\mathrm{vech}(\boldsymbol{\Sigma}
)}(\boldsymbol{X})\boldsymbol{\varsigma}_{\boldsymbol{\mu}}^{\top}(\boldsymbol{X}
) & \boldsymbol{\varsigma}_{(\beta+1)\mathrm{vech}(\boldsymbol{\Sigma}
)}(\boldsymbol{X})\boldsymbol{\varsigma}_{(\beta+1)\mathrm{vech}
(\boldsymbol{\Sigma})}^{\top}(\boldsymbol{X})
\end{pmatrix}
,
\end{aligned}
\]

\[
\boldsymbol{\varsigma}_{\boldsymbol{\mu}}(\boldsymbol{X})\boldsymbol{\varsigma
}_{\boldsymbol{\mu}}^{\top}(\boldsymbol{X})=\frac{1}{(\beta+1)^{2}}
\boldsymbol{s}_{\boldsymbol{\mu}}(\boldsymbol{X})\boldsymbol{s}
_{\boldsymbol{\mu}}^{\top}(\boldsymbol{X}),
\]

\[
\begin{aligned}
&  \boldsymbol{\varsigma}_{(\beta+1)\mathrm{vech}(\boldsymbol{\Sigma}
)}(\boldsymbol{X})\boldsymbol{\varsigma}_{\boldsymbol{\mu}}^{\top}(\boldsymbol{X}
) \\
&  =\frac{1}{(\beta+1)^{2}}\frac{1}{\beta+1}\left(  -\frac{\beta}{2}
\frac{\partial}{\partial\mathrm{vech}(\boldsymbol{\Sigma})}\log\left\vert
\boldsymbol{\Sigma}\right\vert +\boldsymbol{s}_{\mathrm{vech}
(\boldsymbol{\Sigma})}(\boldsymbol{X})\right)  \boldsymbol{s}_{\boldsymbol{\mu
}}^{\top}(\boldsymbol{X}) \\
&  =\frac{1}{(\beta+1)^{2}}\left(  -\frac{\beta}{2(\beta+1)}\frac{\partial
}{\partial\mathrm{vech}(\boldsymbol{\Sigma})}\log\left\vert \boldsymbol{\Sigma
}\right\vert \boldsymbol{s}_{\boldsymbol{\mu}}^{\top}(\boldsymbol{X})+\frac
{1}{\beta+1}\boldsymbol{s}_{\mathrm{vech}(\boldsymbol{\Sigma})}(\boldsymbol{X}
)\boldsymbol{s}_{\boldsymbol{\mu}}^{\top}(\boldsymbol{X})\right)
\end{aligned}
\]

\[
\begin{aligned}
&  \boldsymbol{\varsigma}_{(\beta+1)\mathrm{vech}(\boldsymbol{\Sigma}
)}(\boldsymbol{X})\boldsymbol{\varsigma}_{(\beta+1)\mathrm{vech}
(\boldsymbol{\Sigma})}^{\top}(\boldsymbol{X}) \\
&  =\frac{1}{(\beta+1)^{2}}\frac{1}{(\beta+1)^{2}}\left(  -\frac{\beta}
{2}\frac{\partial}{\partial\mathrm{vech}(\boldsymbol{\Sigma})}\log\left\vert
\boldsymbol{\Sigma}\right\vert +\boldsymbol{s}_{\mathrm{vech}
(\boldsymbol{\Sigma})}(\boldsymbol{X})\right)  \left(  -\frac{\beta}{2}
\frac{\partial}{\partial\mathrm{vech}(\boldsymbol{\Sigma})}\log\left\vert
\boldsymbol{\Sigma}\right\vert +\boldsymbol{s}_{\mathrm{vech}
(\boldsymbol{\Sigma})}(\boldsymbol{X})\right)  ^{\top} \\
&  =\frac{1}{(\beta+1)^{2}}\left(  \frac{\beta^{2}}{4(\beta+1)^{2}}
\frac{\partial}{\partial\mathrm{vech}(\boldsymbol{\Sigma})}\log\left\vert
\boldsymbol{\Sigma}\right\vert \frac{\partial}{\partial\mathrm{vech}
^{\top}(\boldsymbol{\Sigma})}\log\left\vert \boldsymbol{\Sigma}\right\vert
-\frac{\beta}{2(\beta+1)^{2}}\frac{\partial}{\partial\mathrm{vech}
(\boldsymbol{\Sigma})}\log\left\vert \boldsymbol{\Sigma}\right\vert
\boldsymbol{s}_{\mathrm{vech}(\boldsymbol{\Sigma})}^{\top}(\boldsymbol{X})\right. \\
&  \left.  -\frac{\beta}{2(\beta+1)^{2}}\boldsymbol{s}_{\mathrm{vech}
(\boldsymbol{\Sigma})}(\boldsymbol{X})\frac{\partial}{\partial\mathrm{vech}
^{\top}(\boldsymbol{\Sigma})}\log\left\vert \boldsymbol{\Sigma}\right\vert
+\frac{1}{(\beta+1)^{2}}\boldsymbol{s}_{\mathrm{vech}(\boldsymbol{\Sigma}
)}(\boldsymbol{X})\boldsymbol{s}_{\mathrm{vech}(\boldsymbol{\Sigma})}
^{\top}(\boldsymbol{X})\right)  ,
\end{aligned}
\]

its expectation is

\[
\begin{aligned}
E_{\boldsymbol{\theta}}[\boldsymbol{\varsigma}_{\boldsymbol{\vartheta}
}(\boldsymbol{X})\boldsymbol{\varsigma}_{\boldsymbol{\vartheta}}
^{\top}(\boldsymbol{X})]=\frac{1}{(\beta+1)^{2}}
\begin{pmatrix}
E_{\boldsymbol{\theta}}[\boldsymbol{s}_{\boldsymbol{\mu}}(\boldsymbol{X}
)\boldsymbol{s}_{\boldsymbol{\mu}}^{\top}(\boldsymbol{X})] & \frac{1}{\beta
+1}E_{\boldsymbol{\theta}}[\boldsymbol{s}_{\mathrm{vech}(\boldsymbol{\Sigma}
)}(\boldsymbol{X})\boldsymbol{s}_{\boldsymbol{\mu}}^{\top}(\boldsymbol{X})]\\
\frac{1}{\beta+1}E_{\boldsymbol{\theta}}[\boldsymbol{s}_{\mathrm{vech}
(\boldsymbol{\Sigma})}(\boldsymbol{X})\boldsymbol{s}_{\boldsymbol{\mu}}
^{\top}(\boldsymbol{X})] & \frac{1}{(\beta+1)^{2}}\left(  \frac{\beta^{2}}
{4}\boldsymbol{C}_{\boldsymbol{\theta}}+E_{\boldsymbol{\theta}}[\boldsymbol{s}
_{\mathrm{vech}(\boldsymbol{\Sigma})}(\boldsymbol{X})\boldsymbol{s}
_{\mathrm{vech}(\boldsymbol{\Sigma})}^{\top}(\boldsymbol{X})]\right)
\end{pmatrix}
,
\end{aligned}
\]

where

\[
\boldsymbol{C}_{\boldsymbol{\theta}}=\frac{\partial}{\partial\mathrm{vech}
(\boldsymbol{\Sigma})}\log\left\vert \boldsymbol{\Sigma}\right\vert
\frac{\partial}{\partial\mathrm{vech}^{\top}(\boldsymbol{\Sigma})}\log\left\vert
\boldsymbol{\Sigma}\right\vert =\boldsymbol{G}_{p}^{\top}\mathrm{vec}\left(
\boldsymbol{\Sigma}^{-1}\right)  \mathrm{vec}^{\top}\left(  \boldsymbol{\Sigma
}^{-1}\right)  \boldsymbol{G}_{p}
\]

with \(\boldsymbol{G}_{p}\) being the so-called ``duplication matrix''
of order \(p\), i.e., the unique \(p^{2}\times \frac{p(p+1)}{2}\)~matrix
such that
\(\mathrm{vec}(\boldsymbol{\Sigma })=\boldsymbol{G}_{p}\mathrm{vech}(\boldsymbol{\Sigma})\)
and the last derivatives are deduced from McCulloch (1982). Notice that
some terms are cancelled since they appear multiplied by expectations,
which were null,
\(E_{\boldsymbol{\theta}}[\boldsymbol{s}_{\boldsymbol{\theta}}(\boldsymbol{X} )]=\boldsymbol{0}_{p+\frac{p(p+1)}{2}}\).
Hence

\begin{equation}\protect\phantomsection\label{eq-JJJ}{
\begin{aligned}
\boldsymbol{J}_{\beta}^{\ast}(\boldsymbol{\vartheta})  &  =\frac
{(\beta+1)^{-\frac{p}{2}(\beta+1)}}{(2\pi)^{\frac{\beta p}{2}}\left\vert
\boldsymbol{\Sigma}\right\vert ^{\frac{\beta}{2}}}E_{\boldsymbol{\theta}
}[\boldsymbol{\varsigma}_{\boldsymbol{\vartheta}}(\boldsymbol{X}
)\boldsymbol{\varsigma}_{\boldsymbol{\vartheta}}^{\top}(\boldsymbol{X}
)] \\
&  =\frac{(\beta+1)^{-\frac{p}{2}(\beta+1)-2}}{(2\pi)^{\frac{\beta p}{2}
}\left\vert \boldsymbol{\Sigma}\right\vert ^{\frac{\beta}{2}}}
\begin{pmatrix}
E_{\boldsymbol{\theta}}[\boldsymbol{s}_{\boldsymbol{\mu}}(\boldsymbol{X}
)\boldsymbol{s}_{\boldsymbol{\mu}}^{\top}(\boldsymbol{X})] & \frac{1}{\beta
+1}E_{\boldsymbol{\theta}}[\boldsymbol{s}_{\boldsymbol{\mu}}(\boldsymbol{X}
)\boldsymbol{s}_{\mathrm{vech}(\boldsymbol{\Sigma})}^{\top}(\boldsymbol{X})] \\
\frac{1}{\beta+1}E_{\boldsymbol{\theta}}[\boldsymbol{s}_{\mathrm{vech}
(\boldsymbol{\Sigma})}(\boldsymbol{X})\boldsymbol{s}_{\boldsymbol{\mu}}
^{\top}(\boldsymbol{X})] & \frac{1}{(\beta+1)^{2}}\left(  \frac{\beta^{2}}
{4}\boldsymbol{C}_{\boldsymbol{\theta}}+E_{\boldsymbol{\theta}}[\boldsymbol{s}
_{\mathrm{vech}(\boldsymbol{\Sigma})}(\boldsymbol{X})\boldsymbol{s}
_{\mathrm{vech}(\boldsymbol{\Sigma})}^{\top}(\boldsymbol{X})]\right)
\end{pmatrix}
,
\end{aligned}
}\end{equation}

\[
\begin{aligned}
\boldsymbol{s}_{\boldsymbol{\mu}}(\boldsymbol{X})  &  =\frac{\partial
}{\partial\boldsymbol{\mu}}\log f_{\boldsymbol{\theta}}(\boldsymbol{X}
)=\frac{\partial}{\partial\boldsymbol{\mu}}\left(  -\frac{p}{2}\log
(2\pi)-\frac{1}{2}\log\left\vert \boldsymbol{\Sigma}\right\vert -\frac{1}
{2}(\boldsymbol{X}-\boldsymbol{\mu})^{\top}\boldsymbol{\Sigma}^{-1}
(\boldsymbol{X}-\boldsymbol{\mu})\right) \\
&  =-\frac{\partial}{\partial\boldsymbol{\mu}}\frac{1}{2}(\boldsymbol{X}
-\boldsymbol{\mu})^{\top}\boldsymbol{\Sigma}^{-1}(\boldsymbol{X}-\boldsymbol{\mu
}) \\
&  =\boldsymbol{\Sigma}^{-1}(\boldsymbol{X}-\boldsymbol{\mu}),
\end{aligned}
\]

and

\[
\begin{aligned}
E_{\boldsymbol{\theta}}[\boldsymbol{s}_{\boldsymbol{\mu}}(\boldsymbol{X}
)\boldsymbol{s}_{\boldsymbol{\mu}}^{\top}(\boldsymbol{X})]  &
=E_{\boldsymbol{\theta}}\left[  \boldsymbol{\Sigma}^{-1}(\boldsymbol{X}
-\boldsymbol{\mu})(\boldsymbol{X}-\boldsymbol{\mu})^{\top}\boldsymbol{\Sigma
}^{-1}\right] \\
&  =\boldsymbol{\Sigma}^{-1}E_{\boldsymbol{\theta}}\left[  (\boldsymbol{X}
-\boldsymbol{\mu})(\boldsymbol{X}-\boldsymbol{\mu})^{\top}\right]
\boldsymbol{\Sigma}^{-1} \\
&  =\boldsymbol{\Sigma}^{-1}\boldsymbol{\Sigma\Sigma}^{-1} \\
&  =\boldsymbol{\Sigma}^{-1}.
\end{aligned}
\]

We omit the rest of the terms, since they were calculated by McCulloch
(1982) in detailed way. Accordingly, the Fisher information matrix is

\[
\begin{aligned}
\begin{pmatrix}
E_{\boldsymbol{\theta}}[\boldsymbol{s}_{\boldsymbol{\mu}}(\boldsymbol{X}
)\boldsymbol{s}_{\boldsymbol{\mu}}^{\top}(\boldsymbol{X})] &
E_{\boldsymbol{\theta}}[\boldsymbol{s}_{\mathrm{vech}(\boldsymbol{\Sigma}
)}(\boldsymbol{X})\boldsymbol{s}_{\boldsymbol{\mu}}^{\top}(\boldsymbol{X})]\\
(E_{\boldsymbol{\theta}}[\boldsymbol{s}_{\boldsymbol{\mu}}(\boldsymbol{X}
)\boldsymbol{s}_{\mathrm{vech}(\boldsymbol{\Sigma})}^{\top}(\boldsymbol{X})])^{\top}
& E_{\boldsymbol{\theta}}[\boldsymbol{s}_{\mathrm{vech}(\boldsymbol{\Sigma}
)}(\boldsymbol{X})\boldsymbol{s}_{\mathrm{vech}(\boldsymbol{\Sigma})}
^{\top}(\boldsymbol{X})]
\end{pmatrix}
,
\end{aligned}
\]

where

\[
\begin{aligned}
E_{\boldsymbol{\theta}}[\boldsymbol{s}_{\boldsymbol{\mu}}(\boldsymbol{X}
)\boldsymbol{s}_{\boldsymbol{\mu}}^{\top}(\boldsymbol{X})]  &
=\boldsymbol{\Sigma}^{-1}, \\
E_{\boldsymbol{\theta}}[\boldsymbol{s}_{\boldsymbol{\mu}}(\boldsymbol{X}
)\boldsymbol{s}_{\mathrm{vech}(\boldsymbol{\Sigma})}^{\top}(\boldsymbol{X})]  &
=\boldsymbol{0}_{p\times\frac{p(p+1)}{2}}, \\
E_{\boldsymbol{\theta}}[\boldsymbol{s}_{\mathrm{vech}(\boldsymbol{\Sigma}
)}(\boldsymbol{X})\boldsymbol{s}_{\mathrm{vech}(\boldsymbol{\Sigma})}
^{\top}(\boldsymbol{X})]  &  =\frac{1}{2}\boldsymbol{G}_{p}^{\top}\left(
\boldsymbol{\Sigma}^{-1}\otimes\boldsymbol{\Sigma}^{-1}\right)  \boldsymbol{G}
_{p}.
\end{aligned}
\]

Hence

\begin{equation}\protect\phantomsection\label{eq-eq0}{
\begin{aligned}
\boldsymbol{J}_{\beta}^{\ast}(\boldsymbol{\vartheta})  &
=E_{\boldsymbol{\vartheta}}[\boldsymbol{\varsigma}_{\boldsymbol{\vartheta}
}(\boldsymbol{X})\boldsymbol{\varsigma}_{\boldsymbol{\vartheta}}
^{\top}(\boldsymbol{X})f_{\boldsymbol{\vartheta}}^{\beta}(\boldsymbol{X}
)] \\
\boldsymbol{J}_{\beta}^{\ast}(\boldsymbol{\vartheta})  &  =
\begin{pmatrix}
\boldsymbol{J}_{\beta}^{\ast}(\boldsymbol{\mu}) & \boldsymbol{0}_{p\times
\frac{p(p+1)}{2}} \\
\boldsymbol{0}_{\frac{p(p+1)}{2}\times p} & \boldsymbol{J}_{\beta}^{\ast
}((\beta+1)\mathrm{vech}(\boldsymbol{\Sigma}))
\end{pmatrix}
,
\end{aligned}
}\end{equation}

where

\begin{equation}\protect\phantomsection\label{eq-eq1}{
\begin{aligned}
\boldsymbol{J}_{\beta}^{\ast}(\boldsymbol{\mu})  &  =\frac{(\beta
+1)^{-\frac{p}{2}(\beta+1)-2}}{(2\pi)^{\frac{\beta p}{2}}\left\vert
\boldsymbol{\Sigma}\right\vert ^{\frac{\beta}{2}}}\boldsymbol{\Sigma}
^{-1},
\end{aligned}
}\end{equation}

\begin{equation}\protect\phantomsection\label{eq-eq2}{
\begin{aligned}
\boldsymbol{J}_{\beta}^{\ast}((\beta+1)\mathrm{vech}(\boldsymbol{\Sigma}))  &
=\frac{(\beta+1)^{-\frac{p}{2}(\beta+1)-4}}{4(2\pi)^{\frac{\beta p}{2}
}\left\vert \boldsymbol{\Sigma}\right\vert ^{\frac{\beta}{2}}}\left[
\beta^{2}\boldsymbol{C}_{\boldsymbol{\theta}}+2\boldsymbol{G}_{p}^{\top}\left(
\boldsymbol{\Sigma}^{-1}\otimes\boldsymbol{\Sigma}^{-1}\right)  \boldsymbol{G}
_{p}\right]  .
\end{aligned}
}\end{equation}

According to the original problem, from
(\ref{eq-eq0})-(\ref{eq-eq1})-(\ref{eq-eq2}) associated to
(\ref{eq-paramEta}), taking into account

\[
\begin{aligned}
\boldsymbol{J}_{\beta}^{\ast}(\boldsymbol{\mu})  &  =\frac{(\beta
+1)^{-\frac{p}{2}-1}}{(2\pi)^{\frac{\beta p}{2}}\left\vert (\beta
+1)\boldsymbol{\Sigma}\right\vert ^{\frac{\beta}{2}}}\left(  (\beta
+1)\boldsymbol{\Sigma}\right)  ^{-1}, \\
\boldsymbol{J}_{\beta}^{\ast}((\beta+1)\mathrm{vech}(\boldsymbol{\Sigma}))  &
=\frac{(\beta+1)^{-\frac{p}{2}-2}}{4(2\pi)^{\frac{\beta p}{2}}\left\vert
(\beta+1)\boldsymbol{\Sigma}\right\vert ^{\frac{\beta}{2}}}\left[  \beta
^{2}\boldsymbol{C}_{\boldsymbol{\vartheta}}+2\boldsymbol{G}_{p}^{\top}\left(
(\beta+1)\boldsymbol{\Sigma}^{-1}\otimes(\beta+1)\boldsymbol{\Sigma}
^{-1}\right)  \boldsymbol{G}_{p}\right]  ,
\end{aligned}
\]

we deduce the final expression of
\(\boldsymbol{J}_{\beta}(\boldsymbol{\theta })\) for
\(\boldsymbol{\theta}=(\boldsymbol{\mu}^{\top},\mathrm{vech} ^{\top}(\boldsymbol{\Sigma}))^{\top}\).

\[
\frac{(\beta+1)^{-\frac{p}{2}-2}}{4(2\pi)^{\frac{\beta p}{2}}\left\vert
\boldsymbol{\Sigma}\right\vert ^{\frac{\beta}{2}}}\left[  \beta^{2}
\boldsymbol{C}_{\boldsymbol{\theta}}+2\boldsymbol{G}_{p}^{\top}\left(
\boldsymbol{\Sigma}^{-1}\otimes\boldsymbol{\Sigma}^{-1}\right)  \boldsymbol{G}
_{p}\right]
\]

\subsection{\texorpdfstring{Proof of Theorem
\ref{thm-Th-Psi-BS1}}{Proof of Theorem }}\label{sec-Proof-Th-Psi-BS1}

Following the parametrization of the Proof \ref{thm-Th-J-BS1}~(see
\ref{sec-Proof-Th-J-BS1}), we obtain

\[
\begin{aligned}
\xi^{\ast}(\boldsymbol{\mu})  &  =\int_{\mathbb{R}^{p}}\boldsymbol{\varsigma}_{\boldsymbol{\mu}}(\boldsymbol{x}
)f_{\boldsymbol{\vartheta}}^{\beta+1}(\boldsymbol{x})d\boldsymbol{x} \\
&  =\frac{(\beta+1)^{-\frac{p}{2}(\beta+1)-1}}{(2\pi)^{\frac{\beta p}{2}
}\left\vert \boldsymbol{\Sigma}\right\vert ^{\frac{\beta}{2}}}\int_{\mathbb{R}^{p}}\boldsymbol{s}_{\boldsymbol{\mu}}(\boldsymbol{x})f_{\boldsymbol{\theta}
}(\boldsymbol{x})d\boldsymbol{x} \\
&  =\boldsymbol{0}_{p},
\end{aligned}
\]

\[
\begin{aligned}
\boldsymbol{\xi}^{\ast}((\beta+1)\mathrm{vech}(\boldsymbol{\Sigma}))  &
=\int_{\mathbb{R}^{p}}\boldsymbol{\varsigma}_{(\beta+1)\mathrm{vech}(\boldsymbol{\Sigma}
)}(\boldsymbol{x})f_{\boldsymbol{\vartheta}}^{\beta+1}(\boldsymbol{x}
)d\boldsymbol{x} \\
&  =\frac{(\beta+1)^{-\frac{p}{2}(\beta+1)-2}}{(2\pi)^{\frac{\beta p}{2}
}\left\vert \boldsymbol{\Sigma}\right\vert ^{\frac{\beta}{2}}}\int_{\mathbb{R}^{p}}\left(  -\frac{\beta}{2}\frac{\partial}{\partial\mathrm{vech}
(\boldsymbol{\Sigma})}\log\left\vert \boldsymbol{\Sigma}\right\vert
+\boldsymbol{s}_{\mathrm{vech}(\boldsymbol{\Sigma})}(\boldsymbol{x})\right)
f_{\boldsymbol{\theta}}(\boldsymbol{x})d\boldsymbol{x} \\
&  =-\frac{\beta}{2}\frac{(\beta+1)^{-\frac{p}{2}(\beta+1)-2}}{(2\pi
)^{\frac{\beta p}{2}}\left\vert \boldsymbol{\Sigma}\right\vert ^{\frac{\beta
}{2}}}\frac{\partial}{\partial\mathrm{vech}(\boldsymbol{\Sigma})}
\log\left\vert \boldsymbol{\Sigma}\right\vert \\
&  =-\frac{\beta}{2}\frac{(\beta+1)^{-\frac{p}{2}(\beta+1)-2}}{(2\pi
)^{\frac{\beta p}{2}}\left\vert \boldsymbol{\Sigma}\right\vert ^{\frac{\beta
}{2}}}\boldsymbol{G}_{p}^{\top}\mathrm{vec}\left(  \boldsymbol{\Sigma}
^{-1}\right) \\
&  =-\frac{\beta}{2}\frac{(\beta+1)^{-\left(  \frac{p}{2}+1\right)  }}
{(2\pi)^{\frac{\beta p}{2}}\left\vert (\beta+1)\boldsymbol{\Sigma}\right\vert
^{\frac{\beta}{2}}}\boldsymbol{G}_{p}^{\top}\mathrm{vec}\left(  \left(
(\beta+1)\boldsymbol{\Sigma}\right)  ^{-1}\right)
\end{aligned}
\]

and their cross products

\[
\xi^{\ast}(\boldsymbol{\mu})\xi^{\ast T}(\boldsymbol{\mu})=\boldsymbol{0}
_{p\times p},
\]

\[
\begin{aligned}
&  \boldsymbol{\xi}^{\ast}((\beta+1)\mathrm{vech}(\boldsymbol{\Sigma
}))\boldsymbol{\xi}^{\ast T}((\beta+1)\mathrm{vech}(\boldsymbol{\Sigma})) \\
&  =\frac{\beta^{2}}{4}\frac{(\beta+1)^{-p-2}}{(2\pi)^{\beta p}\left\vert
(\beta+1)\boldsymbol{\Sigma}\right\vert ^{\beta}}\boldsymbol{G}_{p}
^{\top}\mathrm{vec}\left(  \left(  (\beta+1)\boldsymbol{\Sigma}\right)
^{-1}\right)  \mathrm{vec}^{\top}\left(  \left(  (\beta+1)\boldsymbol{\Sigma
}\right)  ^{-1}\right)  \boldsymbol{G}_{p} \\
&  =\frac{\beta^{2}}{4}\frac{(\beta+1)^{-(p+2)}}{(2\pi)^{\beta p}\left\vert
(\beta+1)\boldsymbol{\Sigma}\right\vert ^{\beta}}\boldsymbol{C}
_{\boldsymbol{\vartheta}}.
\end{aligned}
\]

Translating these terms to the original parametrization
\(\boldsymbol{\theta }=(\boldsymbol{\mu}^{\top},\mathrm{vech}^{\top}(\boldsymbol{\Sigma}))^{\top}\),
we obtain the final expression of
\(\boldsymbol{\xi}_{\beta}(\boldsymbol{\theta})\) and
\(\boldsymbol{\xi}_{\beta}(\boldsymbol{\theta})\boldsymbol{\xi}_{\beta} ^{\top}(\boldsymbol{\theta})\).

\subsection{\texorpdfstring{Proof of Corollary
\ref{cor-Th-Psi-KS1}}{Proof of Corollary }}\label{sec-ProofTh-Psi-KS1}

From the (\ref{eq-eqJ})-(\ref{eq-eqJ11})-(\ref{eq-eqJ22}) and
((\ref{eq-psipsi}) )-(\ref{eq-psipsi2}) we can obtain

\[
\begin{aligned}
&  \boldsymbol{K}_{\beta}(\mathrm{vech}(\boldsymbol{\Sigma}))=\boldsymbol{J}
_{2\beta}(\mathrm{vech}(\boldsymbol{\Sigma}))-\boldsymbol{\xi}_{\beta
}(\mathrm{vech}(\boldsymbol{\Sigma}))\boldsymbol{\xi}_{\beta}^{\top}
(\mathrm{vech}(\boldsymbol{\Sigma})) \\
&  =\frac{(2\beta+1)^{-\frac{p}{2}-2}}{4(2\pi)^{\beta p}\left\vert
\boldsymbol{\Sigma}\right\vert ^{\beta}}\left[  4\beta^{2}\boldsymbol{C}
_{\boldsymbol{\theta}}+2\boldsymbol{G}_{p}^{\top}\left(  \boldsymbol{\Sigma}
^{-1}\otimes\boldsymbol{\Sigma}^{-1}\right)  \boldsymbol{G}_{p}\right]
-\frac{\beta^{2}(\beta+1)^{-(p+2)}}{4(2\pi)^{\beta p}\left\vert
\boldsymbol{\Sigma}\right\vert ^{\beta}}\boldsymbol{C}_{\boldsymbol{\theta}} \\
&  =\frac{1}{4(2\pi)^{\beta p}\left\vert \boldsymbol{\Sigma}\right\vert
^{\beta}}\left\{  2(2\beta+1)^{-\frac{p}{2}-2}\boldsymbol{G}_{p}^{\top}\left(
\boldsymbol{\Sigma}^{-1}\otimes\boldsymbol{\Sigma}^{-1}\right)  \boldsymbol{G}
_{p}+\beta^{2}\left[  4(2\beta+1)^{-\frac{p}{2}-2}-(\beta+1)^{-(p+2)}\right]
\boldsymbol{C}_{\boldsymbol{\theta}}\right\} \\
&  =\frac{1}{4(2\pi)^{\beta p}\left\vert \boldsymbol{\Sigma}\right\vert
^{\beta}}\boldsymbol{G}_{p}^{\top}\left[  \overline{\boldsymbol{J}}_{2\beta
}(\boldsymbol{\Sigma}^{-1})+\overline{\boldsymbol{\xi}}_{\beta}
(\boldsymbol{\Sigma}^{-1})\overline{\boldsymbol{\xi}}_{\beta}^{\top}
(\boldsymbol{\Sigma}^{-1})\right]  \boldsymbol{G}_{p}.
\end{aligned}
\]

\subsection{\texorpdfstring{Proof of Proposition
\ref{prp-propK}}{Proof of Proposition }}\label{sec-ProofPropK}

For calculating (\ref{eq-inv}) we take into account the following
version of the Woodbury's formula

\[
\left(  \boldsymbol{G}+\boldsymbol{uv}^{\top}\right)  ^{-1}=\boldsymbol{G}
^{-1}-(1+\boldsymbol{v}^{\top}\boldsymbol{G}^{-1}\boldsymbol{u})^{-1}
\boldsymbol{G}^{-1}\boldsymbol{uv}^{\top}\boldsymbol{G}^{-1},
\]

i.e.

\[
\begin{aligned}
&  \left(  \overline{\boldsymbol{J}}_{2\beta}(\boldsymbol{R}_{0}
^{-1})+\overline{\boldsymbol{\xi}}_{\beta}(\boldsymbol{R}_{0}^{-1}
)\overline{\boldsymbol{\xi}}_{\beta}^{\top}(\boldsymbol{R}_{0}^{-1})\right)
^{-1}=\overline{\boldsymbol{J}}_{2\beta}^{-1}(\boldsymbol{R}_{0}^{-1}
)-\frac{\overline{\boldsymbol{J}}_{2\beta}^{-1}(\boldsymbol{R}_{0}
^{-1})\overline{\boldsymbol{\xi}}_{\beta}(\boldsymbol{R}_{0}^{-1}
)\overline{\boldsymbol{\xi}}_{\beta}^{\top}(\boldsymbol{R}_{0}^{-1}
)\overline{\boldsymbol{J}}_{2\beta}^{-1}(\boldsymbol{R}_{0}^{-1})}
{1+\overline{\boldsymbol{\xi}}_{\beta}^{\top}(\boldsymbol{R}_{0}^{-1}
)\overline{\boldsymbol{J}}_{2\beta}^{-1}(\boldsymbol{R}_{0}^{-1}
)\overline{\boldsymbol{\xi}}_{\beta}(\boldsymbol{R}_{0}^{-1})} \\
&  =\kappa_{1}^{-1}(p,\beta)\left(  \left(  \boldsymbol{R}_{0}\otimes
\boldsymbol{R}_{0}\right)  -\frac{\kappa_{3}(p,\beta)\mathrm{vec}
(\boldsymbol{R}_{0})\mathrm{vec}^{\top}(\boldsymbol{R}_{0})}{1+p\kappa
_{3}(p,\beta)}\right)  ,
\end{aligned}
\]

where \(\kappa_{3}(p,\beta)\) is (\ref{eq-kappa3}). In the last equality

\[
\begin{aligned}
\overline{\boldsymbol{J}}_{2\beta}^{-1}(\boldsymbol{R}_{0}^{-1})\overline
{\boldsymbol{\xi}}_{\beta}(\boldsymbol{R}_{0}^{-1})  &  =\kappa_{1}
^{-1}(p,\beta)\kappa_{2}(p,\beta)\left(  \boldsymbol{R}_{0}\otimes
\boldsymbol{R}_{0}\right)  \mathrm{vec}(\boldsymbol{R}_{0}^{-1}) \\
&  =\kappa_{1}^{-1}(p,\beta)\kappa_{2}(p,\beta)\mathrm{vec}(\boldsymbol{R}
_{0}\boldsymbol{R}_{0}^{-1}\boldsymbol{R}_{0}) \\
&  =\kappa_{1}^{-1}(p,\beta)\kappa_{2}(p,\beta)\mathrm{vec}(\boldsymbol{R}
_{0}),
\end{aligned}
\]

and

\[
\begin{aligned}
\overline{\boldsymbol{\xi}}_{\beta}^{\top}(\boldsymbol{R}_{0}^{-1})\overline
{\boldsymbol{J}}_{2\beta}^{-1}(\boldsymbol{R}_{0}^{-1})\overline
{\boldsymbol{\xi}}_{\beta}(\boldsymbol{R}_{0}^{-1})  &  =\kappa_{1}
^{-1}(p,\beta)\kappa_{2}^{2}(p,\beta)\mathrm{vec}^{\top}(\boldsymbol{R}
_{0})\mathrm{vec}(\boldsymbol{R}_{0}^{-1}) \\
&  =\kappa_{1}^{-1}(p,\beta)\kappa_{2}^{2}(p,\beta)\mathrm{vec}^{\top}
(\boldsymbol{R}_{0}^{\frac{1}{2}}\boldsymbol{I}_{p}\boldsymbol{R}_{0}
^{\frac{1}{2}})\mathrm{vec}(\boldsymbol{R}_{0}^{-\frac{1}{2}}\boldsymbol{I}
_{p}\boldsymbol{R}_{0}^{-\frac{1}{2}}) \\
&  =\kappa_{1}^{-1}(p,\beta)\kappa_{2}^{2}(p,\beta)\mathrm{vec}^{\top}
(\boldsymbol{I}_{p})(\boldsymbol{R}_{0}^{\frac{1}{2}}\otimes\boldsymbol{R}
_{0}^{\frac{1}{2}})(\boldsymbol{R}_{0}^{-\frac{1}{2}}\otimes\boldsymbol{R}
_{0}^{-\frac{1}{2}})\mathrm{vec}(\boldsymbol{I}_{p}) \\
&  =\kappa_{1}^{-1}(p,\beta)\kappa_{2}^{2}(p,\beta)\mathrm{vec}^{\top}
(\boldsymbol{I}_{p})\mathrm{vec}(\boldsymbol{I}_{p}) \\
&  =\kappa_{1}^{-1}(p,\beta)\kappa_{2}^{2}(p,\beta)\mathrm{trace}
(\boldsymbol{I}_{p}) \\
&  =\kappa_{3}(p,\beta)p.
\end{aligned}
\]

\subsection{\texorpdfstring{Proof of Theorem
\ref{thm-Th-U-BS1}}{Proof of Theorem }}\label{sec-ProofTh-U-BS1}

The \(\beta\)-score function is

\[
\begin{aligned}
\boldsymbol{U}_{\beta,n}\left(  \boldsymbol{\mu}\right)   &  =\frac{1}{n}
\sum_{i=1}^{n}f_{\theta}^{\beta}(\boldsymbol{X}_{i})\boldsymbol{s}
_{\boldsymbol{\mu}}(\boldsymbol{X}_{i})-\boldsymbol{\xi}_{\beta}
(\boldsymbol{\mu}) \\
&  =\frac{1}{n}\boldsymbol{\Sigma}^{-1}\sum_{i=1}^{n}f_{\theta}^{\beta
}(\boldsymbol{X}_{i})(\boldsymbol{X}_{i}-\boldsymbol{\mu}) \\
&  =\frac{1}{(2\pi)^{\frac{\beta p}{2}}\left\vert \boldsymbol{\Sigma
}\right\vert ^{\frac{\beta}{2}}}\boldsymbol{\Sigma}^{-1}\frac{1}{n}\sum
_{i=1}^{n}W_{i,\beta}(\boldsymbol{\theta})(\boldsymbol{X}_{i}-\boldsymbol{\mu}),
\end{aligned}
\]

with

\[
f_{\theta}^{\beta}(\boldsymbol{X}_{i})=\frac{1}{(2\pi)^{\frac{\beta p}{2}
}\left\vert \boldsymbol{\Sigma}\right\vert ^{\frac{\beta}{2}}}W_{i,\beta
}(\boldsymbol{\theta}),
\]

\(W_{i,\beta}(\boldsymbol{\theta})\) is (\ref{eq-w}) and

\begin{equation}\protect\phantomsection\label{eq-vecU0}{
\begin{aligned}
\boldsymbol{U}_{\beta,n}\left(  \mathrm{vech}(\boldsymbol{\Sigma})\right)   &
=\frac{1}{n}\sum_{i=1}^{n}f_{\theta}^{\beta}(\boldsymbol{X}_{i})\boldsymbol{s}
_{\mathrm{vech}(\boldsymbol{\Sigma})}(\boldsymbol{X}_{i})-\boldsymbol{\xi
}_{\beta}(\boldsymbol{\Sigma}) \\
&  =-\frac{1}{2(2\pi)^{\frac{\beta p}{2}}\left\vert \boldsymbol{\Sigma
}\right\vert ^{\frac{\beta}{2}}}\boldsymbol{G}_{p}^{\top}\mathrm{vec}
(\boldsymbol{\Sigma}^{-1})\frac{1}{n}\sum_{i=1}^{n}W_{i,\beta}
(\boldsymbol{\theta}) \\
&  +\frac{1}{2(2\pi)^{\frac{\beta p}{2}}\left\vert \boldsymbol{\Sigma
}\right\vert ^{\frac{\beta}{2}}}\boldsymbol{G}_{p}^{\top}\left(
\boldsymbol{\Sigma}^{-1}\otimes\boldsymbol{\Sigma}^{-1}\right)  \frac{1}
{n}\sum_{i=1}^{n}W_{i,\beta}(\boldsymbol{\theta})\left(  (\boldsymbol{X}
_{i}-\boldsymbol{\mu})\otimes(\boldsymbol{X}_{i}-\boldsymbol{\mu})\right) \\
&  +\frac{\beta}{2}\frac{(\beta+1)^{-(\frac{p}{2}+1)}}{(2\pi)^{\frac{\beta
p}{2}}\left\vert \boldsymbol{\Sigma}\right\vert ^{\frac{\beta}{2}}
}\boldsymbol{G}_{p}^{\top}\mathrm{vec}\left(  \boldsymbol{\Sigma}^{-1}\right) \\
&  =-\frac{1}{2(2\pi)^{\frac{\beta p}{2}}\left\vert \boldsymbol{\Sigma
}\right\vert ^{\frac{\beta}{2}}}\boldsymbol{G}_{p}^{\top}\mathrm{vec}\left(
\boldsymbol{\Sigma}^{-1}\right)  \left(  \frac{1}{n}\sum_{i=1}^{n}W_{i,\beta
}(\theta)-\beta(\beta+1)^{-(\frac{p}{2}+1)}\right) \\
&  +\frac{1}{2(2\pi)^{\frac{\beta p}{2}}\left\vert \boldsymbol{\Sigma
}\right\vert ^{\frac{\beta}{2}}}\boldsymbol{G}_{p}^{\top}\frac{1}{n}\sum
_{i=1}^{n}W_{i,\beta}(\boldsymbol{\theta})\left(  \boldsymbol{\Sigma}
^{-1}(\boldsymbol{X}_{i}-\boldsymbol{\mu})\right)  \otimes\left(
\boldsymbol{\Sigma}^{-1}(\boldsymbol{X}_{i}-\boldsymbol{\mu})\right) \\
&  =\boldsymbol{G}_{p}^{\top}\boldsymbol{V}_{\beta,n}\left(  \mathrm{vec}
(\boldsymbol{\Sigma})\right)
\end{aligned}
}\end{equation}

with
\(\boldsymbol{V}_{\beta,n}\left(  \mathrm{vec}(\boldsymbol{\Sigma })\right)\)
given in (\ref{eq-vecU}).

\subsection{\texorpdfstring{Proof of Proposition
\ref{prp-PropV}}{Proof of Proposition }}\label{proof-of-proposition--prp-propv}

According with (\ref{eq-vecU0}),

\[
\begin{aligned}
&  2(2\pi)^{\frac{\beta p}{2}}\left\vert \widetilde{\boldsymbol{\Sigma}
}_{\beta}\right\vert ^{\frac{\beta}{2}}\boldsymbol{V}_{\beta,n}
(\widetilde{\boldsymbol{\Lambda}}_{\beta},\boldsymbol{R}_{0}) \\
&  =(\widetilde{\boldsymbol{\Lambda}}_{\beta}^{-1/2}\otimes
\widetilde{\boldsymbol{\Lambda}}_{\beta}^{-1/2})(\boldsymbol{R}_{0}
^{-1}\otimes\boldsymbol{R}_{0}^{-1})\left[  \frac{1}{n}\sum_{i=1}
^{n}\widetilde{W}_{i,\beta}\left(  \widetilde{\boldsymbol{\Lambda}}_{\beta
}^{-1/2}(\boldsymbol{X}_{i}-\widetilde{\boldsymbol{\mu}}_{\beta})\right)
\otimes\left(  \widetilde{\boldsymbol{\Lambda}}_{\beta}^{-1/2}(\boldsymbol{X}
_{i}-\widetilde{\boldsymbol{\mu}}_{\beta})\right)  \right] \\
&  -\widetilde{\kappa}_{0}(p,\beta)(\widetilde{\boldsymbol{\Lambda}}_{\beta
}^{-1/2}\otimes\widetilde{\boldsymbol{\Lambda}}_{\beta}^{-1/2})\mathrm{vec}
\left(  \boldsymbol{R}_{0}^{-1}\right) \\
&  =\widetilde{\kappa}_{0}(p,\beta)(\widetilde{\boldsymbol{\Lambda}}_{\beta
}^{-1/2}\otimes\widetilde{\boldsymbol{\Lambda}}_{\beta}^{-1/2})\left[
(\boldsymbol{R}_{0}^{-1}\otimes\boldsymbol{R}_{0}^{-1})\mathrm{vec}\left(
\widetilde{\boldsymbol{R}}_{\boldsymbol{X},\beta}\right)  -\mathrm{vec}\left(
\boldsymbol{R}_{0}^{-1}\right)  \right] \\
&  =\widetilde{\kappa}_{0}(p,\beta)(\widetilde{\boldsymbol{\Lambda}}_{\beta
}^{-1/2}\otimes\widetilde{\boldsymbol{\Lambda}}_{\beta}^{-1/2})\left[
\mathrm{vec}\left(  \boldsymbol{R}_{0}^{-1}\widetilde{\boldsymbol{R}
}_{\boldsymbol{X},\beta}\boldsymbol{R}_{0}^{-1}\right)  -\mathrm{vec}\left(
\boldsymbol{R}_{0}^{-1}\boldsymbol{R}_{0}\boldsymbol{R}_{0}^{-1}\right)
\right] \\
&  =\widetilde{\kappa}_{0}(p,\beta)(\widetilde{\boldsymbol{\Lambda}}_{\beta
}^{-1/2}\otimes\widetilde{\boldsymbol{\Lambda}}_{\beta}^{-1/2})(\boldsymbol{R}
_{0}^{-1}\otimes\boldsymbol{R}_{0}^{-1})\left[  \mathrm{vec}\left(
\widetilde{\boldsymbol{R}}_{\boldsymbol{X},\beta}\right)  -\mathrm{vec}\left(
\boldsymbol{R}_{0}\right)  \right]  .
\end{aligned}
\]

\subsection{\texorpdfstring{Proof of Lemma
\ref{lem-short}}{Proof of Lemma }}\label{sec-ProofLemShort}

The matrices \(\boldsymbol{J}_{\beta}(\boldsymbol{\theta})\) and
\(\boldsymbol{K}_{\beta}(\boldsymbol{\theta})\) are block diagonal with
respect to \(\boldsymbol{\mu}\) and
\(\mathrm{vech}(\boldsymbol{\Sigma})\), by Theorem~\ref{thm-Th-J-BS1}
and Corollary~\ref{cor-Th-Psi-KS1}, and so are the corresponding
matrices for \(\boldsymbol{\eta}\), so it suffices to consider the block
of the variance-covariance parameters. Write
\(\boldsymbol{A}=\boldsymbol{G}_{p}^{\top}(\boldsymbol{\Sigma}^{-1}\otimes\boldsymbol{\Sigma}^{-1})\boldsymbol{G}_{p}\)
and \(\boldsymbol{s}=\mathrm{vech}(\boldsymbol{\Sigma})\). Since
\(\boldsymbol{G}_{p}^{\top}\mathrm{vec}(\boldsymbol{\Sigma}^{-1})=\boldsymbol{G}_{p}^{\top}(\boldsymbol{\Sigma}^{-1}\otimes\boldsymbol{\Sigma}^{-1})\mathrm{vec}(\boldsymbol{\Sigma})=\boldsymbol{A}\boldsymbol{s}\),
(\ref{eq-eqJ22}), Corollary~\ref{cor-Th-Psi-KS1} and
\(\boldsymbol{C}_{\boldsymbol{\theta}}=\boldsymbol{A}\boldsymbol{s}\boldsymbol{s}^{\top}\boldsymbol{A}\)
give \[
\boldsymbol{J}_{\beta}(\mathrm{vech}(\boldsymbol{\Sigma}))=a_{\beta}\left(2\boldsymbol{A}+\beta^{2}\boldsymbol{A}\boldsymbol{s}\boldsymbol{s}^{\top}\boldsymbol{A}\right),
\qquad
\boldsymbol{K}_{\beta}(\mathrm{vech}(\boldsymbol{\Sigma}))=b_{\beta}\left(\kappa_{1}(p,\beta)\boldsymbol{A}+\kappa_{2}(p,\beta)\boldsymbol{A}\boldsymbol{s}\boldsymbol{s}^{\top}\boldsymbol{A}\right),
\] for positive scalars \(a_{\beta}\) and \(b_{\beta}\). Let
\(\boldsymbol{N}\) and \(\boldsymbol{T}\) be the Jacobian matrices of
\(\mathrm{vech}(\boldsymbol{\Sigma})\) with respect to the nuisance and
to the remaining variance-covariance parameters, so that \[
\boldsymbol{J}_{\beta,21}(\boldsymbol{\eta})=\boldsymbol{T}^{\top}\boldsymbol{J}_{\beta}(\mathrm{vech}(\boldsymbol{\Sigma}))\boldsymbol{N},\qquad
\boldsymbol{J}_{\beta,11}(\boldsymbol{\eta})=\boldsymbol{N}^{\top}\boldsymbol{J}_{\beta}(\mathrm{vech}(\boldsymbol{\Sigma}))\boldsymbol{N},
\] and the same for \(\boldsymbol{K}_{\beta}\). Since
\(\sigma_{ij}=\rho_{ij}\sigma_{i}\sigma_{j}\) is homogeneous of degree
one in the variances, Euler's theorem gives
\(\boldsymbol{s}=\boldsymbol{N}\boldsymbol{c}\), with \(\boldsymbol{c}\)
the vector formed by the variances and zeros in the remaining nuisance
components. With
\(\boldsymbol{A}_{11}=\boldsymbol{N}^{\top}\boldsymbol{A}\boldsymbol{N}\),
which is nonsingular, and
\(\boldsymbol{A}_{21}=\boldsymbol{T}^{\top}\boldsymbol{A}\boldsymbol{N}\),
\[
\boldsymbol{T}^{\top}\left(2\boldsymbol{A}+\beta^{2}\boldsymbol{A}\boldsymbol{s}\boldsymbol{s}^{\top}\boldsymbol{A}\right)\boldsymbol{N}
=2\boldsymbol{A}_{21}+\beta^{2}\boldsymbol{A}_{21}\boldsymbol{c}\boldsymbol{c}^{\top}\boldsymbol{A}_{11}
=\boldsymbol{A}_{21}\boldsymbol{A}_{11}^{-1}\left(2\boldsymbol{A}_{11}+\beta^{2}\boldsymbol{A}_{11}\boldsymbol{c}\boldsymbol{c}^{\top}\boldsymbol{A}_{11}\right),
\] that is,
\(\boldsymbol{J}_{\beta,21}(\boldsymbol{\eta})=\boldsymbol{A}_{21}\boldsymbol{A}_{11}^{-1}\boldsymbol{J}_{\beta,11}(\boldsymbol{\eta})\).
The same computation with \(\kappa_{1}(p,\beta)\) and
\(\kappa_{2}(p,\beta)\) in place of \(2\) and \(\beta^{2}\) gives
\(\boldsymbol{K}_{\beta,21}(\boldsymbol{\eta})=\boldsymbol{A}_{21}\boldsymbol{A}_{11}^{-1}\boldsymbol{K}_{\beta,11}(\boldsymbol{\eta})\),
and both ratios are equal to
\(\boldsymbol{A}_{21}\boldsymbol{A}_{11}^{-1}\). Then
(\ref{eq-Vgeneral}) reduces to
\(\boldsymbol{K}_{\beta,22}(\boldsymbol{\eta})-\boldsymbol{K}_{\beta,21}(\boldsymbol{\eta})\boldsymbol{K}_{\beta,11}^{-1}(\boldsymbol{\eta})\boldsymbol{K}_{\beta,12}(\boldsymbol{\eta})\),
whose inverse is the block of
\(\boldsymbol{K}_{\beta}^{-1}(\boldsymbol{\eta})\) corresponding to
\(\boldsymbol{\vartheta}_{2}\), and the first statement follows from
\(\boldsymbol{U}_{\beta,n,1}(\widetilde{\boldsymbol{\eta}}_{\beta})=\boldsymbol{0}\).
For the second one, if \(\boldsymbol{D}\) is the nonsingular Jacobian
matrix of \(\boldsymbol{\theta}\) with respect to \(\boldsymbol{\eta}\),
then
\(\boldsymbol{U}_{\beta,n}(\boldsymbol{\eta})=\boldsymbol{D}^{\top}\boldsymbol{U}_{\beta,n}(\boldsymbol{\theta})\)
and
\(\boldsymbol{K}_{\beta}(\boldsymbol{\eta})=\boldsymbol{D}^{\top}\boldsymbol{K}_{\beta}(\boldsymbol{\theta})\boldsymbol{D}\),
so that
\(\boldsymbol{U}_{\beta,n}^{\top}(\boldsymbol{\eta})\boldsymbol{K}_{\beta}^{-1}(\boldsymbol{\eta})\boldsymbol{U}_{\beta,n}(\boldsymbol{\eta})=\boldsymbol{U}_{\beta,n}^{\top}(\boldsymbol{\theta})\boldsymbol{K}_{\beta}^{-1}(\boldsymbol{\theta})\boldsymbol{U}_{\beta,n}(\boldsymbol{\theta})\).

\subsection{\texorpdfstring{Proof of Theorem
\ref{thm-RestrMDPDs}}{Proof of Theorem }}\label{sec-ProofRestrMDPDs}

For known correlation matrix, \(\boldsymbol{R=R}_{0}\), the variance
covariance matrix has less parameters to be estimated,

\[
\begin{aligned}
\boldsymbol{\Sigma}(\boldsymbol{\Lambda})  &  =\boldsymbol{\Lambda}
^{1/2}\boldsymbol{R}_{0}\boldsymbol{\Lambda}^{1/2}, \\
\boldsymbol{\Lambda}  &  =\mathrm{diag}\{\sigma_{j}^{2}\}_{j=1}^{p}.
\end{aligned}
\]

The estimating equations must be cautiously calculated, in comparison
with the non-restricted estimators.

The MLEs (\(\beta=0\)) are obtained as solution of

\begin{equation}\protect\phantomsection\label{eq-MLEeq1Aux}{
\begin{aligned}
\sum_{i=1}^{n}\boldsymbol{s}_{\boldsymbol{\mu}}(\boldsymbol{X}_{i})  &
=\boldsymbol{0}_{p},
\end{aligned}
}\end{equation}

\begin{equation}\protect\phantomsection\label{eq-MLEeq2Aux}{
\begin{aligned}
\sum_{i=1}^{n}\boldsymbol{s}_{\boldsymbol{\Lambda}}(\boldsymbol{X}_{i})  &
=\boldsymbol{0}_{p\times p},
\end{aligned}
}\end{equation}

where

\[
\boldsymbol{s}_{\boldsymbol{\Lambda}}(\boldsymbol{X})=-\frac{1}{2}
\frac{\partial}{\partial\boldsymbol{\Lambda}}\left(  \log\left\vert
\boldsymbol{\Lambda}^{1/2}\boldsymbol{R}_{0}\boldsymbol{\Lambda}
^{1/2}\right\vert +(\boldsymbol{X}-\boldsymbol{\mu})^{\top}\boldsymbol{\Lambda
}^{-1/2}\boldsymbol{R}_{0}^{-1}\boldsymbol{\Lambda}^{-1/2}(\boldsymbol{X}
-\boldsymbol{\mu})\right)  .
\]

The last term's derivation is

\[
\begin{aligned}
\boldsymbol{s}_{\boldsymbol{\Lambda}}(\boldsymbol{X})  &  =-\frac{1}{2}\left(
\frac{\partial}{\partial\boldsymbol{\Lambda}}\log\left\vert
\boldsymbol{\Lambda}\right\vert +2\boldsymbol{R}_{0}^{-1}\boldsymbol{\Lambda
}^{-1/2}(\boldsymbol{X}-\boldsymbol{\mu})(\boldsymbol{X}-\boldsymbol{\mu}
)^{\top}\frac{\partial}{\partial\boldsymbol{\Lambda}}\boldsymbol{\Lambda}
^{-1/2}\right) \\
&  =-\frac{1}{2}\left(  \boldsymbol{\Lambda}^{-1}-\boldsymbol{R}_{0}
^{-1}\boldsymbol{\Lambda}^{-1/2}(\boldsymbol{X}-\boldsymbol{\mu}
)(\boldsymbol{X}-\boldsymbol{\mu})^{\top}\boldsymbol{\Lambda}^{-3/2}\right) \\
&  =-\frac{1}{2}\left(  \boldsymbol{I}_{p}-\boldsymbol{R}_{0}^{-1}
\boldsymbol{\Lambda}^{-1/2}(\boldsymbol{X}-\boldsymbol{\mu})(\boldsymbol{X}
-\boldsymbol{\mu})^{\top}\boldsymbol{\Lambda}^{-1/2}\right)  \boldsymbol{\Lambda
}^{-1},
\end{aligned}
\]

hence from (\ref{eq-MLEeq2Aux}) it holds that the MLE of
\(\boldsymbol{\Lambda}\) under \(\boldsymbol{R=R}_{0}\),
\(\widetilde{\boldsymbol{\Lambda}}\), is the solution of

\[
-\frac{n}{2}\left(  \boldsymbol{I}_{p}-\boldsymbol{R}_{0}^{-1}\boldsymbol{R}
_{\boldsymbol{X}}(\boldsymbol{\Lambda})\right)  \boldsymbol{\Lambda}
^{-1}=\boldsymbol{0}_{p\times p},
\]

where

\[
\begin{aligned}
\boldsymbol{R}_{\boldsymbol{X}}(\boldsymbol{\Lambda})  &  =\boldsymbol{\Lambda
}^{-1/2}\boldsymbol{S}_{\boldsymbol{X}}\boldsymbol{\Lambda}^{-1/2}, \\
\boldsymbol{S}_{\boldsymbol{X}}  &  =\frac{1}{n}\sum_{i=1}^{n}(\boldsymbol{X}
_{i}-\boldsymbol{\bar{X}}_{n})(\boldsymbol{X}_{i}-\boldsymbol{\bar{X}}
_{n})^{\top},
\end{aligned}
\]

from which are estimated the diagonal elements of
\(\boldsymbol{\Lambda}\), as solution in \(\boldsymbol{\Lambda}\), of

\begin{equation}\protect\phantomsection\label{eq-eqRR}{
\boldsymbol{1}_{p}=\mathrm{diag}\{\boldsymbol{R}_{0}^{-1}\boldsymbol{R}
_{\boldsymbol{X}}(\boldsymbol{\Lambda})\}\boldsymbol{1}_{p}
}\end{equation}

or equivalently

\begin{equation}\protect\phantomsection\label{eq-eqRR2}{
\boldsymbol{1}_{p}=\mathrm{diag}\{\boldsymbol{R}_{0}^{-1}\Xi
(\boldsymbol{\Lambda})\boldsymbol{R}_{\boldsymbol{X}}\Xi^{\top}
(\boldsymbol{\Lambda})\}\boldsymbol{1}_{p},
}\end{equation}

where

\[
\begin{aligned}
\boldsymbol{R}_{\boldsymbol{X}}  &  =\mathrm{diag}^{-\frac{1}{2}
}\{\boldsymbol{S}_{\boldsymbol{X}}\}\boldsymbol{S}_{\boldsymbol{X}
}\mathrm{diag}^{-\frac{1}{2}}\{\boldsymbol{S}_{\boldsymbol{X}}\}, \\
\Xi(\boldsymbol{\Lambda})  &  =\boldsymbol{\Lambda}^{-\frac{1}{2}
}\mathrm{diag}^{\frac{1}{2}}\{\boldsymbol{S}_{\boldsymbol{X}}\}=\mathrm{diag}\{\tfrac
{S_{j}}{\sigma_{j}}\}_{j=1}^{p}.
\end{aligned}
\]

The minimum DPD estimators are obtained as a solution of the system the
non-linear system of equations

\[
\begin{aligned}
\frac{1}{n}\sum_{i=1}^{n}f_{\theta}^{\beta}(\boldsymbol{X}_{i})\boldsymbol{s}
_{\boldsymbol{\mu}}(\boldsymbol{X}_{i})-\boldsymbol{\xi}_{\beta}
(\boldsymbol{\mu})  &  =\boldsymbol{0}_{p}, \\
\frac{1}{n}\sum_{i=1}^{n}f_{\theta}^{\beta}(\boldsymbol{X}_{i})\boldsymbol{s}
_{\boldsymbol{\Lambda}}(\boldsymbol{X}_{i})-\boldsymbol{\xi}_{\beta
}(\boldsymbol{\Lambda})  &  =\boldsymbol{0}_{p\times p},
\end{aligned}
\]

where

\[
\begin{aligned}
\boldsymbol{s}_{\boldsymbol{\mu}}(\boldsymbol{X})  &  =\boldsymbol{\Lambda
}^{-1/2}\boldsymbol{R}_{0}^{-1}\boldsymbol{\Lambda}^{-1/2}(\boldsymbol{X}
-\boldsymbol{\mu}), \\
\boldsymbol{\xi}_{\beta}(\boldsymbol{\mu})  &  =E\left[  f_{\theta}^{\beta
}(\boldsymbol{X})\boldsymbol{s}_{\boldsymbol{\mu}}(\boldsymbol{X})\right]
=\boldsymbol{0}_{p}, \\
\boldsymbol{s}_{\boldsymbol{\Lambda}}(\boldsymbol{X})  &  =-\frac{1}{2}\left(
\boldsymbol{I}_{p}-\boldsymbol{R}_{0}^{-1}\boldsymbol{\Lambda}^{-1/2}
(\boldsymbol{X}-\boldsymbol{\mu})(\boldsymbol{X}-\boldsymbol{\mu}
)^{\top}\boldsymbol{\Lambda}^{-1/2}\right)  \boldsymbol{\Lambda}^{-1}, \\
\boldsymbol{\xi}_{\beta}(\boldsymbol{\Lambda})  &  =E\left[  f_{\theta}
^{\beta}(\boldsymbol{X})\boldsymbol{s}_{\boldsymbol{\Lambda}}(\boldsymbol{X}
)\right]  , \\
\boldsymbol{\theta}  &  =(\boldsymbol{\mu}^{\top},\boldsymbol{1}^{\top}
\boldsymbol{\Lambda})^{\top}.
\end{aligned}
\]

Since

\[
\begin{aligned}
\mathrm{E}\left[  f_{\theta}^{\beta}(\boldsymbol{X})\boldsymbol{h}
(\boldsymbol{X})\right]   &  =\int_{\mathbb{R}^{p}}f_{\theta}^{\beta+1}(\boldsymbol{x})\boldsymbol{h}(\boldsymbol{x}
)d\boldsymbol{x} \\
&  =\frac{(\beta+1)^{-\frac{p}{2}}}{(2\pi)^{\frac{\beta p}{2}}\left\vert
\boldsymbol{\Sigma}(\boldsymbol{\Lambda})\right\vert ^{\frac{\beta}{2}}}\int_{\mathbb{R}^{p}}f_{\mathcal{N}_{p}(\boldsymbol{\mu},\frac{1}{\beta+1}\boldsymbol{\Sigma
}(\boldsymbol{\Lambda}))}(\boldsymbol{x})\boldsymbol{h}(\boldsymbol{x}
)d\boldsymbol{x},
\end{aligned}
\]

it holds

\[
\begin{aligned}
\mathrm{E}\left[  f_{\theta}^{\beta}(\boldsymbol{X})(\boldsymbol{X}
-\boldsymbol{\mu})\right]   &  =\frac{(\beta+1)^{-\frac{p}{2}}}{(2\pi
)^{\frac{\beta p}{2}}\left\vert \boldsymbol{\Sigma}(\boldsymbol{\Lambda
})\right\vert ^{\frac{\beta}{2}}}\left(  \mathrm{E}[\mathcal{N}_{p}
(\boldsymbol{\mu},\tfrac{1}{\beta+1}\boldsymbol{\Sigma}(\boldsymbol{\Lambda
}))]-\boldsymbol{\mu}\right)  =\boldsymbol{0}_{p}, \\
\boldsymbol{\xi}_{\beta}(\boldsymbol{\mu})  &  =E\left[  f_{\theta}^{\beta
}(\boldsymbol{X})\boldsymbol{s}_{\boldsymbol{\mu}}(\boldsymbol{X})\right]
=\boldsymbol{\Lambda}^{-1/2}\boldsymbol{R}_{0}^{-1}\boldsymbol{\Lambda}
^{-1/2}\mathrm{E}\left[  f_{\theta}^{\beta}(\boldsymbol{X})(\boldsymbol{X}
-\boldsymbol{\mu})\right]  =\boldsymbol{0}_{p},
\end{aligned}
\]

and

\[
\begin{aligned}
\mathrm{E}\left[  f_{\theta}^{\beta}(\boldsymbol{X})(\boldsymbol{X}
-\boldsymbol{\mu})(\boldsymbol{X}-\boldsymbol{\mu})^{\top}\right]   &
=\frac{(\beta+1)^{-\frac{p}{2}}}{(2\pi)^{\frac{\beta p}{2}}\left\vert
\boldsymbol{\Sigma}(\boldsymbol{\Lambda})\right\vert ^{\frac{\beta}{2}}
}\mathrm{Var}[\mathcal{N}_{p}(\boldsymbol{\mu},\tfrac{1}{\beta+1}
\boldsymbol{\Sigma}(\boldsymbol{\Lambda}))] \\
&  =\frac{(\beta+1)^{-(\frac{p}{2}+1)}}{(2\pi)^{\frac{\beta p}{2}}\left\vert
\boldsymbol{\Sigma}(\boldsymbol{\Lambda})\right\vert ^{\frac{\beta}{2}}
}\boldsymbol{\Sigma}(\boldsymbol{\Lambda}),
\end{aligned}
\]

\[
\mathrm{E}\left[  f_{\theta}^{\beta}(\boldsymbol{X})\right]  =\frac
{(\beta+1)^{-\frac{p}{2}}}{(2\pi)^{\frac{\beta p}{2}}\left\vert
\boldsymbol{\Sigma}(\boldsymbol{\Lambda})\right\vert ^{\frac{\beta}{2}}},
\]

\[
\begin{aligned}
\boldsymbol{\xi}_{\beta}(\boldsymbol{\Lambda})  &  =E\left[  f_{\theta}
^{\beta}(\boldsymbol{X})\boldsymbol{s}_{\boldsymbol{\Lambda}}(\boldsymbol{X}
)\right] \\
&  =-\frac{1}{2}\left(  \boldsymbol{I}_{p}\mathrm{E}\left[  f_{\theta}^{\beta
}(\boldsymbol{X})\right]  -\boldsymbol{R}_{0}^{-1}\boldsymbol{\Lambda}
^{-1/2}\mathrm{E}\left[  f_{\theta}^{\beta}(\boldsymbol{X})(\boldsymbol{X}
-\boldsymbol{\mu})(\boldsymbol{X}-\boldsymbol{\mu})^{\top}\right]
\boldsymbol{\Lambda}^{-1/2}\right)  \boldsymbol{\Lambda}^{-1} \\
&  =-\frac{1}{2}\left(  \boldsymbol{I}_{p}\frac{(\beta+1)^{-\frac{p}{2}}
}{(2\pi)^{\frac{\beta p}{2}}\left\vert \boldsymbol{\Sigma}(\boldsymbol{\Lambda
})\right\vert ^{\frac{\beta}{2}}}-\frac{(\beta+1)^{-(\frac{p}{2}+1)}}
{(2\pi)^{\frac{\beta p}{2}}\left\vert \boldsymbol{\Sigma}(\boldsymbol{\Lambda
})\right\vert ^{\frac{\beta}{2}}}\boldsymbol{R}_{0}^{-1}\boldsymbol{\Lambda
}^{-1/2}\boldsymbol{\Sigma}(\boldsymbol{\Lambda})\boldsymbol{\Lambda}
^{-1/2}\right)  \boldsymbol{\Lambda}^{-1} \\
&  =-\frac{1}{2}\frac{(\beta+1)^{-\frac{p}{2}}}{(2\pi)^{\frac{\beta p}{2}
}\left\vert \boldsymbol{\Sigma}(\boldsymbol{\Lambda})\right\vert ^{\frac
{\beta}{2}}}\left(  \boldsymbol{I}_{p}-\tfrac{1}{\beta+1}\boldsymbol{R}
_{0}^{-1}\boldsymbol{R}_{0}\right)  \boldsymbol{\Lambda}^{-1} \\
&  =-\frac{1}{2}\frac{\beta(\beta+1)^{-(\frac{p}{2}+1)}}{(2\pi)^{\frac{\beta
p}{2}}\left\vert \boldsymbol{\Sigma}(\boldsymbol{\Lambda})\right\vert
^{\frac{\beta}{2}}}\boldsymbol{\Lambda}^{-1}.
\end{aligned}
\]

Hence,

\[
\begin{aligned}
-\frac{1}{2}\boldsymbol{\Sigma}^{-1}(\boldsymbol{\Lambda})\frac{1}{n}
\sum_{i=1}^{n}f_{\theta}^{\beta}(\boldsymbol{X}_{i})(\boldsymbol{X}
_{i}-\boldsymbol{\mu})  &  =\boldsymbol{0}_{p}, \\
-\frac{1}{2n}\sum_{i=1}^{n}f_{\theta}^{\beta}(\boldsymbol{X}_{i})\left[
\boldsymbol{I}_{p}-\boldsymbol{R}_{0}^{-1}\boldsymbol{\Lambda}^{-1/2}
(\boldsymbol{X}_{i}-\boldsymbol{\mu})(\boldsymbol{X}_{i}-\boldsymbol{\mu}
)^{\top}\boldsymbol{\Lambda}^{-1/2}\right]  \boldsymbol{\Lambda}^{-1}+\frac
{\beta}{2}\frac{(\beta+1)^{-(\frac{p}{2}+1)}}{(2\pi)^{\frac{\beta p}{2}
}\left\vert \boldsymbol{\Sigma}(\boldsymbol{\Lambda})\right\vert ^{\frac
{\beta}{2}}}\boldsymbol{\Lambda}^{-1}  &  =\boldsymbol{0}_{p\times p},
\end{aligned}
\]

equivalent to

\[
\begin{aligned}
\sum_{i=1}^{n}f_{\theta}^{\beta}(\boldsymbol{X}_{i})(\boldsymbol{X}
_{i}-\boldsymbol{\mu})  &  =\boldsymbol{0}_{p}, \\
\boldsymbol{R}_{0}^{-1}\boldsymbol{\Lambda}^{-1/2}\left(  \frac{1}{n}
\sum_{i=1}^{n}f_{\theta}^{\beta}(\boldsymbol{X}_{i})(\boldsymbol{X}
_{i}-\boldsymbol{\mu})(\boldsymbol{X}_{i}-\boldsymbol{\mu})^{\top}\right)
\boldsymbol{\Lambda}^{-1/2}  &  =\left(  \frac{1}{n}\sum_{i=1}^{n}f_{\theta
}^{\beta}(\boldsymbol{X}_{i})-\beta\frac{(\beta+1)^{-(\frac{p}{2}+1)}}
{(2\pi)^{\frac{\beta p}{2}}\left\vert \boldsymbol{\Sigma}(\boldsymbol{\Lambda
})\right\vert ^{\frac{\beta}{2}}}\right)  \boldsymbol{I}_{p},
\end{aligned}
\]

i.e.

\[
\begin{aligned}
\frac{\sum_{i=1}^{n}f_{\theta}^{\beta}(\boldsymbol{X}_{i})\boldsymbol{X}_{i}
}{\sum_{i=1}^{n}f_{\theta}^{\beta}(\boldsymbol{X}_{i})}  &  =\boldsymbol{\mu
}, \\
\boldsymbol{R}_{0}^{-1}\boldsymbol{\Lambda}^{-1/2}\frac{\frac{1}{n}\sum
_{i=1}^{n}f_{\theta}^{\beta}(\boldsymbol{X}_{i})(\boldsymbol{X}_{i}
-\boldsymbol{\mu})(\boldsymbol{X}_{i}-\boldsymbol{\mu})^{\top}}{\left(  \frac
{1}{n}\sum_{i=1}^{n}f_{\theta}^{\beta}(\boldsymbol{X}_{i})-\beta\frac
{(\beta+1)^{-(\frac{p}{2}+1)}}{(2\pi)^{\frac{\beta p}{2}}\left\vert
\boldsymbol{\Sigma}\right\vert ^{\frac{\beta}{2}}}\right)  }
\boldsymbol{\Lambda}^{-1/2}  &  =\boldsymbol{I}_{p}.
\end{aligned}
\]

From the previous expression it is obtained the system of equations
(\ref{eq-eqMuBeta})-(\ref{eq-eqRRBeta}).

\subsection{\texorpdfstring{Proof of Theorem
\ref{thm-RestrMDPDsEqui}}{Proof of Theorem }}\label{proof-of-theorem--thm-restrmdpdsequi}

From the particular case of the Woodbury's formula

\[
\left(  \boldsymbol{I}_{p}+\boldsymbol{UV}\right)  ^{-1}=\boldsymbol{I}
_{p}-\boldsymbol{U}(\boldsymbol{I}_{r}+\boldsymbol{VU})^{-1}\boldsymbol{V},
\]

it holds

\[
\begin{aligned}
\frac{1}{1-\rho_{0}}\boldsymbol{R}_{0}(\rho_{0})  &  =\boldsymbol{I}_{p}
+\frac{\rho_{0}}{1-\rho_{0}}\boldsymbol{1}_{p}\boldsymbol{1}_{p}^{\top} \\
(1-\rho_{0})\boldsymbol{R}_{0}^{-1}(\rho_{0})  &  =\left[  \boldsymbol{U}
=\frac{\rho_{0}}{1-\rho_{0}}\boldsymbol{1}_{p};\boldsymbol{V}=\boldsymbol{1}
_{p}^{\top}\right] \\
&  =\boldsymbol{I}_{p}-\frac{\rho_{0}}{1-\rho_{0}}\boldsymbol{1}_{p}
(1+p\tfrac{\rho_{0}}{1-\rho_{0}})^{-1}\boldsymbol{1}_{p}^{\top} \\
&  =\boldsymbol{I}_{p}-\frac{\rho_{0}}{1-\rho_{0}}\boldsymbol{1}_{p}
(\tfrac{1+(p-1)\rho_{0}}{1-\rho_{0}})^{-1}\boldsymbol{1}_{p}^{\top} \\
&  =\boldsymbol{I}_{p}-\frac{\rho_{0}}{1+(p-1)\rho_{0}}\boldsymbol{1}
_{p}\boldsymbol{1}_{p}^{\top},
\end{aligned}
\]

and hence

\begin{equation}\protect\phantomsection\label{eq-R-minus1}{
\begin{aligned}
\boldsymbol{R}_{0}^{-1}(\rho_{0})  &  =\frac{1}{1-\rho_{0}}\left(
\boldsymbol{I}_{p}-\frac{\rho_{0}}{1+(p-1)\rho_{0}}\boldsymbol{1}
_{p}\boldsymbol{1}_{p}^{\top}\right)  ,
\end{aligned}
}\end{equation}

\begin{equation}\protect\phantomsection\label{eq-R-1Rt}{
\begin{aligned}
\boldsymbol{R}_{0}^{-1}(\rho_{0})\widetilde{\boldsymbol{R}}_{\boldsymbol{X}
,\beta}  &  =\frac{1}{1-\rho_{0}}\left(  \widetilde{\boldsymbol{R}
}_{\boldsymbol{X},\beta}-\frac{\rho_{0}}{1+(p-1)\rho_{0}}\boldsymbol{1}
_{p}\boldsymbol{1}_{p}^{\top}\widetilde{\boldsymbol{R}}_{\boldsymbol{X},\beta
}\right) \\
&  =\frac{1}{1-\rho_{0}}\left(  \widetilde{\boldsymbol{R}}_{\boldsymbol{X}
,\beta}-\frac{\rho_{0}}{1+(p-1)\rho_{0}}\boldsymbol{1}_{p}\otimes
(\widetilde{R}_{\cdot1},\ldots,\widetilde{R}_{\cdot p})\right) \\
\mathrm{diag}\{\boldsymbol{R}_{0}^{-1}(\rho_{0})\widetilde{\boldsymbol{R}
}_{\boldsymbol{X},\beta}\}  &  =\frac{1}{1-\rho_{0}}\mathrm{diag}
\{\widetilde{R}_{jj,\beta}-\tfrac{\rho_{0}}{1+(p-1)\rho_{0}}\widetilde{R}
_{\cdot j,\beta}\}_{j=1}^{p}.
\end{aligned}
}\end{equation}

From (\ref{eq-R-1Rt}), taking into account (\ref{eq-eqRRBeta}), can be
constructed the estimating equations

\[
\widetilde{R}_{jj,\beta}-\frac{\rho_{0}}{1+(p-1)\rho_{0}}\widetilde{R}_{\cdot
j,\beta}=1-\rho_{0},
\]

or equivalently (\ref{eq-estimEqui}).

\subsection{\texorpdfstring{Proof of Theorem
\ref{thm-ThEstimRho}}{Proof of Theorem }}\label{proof-of-theorem--thm-thestimrho}

The estimating equations of \(\boldsymbol{\mu}\) and of the variances
are obtained as in the proof of Theorem~\ref{thm-RestrMDPDs}, with
\(\boldsymbol{R}_{0}\) replaced by \(\boldsymbol{R}(\rho_{12})\), which
gives (\ref{eq-eqMuBeta}) and (\ref{eq-rhoEstVar}). For \(\rho_{12}\),
the \(\beta\)-score in the direction
\(\frac{\partial}{\partial\rho_{12}}\mathrm{vec}(\boldsymbol{\Sigma})=\mathrm{vec}(\boldsymbol{\Lambda}^{1/2}\boldsymbol{E}\boldsymbol{\Lambda}^{1/2})\),
with
\(\boldsymbol{E}=\boldsymbol{1}_{p}\boldsymbol{1}_{p}^{\top}-\boldsymbol{I}_{p}\),
is, by (\ref{eq-V}) with \(\boldsymbol{R}_{0}\) replaced by
\(\boldsymbol{R}=\boldsymbol{R}(\rho_{12})\), proportional to \[
\begin{aligned}
&\mathrm{vec}^{\top}(\boldsymbol{\Lambda}^{1/2}\boldsymbol{E}\boldsymbol{\Lambda}^{1/2})(\boldsymbol{\Lambda}^{-1/2}\otimes\boldsymbol{\Lambda}^{-1/2})(\boldsymbol{R}^{-1}\otimes\boldsymbol{R}^{-1})\mathrm{vec}\left(\widetilde{\boldsymbol{R}}_{\boldsymbol{X},\beta}-\boldsymbol{R}\right)\\
&\qquad=\mathrm{vec}^{\top}(\boldsymbol{E})\mathrm{vec}\left(\boldsymbol{R}^{-1}(\widetilde{\boldsymbol{R}}_{\boldsymbol{X},\beta}-\boldsymbol{R})\boldsymbol{R}^{-1}\right)
=\mathrm{trace}(\boldsymbol{E}\boldsymbol{M}),
\end{aligned}
\] with
\(\boldsymbol{M}=\boldsymbol{R}^{-1}(\widetilde{\boldsymbol{R}}_{\boldsymbol{X},\beta}-\boldsymbol{R})\boldsymbol{R}^{-1}\),
and
\(\mathrm{trace}(\boldsymbol{E}\boldsymbol{M})=\boldsymbol{1}_{p}^{\top}\boldsymbol{M}\boldsymbol{1}_{p}-\mathrm{trace}(\boldsymbol{M})\),
which gives (\ref{eq-rhoEst}). Finally, by (\ref{eq-R-minus1}),
\(\boldsymbol{1}_{p}^{\top}\boldsymbol{R}^{-1}(\rho_{12})=\{1+(p-1)\rho_{12}\}^{-1}\boldsymbol{1}_{p}^{\top}\)
and
\(\boldsymbol{1}_{p}^{\top}\boldsymbol{R}(\rho_{12})\boldsymbol{1}_{p}=p\{1+(p-1)\rho_{12}\}\),
so that
\(\boldsymbol{1}_{p}^{\top}\boldsymbol{M}\boldsymbol{1}_{p}=\{1+(p-1)\rho_{12}\}^{-2}\{\boldsymbol{1}_{p}^{\top}\widetilde{\boldsymbol{R}}_{\boldsymbol{X},\beta}\boldsymbol{1}_{p}-p(1+(p-1)\rho_{12})\}\).

\subsection{\texorpdfstring{Proof of Theorem
\ref{thm-ThTestGeneral}}{Proof of Theorem }}\label{sec-ProofThTestGeneral}

By Lemma~\ref{lem-short}, the Rao's \(\beta\)-score test-statistic for
(\ref{eq-nullVSalt}) is
\begin{equation}\protect\phantomsection\label{eq-R-dif-par}{
\widetilde{R}_{\beta,n}=n\boldsymbol{U}_{\beta,n}^{\top}(\widetilde{\boldsymbol{\theta}}_{2,\beta})\boldsymbol{K}_{\beta}^{-1}(\widetilde{\boldsymbol{\theta}}_{2,\beta})\boldsymbol{U}_{\beta,n}(\widetilde{\boldsymbol{\theta}}_{2,\beta}),
}\end{equation}

where

\[
\boldsymbol{K}_{\beta}^{-1}(\widetilde{\boldsymbol{\theta}}_{2,\beta}
)=4(2\pi)^{\beta p}\left\vert \widetilde{\boldsymbol{\Sigma}}_{\beta
}\right\vert ^{\beta}\boldsymbol{G}_{p}^{+}\left[  \overline{\boldsymbol{J}
}_{2\beta}(\widetilde{\boldsymbol{\Sigma}}_{\beta}^{-1})+\overline
{\boldsymbol{\xi}}_{\beta}(\widetilde{\boldsymbol{\Sigma}}_{\beta}
^{-1})\overline{\boldsymbol{\xi}}_{\beta}^{\top}(\widetilde{\boldsymbol{\Sigma}
}_{\beta}^{-1})\right]  ^{-1}\boldsymbol{G}_{p}^{+\top},
\]

\(\boldsymbol{G}_{p}^{+}=(\boldsymbol{G}_{p}^{\top}\boldsymbol{G}_{p})^{-1}\boldsymbol{G}_{p}^{\top}\)
is the Moore-Penrose inverse of \(\boldsymbol{G}_{p}^{\top}\)
(\(\boldsymbol{G}_{p}^{+}\boldsymbol{G}_{p}=\boldsymbol{I}_{p(p+1)/2}\))
and
\(\widetilde{\boldsymbol{\theta}}_{2,\beta}=\mathrm{vech} (\widetilde{\boldsymbol{\Sigma}}_{\beta})\),
with
\(\widetilde{\boldsymbol{\Sigma}}_{\beta}\boldsymbol{=} \widetilde{\boldsymbol{\Lambda}}_{\beta}^{1/2}\boldsymbol{R}_{0} \widetilde{\boldsymbol{\Lambda}}_{\beta}^{1/2}\).
For the calculation of this inverse see Magnus and Neudecker (1980) and
Browne (1974) and according to Proposition \ref{prp-propK}

\[
\begin{aligned}
&  \left[  \overline{\boldsymbol{J}}_{2\beta}(\widetilde{\boldsymbol{\Sigma}
}_{\beta}^{-1})+\overline{\boldsymbol{\xi}}_{\beta}
(\widetilde{\boldsymbol{\Sigma}}_{\beta}^{-1})\overline{\boldsymbol{\xi}
}_{\beta}^{\top}(\widetilde{\boldsymbol{\Sigma}}_{\beta}^{-1})\right]
^{-1}=\kappa_{1}^{-1}(p,\beta)(\widetilde{\boldsymbol{\Lambda}}_{\beta}
^{1/2}\otimes\widetilde{\boldsymbol{\Lambda}}_{\beta}^{1/2})\left(
\boldsymbol{R}_{0}\otimes\boldsymbol{R}_{0}\right)
(\widetilde{\boldsymbol{\Lambda}}_{\beta}^{1/2}\otimes
\widetilde{\boldsymbol{\Lambda}}_{\beta}^{1/2}) \\
&  -\kappa_{1}^{-1}(p,\beta)(\widetilde{\boldsymbol{\Lambda}}_{\beta}
^{1/2}\otimes\widetilde{\boldsymbol{\Lambda}}_{\beta}^{1/2})\frac{\kappa
_{3}(p,\beta)\mathrm{vec}(\boldsymbol{R}_{0})\mathrm{vec}^{\top}(\boldsymbol{R}
_{0})}{1+p\kappa_{3}(p,\beta)}(\widetilde{\boldsymbol{\Lambda}}_{\beta}
^{1/2}\otimes\widetilde{\boldsymbol{\Lambda}}_{\beta}^{1/2}).
\end{aligned}
\]

From the expression of (\ref{eq-V}), it holds

\[
\begin{aligned}
\widetilde{R}_{\beta,n}  &  =n\frac{\widetilde{\kappa}_{0}^{2}(p,\beta
)}{\kappa_{1}(p,\beta)}\mathrm{vec}^{\top}\left(  \widetilde{\boldsymbol{R}
}_{\boldsymbol{X},\beta}-\boldsymbol{R}_{0}\right)  (\boldsymbol{R}_{0}
^{-1}\otimes\boldsymbol{R}_{0}^{-1})\mathrm{vec}\left(
\widetilde{\boldsymbol{R}}_{\boldsymbol{X},\beta}-\boldsymbol{R}_{0}\right) \\
&  -n\frac{\widetilde{\kappa}_{0}^{2}(p,\beta)}{\kappa_{1}(p,\beta
)}\mathrm{vec}^{\top}\left(  \widetilde{\boldsymbol{R}}_{\boldsymbol{X},\beta
}-\boldsymbol{R}_{0}\right)  \frac{\kappa_{3}(p,\beta)\mathrm{vec}
(\boldsymbol{R}_{0}^{-1})\mathrm{vec}^{\top}(\boldsymbol{R}_{0}^{-1})}
{1+p\kappa_{3}(p,\beta)}\mathrm{vec}\left(  \widetilde{\boldsymbol{R}
}_{\boldsymbol{X},\beta}-\boldsymbol{R}_{0}\right) \\
&  =n\frac{\widetilde{\kappa}_{0}^{2}(p,\beta)}{\kappa_{1}(p,\beta
)}\mathrm{vec}^{\top}\left(  \boldsymbol{R}_{0}^{-1/2}\widetilde{\boldsymbol{R}
}_{\boldsymbol{X},\beta}\boldsymbol{R}_{0}^{-1/2}-\boldsymbol{I}_{p}\right)
\mathrm{vec}\left(  \boldsymbol{R}_{0}^{-1/2}\widetilde{\boldsymbol{R}
}_{\boldsymbol{X},\beta}\boldsymbol{R}_{0}^{-1/2}-\boldsymbol{I}_{p}\right) \\
&  =n\frac{\widetilde{\kappa}_{0}^{2}(p,\beta)}{\kappa_{1}(p,\beta
)}\mathrm{trace}\left(  \left(  \boldsymbol{R}_{0}^{-1/2}
\widetilde{\boldsymbol{R}}_{\boldsymbol{X},\beta}\boldsymbol{R}_{0}
^{-1/2}-\boldsymbol{I}_{p}\right)  ^{2}\right) \\
&  =n\frac{\widetilde{\kappa}_{0}^{2}(p,\beta)}{\kappa_{1}(p,\beta
)}\mathrm{trace}\left(  \left(  \boldsymbol{R}_{0}^{-1}
\widetilde{\boldsymbol{R}}_{\boldsymbol{X},\beta}-\boldsymbol{I}_{p}\right)
^{2}\right)  .
\end{aligned}
\]

The expression in the second row vanishes since

\[
\begin{aligned}
\mathrm{vec}^{\top}(\boldsymbol{R}_{0}^{-1})\mathrm{vec}\left(
\widetilde{\boldsymbol{R}}_{\boldsymbol{X},\beta}-\boldsymbol{R}_{0}\right)
&  =\mathrm{vec}^{\top}(\boldsymbol{R}_{0}^{-1})\mathrm{vec}\left(
\widetilde{\boldsymbol{R}}_{\boldsymbol{X},\beta}\right)  -\mathrm{vec}
^{\top}(\boldsymbol{R}_{0}^{-1})\mathrm{vec}\left(  \boldsymbol{R}_{0}\right) \\
&  =\mathrm{trace}(\boldsymbol{R}_{0}^{-1}\widetilde{\boldsymbol{R}
}_{\boldsymbol{X},\beta})-\mathrm{trace}(\boldsymbol{I}_{p}) \\
&  =p-p \\
&  =0,
\end{aligned}
\]

and
\(\mathrm{trace}(\boldsymbol{R}_{0}^{-1}\widetilde{\boldsymbol{R} }_{\boldsymbol{X},\beta})=p\)
from the corresponding estimating equation.

\subsection{\texorpdfstring{Proof of Corollary
\ref{cor-ThTestEquiF}}{Proof of Corollary }}\label{sec-ProofThTestEquiF}

Let
\(\widetilde{\boldsymbol{D}}_{\beta}=\boldsymbol{R}_{0}^{-1}\widetilde{\boldsymbol{R}}_{\boldsymbol{X},\beta}-\boldsymbol{I}_{p}\).
By (\ref{eq-trace-vecl}),
\(\mathrm{trace}(\widetilde{\boldsymbol{D}}_{\beta}^{2})=2\sum_{i<j}\widetilde{D}_{ij,\beta}\widetilde{D}_{ji,\beta}\).
By (\ref{eq-R-minus1}), for \(i\neq j\),
\(\widetilde{D}_{ij,\beta}=(1-\rho_{0})^{-1}\{\widetilde{R}_{ij,\beta}-\frac{\rho_{0}}{1+(p-1)\rho_{0}}\widetilde{R}_{\cdot j,\beta}\}\),
and by (\ref{eq-estimEqui}) the second term in braces is
\(\widetilde{R}_{jj,\beta}-(1-\rho_{0})\), so that \[
\widetilde{D}_{ij,\beta}=\frac{\widetilde{R}_{ij,\beta}-\widetilde{R}_{jj,\beta}+1-\rho_{0}}{1-\rho_{0}},\qquad
\widetilde{D}_{ji,\beta}=\frac{\widetilde{R}_{ij,\beta}-\widetilde{R}_{ii,\beta}+1-\rho_{0}}{1-\rho_{0}},
\] and (\ref{eq-RaoR0f}) follows from (\ref{eq-RaoR0}).

\subsection{\texorpdfstring{Proof of Corollary
\ref{cor-ThTestInd}}{Proof of Corollary }}\label{proof-of-corollary--cor-thtestind}

The particularization for the uncorrelatedness or independence test with
respect to Theorem \ref{thm-ThTestGeneral}, where
\(\boldsymbol{R}_{0} =\boldsymbol{I}_{p}\), gives

\begin{equation}\protect\phantomsection\label{eq-RRInd}{
\begin{aligned}
\widetilde{R}_{\beta,n}  &  =n\frac{\widetilde{\kappa}_{0}^{2}(p,\beta
)}{\kappa_{1}(p,\beta)}\mathrm{trace}\left(  \left(  \widetilde{\boldsymbol{R}
}_{\boldsymbol{X},\beta}-\boldsymbol{I}_{p}\right)  ^{2}\right) \\
&  =n\frac{\widetilde{\kappa}_{0}^{2}(p,\beta)}{\kappa_{1}(p,\beta
)}\mathrm{vec}^{\top}\left(  \widetilde{\boldsymbol{R}}_{\boldsymbol{X},\beta
}-\boldsymbol{I}_{p}\right)  \mathrm{vec}\left(  \widetilde{\boldsymbol{R}
}_{\boldsymbol{X},\beta}-\boldsymbol{I}_{p}\right) \\
&  =2n\frac{\widetilde{\kappa}_{0}^{2}(p,\beta)}{\kappa_{1}(p,\beta
)}\mathrm{vecl}^{\top}(\widetilde{\boldsymbol{R}}_{\boldsymbol{X},\beta
})\mathrm{vecl}(\widetilde{\boldsymbol{R}}_{\boldsymbol{X},\beta}) \\
&  =2n\frac{\widetilde{\kappa}_{0}^{2}(p,\beta)}{\kappa_{1}(p,\beta)}
\sum_{i<j}\widetilde{R}_{ij,\beta}^{2}.
\end{aligned}
}\end{equation}

In addition, from Theorem \ref{thm-RestrMDPDsEqui} it is concluded that
\(S_{j,\beta}^{2}=\widetilde{\sigma}_{j,\beta}^{2}\), \(j=1,...,p\), and
hence \(\widetilde{R}_{ij,\beta}^{2}=R_{ij,\beta}^{2}\) for all pairs
such that \(i<j\).

\subsection{\texorpdfstring{Proof of Corollary
\ref{cor-ThTestEquiNF}}{Proof of Corollary }}\label{sec-ProofThTestEquiNF}

The parameter \(\rho_{12}\) is one of the nuisance parameters of
Lemma~\ref{lem-short}, so the Rao's \(\beta\)-score test-statistic is
the quadratic form of the proof of Theorem~\ref{thm-ThTestGeneral}
evaluated at the restricted estimators of Theorem~\ref{thm-ThEstimRho}.
That proof only uses that the components of the \(\beta\)-score function
corresponding to \(\boldsymbol{\mu}\) and to the variances are zero at
the restricted estimators, which holds by (\ref{eq-eqMuBeta}) and
(\ref{eq-rhoEstVar}), so the same computation with
\(\boldsymbol{R}_{0}\) replaced by
\(\boldsymbol{R}(\widetilde{\rho}_{12,\beta})\) gives (\ref{eq-RaoR0g}).
The number of degrees of freedom is the number of correlation
parameters, \(\frac{p(p-1)}{2}\), minus the dimension of the correlation
parameter under the null hypothesis, one.

\subsection{\texorpdfstring{Auxiliary results for the proofs of Section
\ref{sec-Asym}}{Auxiliary results for the proofs of Section }}\label{sec-ProofAux}

We use the notation of the proof of Lemma~\ref{lem-short}, with all the
matrices evaluated at \(\boldsymbol{\theta}_{0}\),
\(\boldsymbol{A}=\boldsymbol{G}_{p}^{\top}(\boldsymbol{\Sigma}_{0}^{-1}\otimes\boldsymbol{\Sigma}_{0}^{-1})\boldsymbol{G}_{p}\),
\(\boldsymbol{s}=\mathrm{vech}(\boldsymbol{\Sigma}_{0})\), and
\(\boldsymbol{N}\) and \(\boldsymbol{T}\) the Jacobian matrices of
\(\mathrm{vech}(\boldsymbol{\Sigma})\) with respect to
\(\boldsymbol{\eta}_{(1)}\) and to the remaining correlation parameters.
By (\ref{eq-eqJ22}) and Corollary~\ref{cor-Th-Psi-KS1}, the scalars of
that proof are \begin{equation}\protect\phantomsection\label{eq-ab}{
a_{\beta}=\frac{(\beta+1)^{-\frac{p}{2}-2}}{4(2\pi)^{\frac{\beta p}{2}}\left\vert \boldsymbol{\Sigma}_{0}\right\vert ^{\frac{\beta}{2}}},
\qquad
b_{\beta}=\frac{1}{4(2\pi)^{\beta p}\left\vert \boldsymbol{\Sigma}_{0}\right\vert ^{\beta}},
\qquad\text{so that}\qquad
\frac{4a_{\beta}^{2}}{b_{\beta}\kappa_{1}(p,\beta)}=\frac{1}{2}\mathrm{ARE}(\beta,p).
}\end{equation} Let also
\(\boldsymbol{A}_{\perp}=\boldsymbol{A}-\boldsymbol{A}\boldsymbol{N}(\boldsymbol{N}^{\top}\boldsymbol{A}\boldsymbol{N})^{-1}\boldsymbol{N}^{\top}\boldsymbol{A}\),
and, for a symmetric matrix \(\boldsymbol{M}\) of order \(p\),
\(\boldsymbol{h}(\boldsymbol{M})=\mathrm{vech}(\boldsymbol{\Lambda}_{0}^{1/2}\boldsymbol{M}\boldsymbol{\Lambda}_{0}^{1/2})\).

\begin{lemma}[]\protect\hypertarget{lem-perp}{}\label{lem-perp}

The following statements hold.

\begin{enumerate}
\def\labelenumi{(\alph{enumi})}
\item
  The Schur complements of the blocks of the nuisance parameters are \[
  \begin{aligned}
  \boldsymbol{J}_{\beta,22}(\boldsymbol{\eta}_{0})-\boldsymbol{J}_{\beta,21}(\boldsymbol{\eta}_{0})\boldsymbol{J}_{\beta,11}^{-1}(\boldsymbol{\eta}_{0})\boldsymbol{J}_{\beta,12}(\boldsymbol{\eta}_{0})&=2a_{\beta}\boldsymbol{T}^{\top}\boldsymbol{A}_{\perp}\boldsymbol{T},\\
  \boldsymbol{K}_{\beta,22}(\boldsymbol{\eta}_{0})-\boldsymbol{K}_{\beta,21}(\boldsymbol{\eta}_{0})\boldsymbol{K}_{\beta,11}^{-1}(\boldsymbol{\eta}_{0})\boldsymbol{K}_{\beta,12}(\boldsymbol{\eta}_{0})&=b_{\beta}\kappa_{1}(p,\beta)\boldsymbol{T}^{\top}\boldsymbol{A}_{\perp}\boldsymbol{T}.
  \end{aligned}
  \]
\item
  \(\boldsymbol{A}_{\perp}\boldsymbol{T}(\boldsymbol{T}^{\top}\boldsymbol{A}_{\perp}\boldsymbol{T})^{-1}\boldsymbol{T}^{\top}\boldsymbol{A}_{\perp}=\boldsymbol{A}_{\perp}\).
\item
  For every symmetric matrix \(\boldsymbol{M}\) of order \(p\) the
  matrix \(\boldsymbol{M}^{\ast}\) of (\ref{eq-adjust}) exists and is
  unique,
  \(\boldsymbol{h}(\boldsymbol{M}^{\ast})=\boldsymbol{h}(\boldsymbol{M})-\boldsymbol{A}\boldsymbol{N}(\boldsymbol{N}^{\top}\boldsymbol{A}\boldsymbol{N})^{-1}\boldsymbol{N}^{\top}\boldsymbol{A}\boldsymbol{h}(\boldsymbol{M})\)
  up to the factor \(\boldsymbol{A}\), and, for symmetric matrices
  \(\boldsymbol{M}_{1}\) and \(\boldsymbol{M}_{2}\), \[
  \boldsymbol{h}^{\top}(\boldsymbol{M}_{1})\boldsymbol{A}_{\perp}\boldsymbol{h}(\boldsymbol{M}_{2})
  =\mathrm{trace}\left(\boldsymbol{R}_{0}^{-1}\boldsymbol{M}_{1}^{\ast}\boldsymbol{R}_{0}^{-1}\boldsymbol{M}_{2}^{\ast}\right)
  =\mathrm{trace}\left(\boldsymbol{R}_{0}^{-1}\boldsymbol{M}_{1}^{\ast}\boldsymbol{R}_{0}^{-1}\boldsymbol{M}_{2}\right).
  \] In particular
  \(\mathrm{trace}((\boldsymbol{R}_{0}^{-1}\boldsymbol{M}^{\ast})^{2})\geq0\),
  with equality if and only if
  \(\boldsymbol{M}^{\ast}=\boldsymbol{0}_{p\times p}\).
\end{enumerate}

\end{lemma}

Proof of Lemma~\ref{lem-perp}. (a) Let
\(\boldsymbol{M}(\alpha_{1},\alpha_{2})=\alpha_{1}\boldsymbol{A}+\alpha_{2}\boldsymbol{A}\boldsymbol{s}\boldsymbol{s}^{\top}\boldsymbol{A}\).
Since \(\boldsymbol{s}=\boldsymbol{N}\boldsymbol{c}\), as in the proof
of Lemma~\ref{lem-short},
\(\boldsymbol{T}^{\top}\boldsymbol{M}(\alpha_{1},\alpha_{2})\boldsymbol{N}=\boldsymbol{A}_{21}\boldsymbol{A}_{11}^{-1}\boldsymbol{N}^{\top}\boldsymbol{M}(\alpha_{1},\alpha_{2})\boldsymbol{N}\),
and then \[
\begin{aligned}
&\boldsymbol{T}^{\top}\boldsymbol{M}(\alpha_{1},\alpha_{2})\boldsymbol{T}-\boldsymbol{A}_{21}\boldsymbol{A}_{11}^{-1}\boldsymbol{N}^{\top}\boldsymbol{M}(\alpha_{1},\alpha_{2})\boldsymbol{T}\\
&\qquad=\alpha_{1}(\boldsymbol{A}_{22}-\boldsymbol{A}_{21}\boldsymbol{A}_{11}^{-1}\boldsymbol{A}_{12})+\alpha_{2}(\boldsymbol{A}_{21}\boldsymbol{c}\boldsymbol{c}^{\top}\boldsymbol{A}_{12}-\boldsymbol{A}_{21}\boldsymbol{c}\boldsymbol{c}^{\top}\boldsymbol{A}_{12})
=\alpha_{1}\boldsymbol{T}^{\top}\boldsymbol{A}_{\perp}\boldsymbol{T},
\end{aligned}
\] with
\(\boldsymbol{A}_{22}=\boldsymbol{T}^{\top}\boldsymbol{A}\boldsymbol{T}\)
and \(\boldsymbol{A}_{12}=\boldsymbol{A}_{21}^{\top}\). The statement
follows with \((\alpha_{1},\alpha_{2})\) equal to
\(a_{\beta}(2,\beta^{2})\) and to
\(b_{\beta}(\kappa_{1}(p,\beta),\kappa_{2}(p,\beta))\), the blocks which
correspond to \(\boldsymbol{\mu}\) playing no role because
\(\boldsymbol{J}_{\beta}(\boldsymbol{\theta}_{0})\) and
\(\boldsymbol{K}_{\beta}(\boldsymbol{\theta}_{0})\) are block diagonal.

\begin{enumerate}
\def\labelenumi{(\alph{enumi})}
\setcounter{enumi}{1}
\item
  The matrix \((\boldsymbol{N},\boldsymbol{T})\) is nonsingular, so
  every vector is
  \(\boldsymbol{h}=\boldsymbol{N}\boldsymbol{d}_{1}+\boldsymbol{T}\boldsymbol{d}_{2}\),
  and \(\boldsymbol{A}_{\perp}\boldsymbol{N}=\boldsymbol{0}\) gives
  \(\boldsymbol{A}_{\perp}\boldsymbol{h}=\boldsymbol{A}_{\perp}\boldsymbol{T}\boldsymbol{d}_{2}\)
  and
  \(\boldsymbol{A}_{\perp}\boldsymbol{T}(\boldsymbol{T}^{\top}\boldsymbol{A}_{\perp}\boldsymbol{T})^{-1}\boldsymbol{T}^{\top}\boldsymbol{A}_{\perp}\boldsymbol{h}=\boldsymbol{A}_{\perp}\boldsymbol{T}\boldsymbol{d}_{2}\).
\item
  Since
  \(\boldsymbol{G}_{p}\mathrm{vech}(\boldsymbol{H})=\mathrm{vec}(\boldsymbol{H})\)
  for a symmetric matrix \(\boldsymbol{H}\),
  \begin{equation}\protect\phantomsection\label{eq-hAh}{
  \boldsymbol{h}^{\top}(\boldsymbol{M}_{1})\boldsymbol{A}\boldsymbol{h}(\boldsymbol{M}_{2})=\mathrm{trace}(\boldsymbol{\Sigma}_{0}^{-1}\boldsymbol{\Lambda}_{0}^{1/2}\boldsymbol{M}_{1}\boldsymbol{\Lambda}_{0}^{1/2}\boldsymbol{\Sigma}_{0}^{-1}\boldsymbol{\Lambda}_{0}^{1/2}\boldsymbol{M}_{2}\boldsymbol{\Lambda}_{0}^{1/2})=\mathrm{trace}(\boldsymbol{R}_{0}^{-1}\boldsymbol{M}_{1}\boldsymbol{R}_{0}^{-1}\boldsymbol{M}_{2}).
  }\end{equation} The column of \(\boldsymbol{N}\) which corresponds to
  \(\sigma_{j}^{2}\) is
  \(\boldsymbol{h}(\frac{1}{2\sigma_{j,0}^{2}}(\boldsymbol{e}_{j}\boldsymbol{e}_{j}^{\top}\boldsymbol{R}_{0}+\boldsymbol{R}_{0}\boldsymbol{e}_{j}\boldsymbol{e}_{j}^{\top}))\),
  with \(\boldsymbol{e}_{j}\) the \(j\)-th unit vector, and, for the
  test of Section \ref{sec-test4}, \(\boldsymbol{N}\) has the additional
  column
  \(\boldsymbol{h}(\boldsymbol{1}_{p}\boldsymbol{1}_{p}^{\top}-\boldsymbol{I}_{p})\),
  because
  \(\frac{\partial}{\partial\rho_{12}}\boldsymbol{\Sigma}=\boldsymbol{\Lambda}^{1/2}(\boldsymbol{1}_{p}\boldsymbol{1}_{p}^{\top}-\boldsymbol{I}_{p})\boldsymbol{\Lambda}^{1/2}\).
  Therefore the matrices
  \(\boldsymbol{D}\boldsymbol{R}_{0}+\boldsymbol{R}_{0}\boldsymbol{D}\)
  with \(\boldsymbol{D}\) diagonal, plus the multiples of
  \(\boldsymbol{1}_{p}\boldsymbol{1}_{p}^{\top}-\boldsymbol{I}_{p}\) in
  the second case, are exactly the symmetric matrices \(\boldsymbol{M}\)
  with \(\boldsymbol{h}(\boldsymbol{M})\) in the column space of
  \(\boldsymbol{N}\), and by
\end{enumerate}

\begin{enumerate}
\def\labelenumi{(\arabic{enumi})}
\tightlist
\item
  the conditions (\ref{eq-adjustvar}) and (\ref{eq-adjustrho}) on
  \(\boldsymbol{M}^{\ast}\) are
  \(\boldsymbol{N}^{\top}\boldsymbol{A}\boldsymbol{h}(\boldsymbol{M}^{\ast})=\boldsymbol{0}\).
  Since \(\boldsymbol{A}\) is positive definite and \(\boldsymbol{N}\)
  has full column rank, \(\boldsymbol{h}(\boldsymbol{M}^{\ast})\) is the
  residual of the projection of \(\boldsymbol{h}(\boldsymbol{M})\) on
  the column space of \(\boldsymbol{N}\) in the inner product defined by
  \(\boldsymbol{A}\), which exists and is unique, and
  \(\boldsymbol{h}^{\top}(\boldsymbol{M}_{1})\boldsymbol{A}_{\perp}\boldsymbol{h}(\boldsymbol{M}_{2})=\boldsymbol{h}^{\top}(\boldsymbol{M}_{1}^{\ast})\boldsymbol{A}\boldsymbol{h}(\boldsymbol{M}_{2}^{\ast})=\boldsymbol{h}^{\top}(\boldsymbol{M}_{1}^{\ast})\boldsymbol{A}\boldsymbol{h}(\boldsymbol{M}_{2})\),
  which is the statement by (1). Finally
  \(\mathrm{trace}((\boldsymbol{R}_{0}^{-1}\boldsymbol{M}^{\ast})^{2})=\mathrm{trace}((\boldsymbol{R}_{0}^{-1/2}\boldsymbol{M}^{\ast}\boldsymbol{R}_{0}^{-1/2})^{2})\geq0\),
  with equality if and only if
  \(\boldsymbol{R}_{0}^{-1/2}\boldsymbol{M}^{\ast}\boldsymbol{R}_{0}^{-1/2}=\boldsymbol{0}_{p\times p}\),
  since this matrix is symmetric.
\end{enumerate}

\subsection{\texorpdfstring{Proof of Theorem
\ref{thm-root}}{Proof of Theorem }}\label{sec-ProofRoot}

Under the null hypothesis the model is the parametric family indexed by
\(\boldsymbol{\eta}_{(1)}\), with the correlation parameters fixed,
which contains the true distribution, and its minimum DPD estimating
equations are the restricted ones. This family verifies the regularity
conditions (D1)-(D5) of Basu et al.~(2011, Section 9.2.1): the support
does not depend on the parameter; the parameter space is open and the
density is three times continuously differentiable with respect to
\(\boldsymbol{\eta}_{(1)}\); the derivatives of
\(f_{\boldsymbol{\eta}}^{\beta}(\boldsymbol{x})\) up to the third order
are the product of \(f_{\boldsymbol{\eta}}^{\beta}(\boldsymbol{x})\) and
a polynomial in \(\boldsymbol{x}\) whose coefficients are continuous
functions of the parameter, so that on a compact neighbourhood of
\(\boldsymbol{\eta}_{(1),0}\) they are bounded by a multiple of
\(1+\left\Vert \boldsymbol{x}\right\Vert ^{6}\), which is integrable,
and the integral
\(\int_{\mathbb{R}^{p}}f_{\boldsymbol{\eta}}^{1+\beta}(\boldsymbol{x})d\boldsymbol{x}\)
can be differentiated under the integral sign; and
\(\boldsymbol{J}_{\beta,11}(\boldsymbol{\eta}_{0})\) is positive
definite, because \(\boldsymbol{J}_{\beta}(\boldsymbol{\theta}_{0})\) is
positive definite by (\ref{eq-eqJ11}) and (\ref{eq-eqJ22}) and the
Jacobian matrix \(\boldsymbol{N}\) has full column rank. Statements (i)
and (ii) are then Theorem 9.2 of Basu et al.~(2011), which is Theorem 2
of Basu et al.~(1998), applied to this family, (\ref{eq-linear}) being
the asymptotic representation obtained in its proof. For (iii), the
matrix of second derivatives of
\(H_{\beta,n}(\boldsymbol{\eta}_{(1)},\mathrm{vecl}(\boldsymbol{R}_{0}))\)
with respect to \(\boldsymbol{\eta}_{(1)}\) is a mean of independent and
identically distributed terms, continuous in \(\boldsymbol{\eta}_{(1)}\)
and dominated as above, so it converges in probability, uniformly on a
compact neighbourhood of \(\boldsymbol{\eta}_{(1),0}\), to its
expectation, which is continuous in \(\boldsymbol{\eta}_{(1)}\) and
equal to \((\beta+1)\boldsymbol{J}_{\beta,11}(\boldsymbol{\eta}_{0})\)
at \(\boldsymbol{\eta}_{(1),0}\). Hence there is a ball centred at
\(\boldsymbol{\eta}_{(1),0}\) on which, with probability tending to one,
this matrix is positive definite, the function is strictly convex and
has at most one stationary point, and the consistent root of (i) belongs
to the ball with probability tending to one. For (iv), the function
(\ref{eq-Hn}) is continuous in the parameter and, on \(C\), it is a mean
of independent and identically distributed terms dominated by an
integrable function, bounded for \(\beta>0\) because
\(0<W_{i,\beta}(\boldsymbol{\theta})\leq1\), so it converges in
probability, uniformly on \(C\), to its expectation, which is the
density power divergence between the true distribution and the model up
to a term which does not depend on the parameter, and whose unique
minimizer is \(\boldsymbol{\eta}_{(1),0}\) because the family is
identifiable (Basu et al., 1998, Theorem 1). Hence the minimizer of
(\ref{eq-Hn}) on \(C\) converges in probability to
\(\boldsymbol{\eta}_{(1),0}\), it is an interior point of \(C\) and
therefore a solution of the estimating equations with probability
tending to one, no other solution in \(C\) has a smaller value of
(\ref{eq-Hn}), and by (iii) it is the root of (i). The last statement
holds because the proofs of Theorem~\ref{thm-ThTestGeneral} and of its
corollaries only use the consistency of the restricted estimators and
(\ref{eq-linear}).

\subsection{\texorpdfstring{Proof of Theorem
\ref{thm-contig}}{Proof of Theorem }}\label{sec-ProofContig}

Let \(\boldsymbol{h}=\boldsymbol{h}(\boldsymbol{\Delta})\), so that
\(\mathrm{vech}(\boldsymbol{\Sigma})=\boldsymbol{s}+n^{-1/2}\boldsymbol{h}\)
under (\ref{eq-H1n}). The multivariate normal family is differentiable
in quadratic mean, so by Theorem 7.2 of van der Vaart (1998) the
logarithm of the likelihood ratio is asymptotically normal,
(\ref{eq-H1n}) is contiguous to the null hypothesis and, under the null
hypothesis,
\(\sqrt{n}\boldsymbol{U}_{\beta,n}(\mathrm{vech}(\boldsymbol{\Sigma}_{0}))\)
and the logarithm of the likelihood ratio are jointly asymptotically
normal, with asymptotic covariance
\(\mathrm{E}_{\boldsymbol{\theta}_{0}}[\boldsymbol{u}_{\beta}(\boldsymbol{X},\boldsymbol{\theta}_{0})\boldsymbol{s}_{\boldsymbol{\theta}_{0}}^{\top}(\boldsymbol{X})]\)
times the direction, which is
\(\boldsymbol{J}_{\beta}(\mathrm{vech}(\boldsymbol{\Sigma}_{0}))\boldsymbol{h}\)
because
\(\mathrm{E}_{\boldsymbol{\theta}_{0}}[\boldsymbol{s}_{\boldsymbol{\theta}_{0}}(\boldsymbol{X})]=\boldsymbol{0}\).
By Le Cam's third lemma (van der Vaart, 1998, Example 6.7), under
(\ref{eq-H1n}), \begin{equation}\protect\phantomsection\label{eq-LeCam}{
\sqrt{n}\boldsymbol{U}_{\beta,n}(\mathrm{vech}(\boldsymbol{\Sigma}_{0}))\underset{n\rightarrow\infty}{\overset{\mathcal{L}}{\longrightarrow}}\mathcal{N}\left(\boldsymbol{J}_{\beta}(\mathrm{vech}(\boldsymbol{\Sigma}_{0}))\boldsymbol{h},\boldsymbol{K}_{\beta}(\mathrm{vech}(\boldsymbol{\Sigma}_{0}))\right),
}\end{equation} and by contiguity the consistency of the restricted
estimators and (\ref{eq-linear}) also hold under (\ref{eq-H1n}). A
Taylor expansion of \(\boldsymbol{U}_{\beta,n,2}\) around
\(\boldsymbol{\eta}_{0}\), (\ref{eq-linear}) and the equality
\(\boldsymbol{J}_{\beta,21}(\boldsymbol{\eta}_{0})\boldsymbol{J}_{\beta,11}^{-1}(\boldsymbol{\eta}_{0})=\boldsymbol{A}_{21}\boldsymbol{A}_{11}^{-1}\)
of the proof of Lemma~\ref{lem-short} give
\begin{equation}\protect\phantomsection\label{eq-projscore}{
\sqrt{n}\boldsymbol{U}_{\beta,n,2}(\widetilde{\boldsymbol{\eta}}_{\beta})=\boldsymbol{T}^{\top}\left(\boldsymbol{I}-\boldsymbol{A}\boldsymbol{N}(\boldsymbol{N}^{\top}\boldsymbol{A}\boldsymbol{N})^{-1}\boldsymbol{N}^{\top}\right)\sqrt{n}\boldsymbol{U}_{\beta,n}(\mathrm{vech}(\boldsymbol{\Sigma}_{0}))+\boldsymbol{r}_{n},
\qquad\boldsymbol{r}_{n}\underset{n\rightarrow\infty}{\overset{P}{\longrightarrow}}\boldsymbol{0}.
}\end{equation} Since
\(\boldsymbol{J}_{\beta}(\mathrm{vech}(\boldsymbol{\Sigma}_{0}))\boldsymbol{h}=a_{\beta}(2\boldsymbol{A}\boldsymbol{h}+\beta^{2}\boldsymbol{A}\boldsymbol{N}\boldsymbol{c}\boldsymbol{s}^{\top}\boldsymbol{A}\boldsymbol{h})\)
and
\((\boldsymbol{I}-\boldsymbol{A}\boldsymbol{N}(\boldsymbol{N}^{\top}\boldsymbol{A}\boldsymbol{N})^{-1}\boldsymbol{N}^{\top})\boldsymbol{A}\boldsymbol{N}=\boldsymbol{0}\),
by (\ref{eq-LeCam}), (\ref{eq-projscore}) and Lemma~\ref{lem-perp} (a),
under (\ref{eq-H1n}), \[
\sqrt{n}\boldsymbol{U}_{\beta,n,2}(\widetilde{\boldsymbol{\eta}}_{\beta})\underset{n\rightarrow\infty}{\overset{\mathcal{L}}{\longrightarrow}}\mathcal{N}\left(2a_{\beta}\boldsymbol{T}^{\top}\boldsymbol{A}_{\perp}\boldsymbol{h},b_{\beta}\kappa_{1}(p,\beta)\boldsymbol{T}^{\top}\boldsymbol{A}_{\perp}\boldsymbol{T}\right).
\] By Lemma~\ref{lem-short}, \(\widetilde{R}_{\beta,n}\) is the
quadratic form of
\(\sqrt{n}\boldsymbol{U}_{\beta,n,2}(\widetilde{\boldsymbol{\eta}}_{\beta})\)
in the inverse of a consistent estimator of its asymptotic
variance-covariance matrix, also consistent under (\ref{eq-H1n}) by
contiguity, so it converges in distribution to a
\(\chi_{r}^{2}(\delta)\) with \[
\delta=\frac{4a_{\beta}^{2}}{b_{\beta}\kappa_{1}(p,\beta)}\boldsymbol{h}^{\top}\boldsymbol{A}_{\perp}\boldsymbol{T}(\boldsymbol{T}^{\top}\boldsymbol{A}_{\perp}\boldsymbol{T})^{-1}\boldsymbol{T}^{\top}\boldsymbol{A}_{\perp}\boldsymbol{h}
=\frac{1}{2}\mathrm{ARE}(\beta,p)\,\mathrm{trace}\left(\left(\boldsymbol{R}_{0}^{-1}\boldsymbol{\Delta}^{\ast}\right)^{2}\right),
\] by (\ref{eq-ab}) and Lemma~\ref{lem-perp} (b) and (c). For
\(\boldsymbol{R}_{0}=\boldsymbol{I}_{p}\) the diagonal elements of
\(\boldsymbol{R}_{0}^{-1}\boldsymbol{\Delta}=\boldsymbol{\Delta}\) are
zero, so \(\boldsymbol{\Delta}^{\ast}=\boldsymbol{\Delta}\) and
\(\mathrm{trace}(\boldsymbol{\Delta}^{2})=2\sum_{i<j}\Delta_{ij}^{2}\).

\subsection{\texorpdfstring{Proof of Corollary
\ref{cor-are}}{Proof of Corollary }}\label{sec-ProofARE}

Two tests whose test-statistics converge to
\(\chi_{r}^{2}(\delta_{\beta}(\boldsymbol{\Delta}))\) and
\(\chi_{r}^{2}(\delta_{0}(\boldsymbol{\Delta}))\) under (\ref{eq-H1n})
have the same asymptotic power when the sample sizes \(n_{\beta}\) and
\(n_{0}\) verify
\(n_{0}/n_{\beta}\rightarrow\delta_{\beta}(\boldsymbol{\Delta})/\delta_{0}(\boldsymbol{\Delta})=\mathrm{ARE}(\beta,p)\),
since the noncentrality parameter is proportional to the sample size for
a fixed alternative. Moreover,
\(\mathrm{ARE}(\beta,p)=\{(1+2\beta)(1+\beta)^{-2}\}^{\frac{p}{2}+2}\)
and \((1+2\beta)(1+\beta)^{-2}\) is equal to one for \(\beta=0\) and
strictly decreasing in \(\beta\geq0\).

\subsection{\texorpdfstring{Proof of Theorem
\ref{thm-IFest}}{Proof of Theorem }}\label{sec-ProofIFest}

Differentiating
\(\mathrm{E}_{F_{\epsilon,\boldsymbol{x}}}[\boldsymbol{u}_{\beta,1}(\boldsymbol{X},(\widetilde{\boldsymbol{T}}_{\beta}(F_{\epsilon,\boldsymbol{x}}),\mathrm{vecl}(\boldsymbol{R}_{0})))]=\boldsymbol{0}\)
with respect to \(\epsilon\) at \(\epsilon=0\), and taking into account
that
\(\mathrm{E}_{\boldsymbol{\theta}_{0}}[\frac{\partial}{\partial\boldsymbol{\eta}_{(1)}^{\top}}\boldsymbol{u}_{\beta,1}(\boldsymbol{X},\boldsymbol{\eta}_{0})]=-\boldsymbol{J}_{\beta,11}(\boldsymbol{\eta}_{0})\),
we obtain (\ref{eq-IFgeneral}). For \(\boldsymbol{\mu}\), by
Theorem~\ref{thm-Th-U-BS1} and (\ref{eq-eqJ11}),
\(\boldsymbol{J}_{\beta}^{-1}(\boldsymbol{\mu}_{0})\boldsymbol{u}_{\beta}(\boldsymbol{x},\boldsymbol{\mu}_{0})=(\beta+1)^{\frac{p}{2}+1}w_{\beta}(\boldsymbol{x})(\boldsymbol{x}-\boldsymbol{\mu}_{0})\).
For the variances, by (\ref{eq-vecU}), the \(\beta\)-score function of
\(\mathrm{vech}(\boldsymbol{\Sigma})\) at \(\boldsymbol{\theta}_{0}\) is
\begin{equation}\protect\phantomsection\label{eq-uPsi}{
\boldsymbol{u}_{\beta}(\boldsymbol{x},\mathrm{vech}(\boldsymbol{\Sigma}_{0}))=\frac{1}{2(2\pi)^{\frac{\beta p}{2}}\left\vert \boldsymbol{\Sigma}_{0}\right\vert ^{\frac{\beta}{2}}}\boldsymbol{A}\boldsymbol{h}(\boldsymbol{\Psi}_{\beta}(\boldsymbol{x})),
}\end{equation} and, by (1) and the expression of the columns of
\(\boldsymbol{N}\) in the proof of Lemma~\ref{lem-perp},
\(\boldsymbol{N}^{\top}\boldsymbol{A}\boldsymbol{h}(\boldsymbol{M})=\boldsymbol{\Lambda}_{0}^{-1}\mathrm{diag}\{\boldsymbol{R}_{0}^{-1}\boldsymbol{M}\}\boldsymbol{1}_{p}\)
and
\(\boldsymbol{N}^{\top}\boldsymbol{A}\boldsymbol{N}\boldsymbol{\Lambda}_{0}\boldsymbol{d}(\boldsymbol{M})=\frac{1}{2}\boldsymbol{\Lambda}_{0}^{-1}\mathrm{diag}\{\boldsymbol{R}_{0}^{-1}\boldsymbol{M}\}\boldsymbol{1}_{p}\),
the last equality because
\(\boldsymbol{N}^{\top}\boldsymbol{A}\boldsymbol{h}(\boldsymbol{M}^{\ast})=\boldsymbol{0}\)
and
\(\boldsymbol{h}(\boldsymbol{D}(\boldsymbol{M})\boldsymbol{R}_{0}+\boldsymbol{R}_{0}\boldsymbol{D}(\boldsymbol{M}))=2\boldsymbol{N}\boldsymbol{\Lambda}_{0}\boldsymbol{d}(\boldsymbol{M})\).
With \(\boldsymbol{c}=\boldsymbol{\Lambda}_{0}\boldsymbol{1}_{p}\) and
\(\boldsymbol{c}^{\top}\boldsymbol{N}^{\top}\boldsymbol{A}\boldsymbol{N}\boldsymbol{c}=\boldsymbol{s}^{\top}\boldsymbol{A}\boldsymbol{s}=p\),
the Sherman-Morrison formula gives \[
\begin{aligned}
\boldsymbol{J}_{\beta,11}^{-1}(\boldsymbol{\Lambda}_{0}\boldsymbol{1}_{p})\boldsymbol{N}^{\top}\boldsymbol{A}\boldsymbol{h}(\boldsymbol{M})
&=\frac{1}{a_{\beta}}\left\{2\boldsymbol{N}^{\top}\boldsymbol{A}\boldsymbol{N}+\beta^{2}\boldsymbol{N}^{\top}\boldsymbol{A}\boldsymbol{N}\boldsymbol{c}\boldsymbol{c}^{\top}\boldsymbol{N}^{\top}\boldsymbol{A}\boldsymbol{N}\right\}^{-1}\boldsymbol{N}^{\top}\boldsymbol{A}\boldsymbol{h}(\boldsymbol{M})\\
&=\frac{2}{a_{\beta}}\boldsymbol{\Lambda}_{0}\left[\boldsymbol{d}(\boldsymbol{M})-\frac{\beta^{2}\,\mathrm{trace}(\boldsymbol{R}_{0}^{-1}\boldsymbol{M})}{2(2+\beta^{2}p)}\boldsymbol{1}_{p}\right],
\end{aligned}
\] because
\(\boldsymbol{c}^{\top}\boldsymbol{N}^{\top}\boldsymbol{A}\boldsymbol{h}(\boldsymbol{M})=\boldsymbol{1}_{p}^{\top}\mathrm{diag}\{\boldsymbol{R}_{0}^{-1}\boldsymbol{M}\}\boldsymbol{1}_{p}=\mathrm{trace}(\boldsymbol{R}_{0}^{-1}\boldsymbol{M})\).
Taking \(\boldsymbol{M}=\boldsymbol{\Psi}_{\beta}(\boldsymbol{x})\) and
multiplying by the scalar of (\ref{eq-uPsi}), which is
\(2a_{\beta}(\beta+1)^{\frac{p}{2}+2}\), gives (\ref{eq-IFexplicit}).
Finally,
\(\boldsymbol{u}_{\beta,1}(\boldsymbol{x},\boldsymbol{\eta}_{0})\) is a
linear function of
\(w_{\beta}(\boldsymbol{x})(\boldsymbol{x}-\boldsymbol{\mu}_{0})\) and
\(\boldsymbol{\Psi}_{\beta}(\boldsymbol{x})\), which are bounded if and
only if \(\beta>0\), in which case they tend to \(\boldsymbol{0}_{p}\)
and to \(\beta(\beta+1)^{-(\frac{p}{2}+1)}\boldsymbol{R}_{0}\) as
\(\left\Vert \boldsymbol{x}\right\Vert \rightarrow\infty\).

\subsection{\texorpdfstring{Proof of Theorem
\ref{thm-IF2}}{Proof of Theorem }}\label{sec-ProofIF2}

Let
\(\boldsymbol{\upsilon}(\epsilon)=\mathrm{E}_{F_{\epsilon,\boldsymbol{x}}}[\boldsymbol{u}_{\beta,2}(\boldsymbol{X},\widetilde{\boldsymbol{\eta}}_{\beta}(F_{\epsilon,\boldsymbol{x}}))]\),
the components of
\(\overline{\boldsymbol{U}}_{\beta}(F_{\epsilon,\boldsymbol{x}})\) which
correspond to the correlation parameters, the remaining ones being zero
by the definition of \(\widetilde{\boldsymbol{T}}_{\beta}\). Then
\(\boldsymbol{\upsilon}(0)=\boldsymbol{0}\) and, by
Theorem~\ref{thm-IFest} and (\ref{eq-uPsi}),
\begin{equation}\protect\phantomsection\label{eq-dupsilon}{
\left.\frac{\partial}{\partial\epsilon}\boldsymbol{\upsilon}(\epsilon)\right\vert _{\epsilon=0}
=\boldsymbol{u}_{\beta,2}(\boldsymbol{x},\boldsymbol{\eta}_{0})-\boldsymbol{J}_{\beta,21}(\boldsymbol{\eta}_{0})\boldsymbol{J}_{\beta,11}^{-1}(\boldsymbol{\eta}_{0})\boldsymbol{u}_{\beta,1}(\boldsymbol{x},\boldsymbol{\eta}_{0})
=\frac{\boldsymbol{T}^{\top}\boldsymbol{A}_{\perp}\boldsymbol{h}(\boldsymbol{\Psi}_{\beta}(\boldsymbol{x}))}{2(2\pi)^{\frac{\beta p}{2}}\left\vert \boldsymbol{\Sigma}_{0}\right\vert ^{\frac{\beta}{2}}}.
}\end{equation} Since
\(\widetilde{R}_{\beta}(F_{\epsilon,\boldsymbol{x}})\) is a quadratic
form in \(\boldsymbol{\upsilon}(\epsilon)\), its first derivative at
\(\epsilon=0\) is zero and one half of its second derivative is the
quadratic form of (\ref{eq-dupsilon}) in the block of
\(\boldsymbol{K}_{\beta}^{-1}(\boldsymbol{\eta}_{0})\) which corresponds
to the correlation parameters, which is
\(\{b_{\beta}\kappa_{1}(p,\beta)\boldsymbol{T}^{\top}\boldsymbol{A}_{\perp}\boldsymbol{T}\}^{-1}\)
by Lemma~\ref{lem-perp} (a). By Lemma~\ref{lem-perp} (b) and (c) and
\((2\pi)^{-\beta p}\left\vert \boldsymbol{\Sigma}_{0}\right\vert ^{-\beta}=4b_{\beta}\),
\[
\mathrm{IF}^{(2)}(\boldsymbol{x},\widetilde{R}_{\beta},F_{\boldsymbol{\theta}_{0}})=\frac{1}{\kappa_{1}(p,\beta)}\mathrm{trace}\left(\left(\boldsymbol{R}_{0}^{-1}\boldsymbol{\Psi}_{\beta}^{\ast}(\boldsymbol{x})\right)^{2}\right),
\] which is (\ref{eq-IF2}). For \(\beta>0\) it is bounded because
\(\boldsymbol{\Psi}_{\beta}(\boldsymbol{x})\) is bounded and
\(\boldsymbol{M}\mapsto\boldsymbol{M}^{\ast}\) is linear. For
\(\beta=0\),
\(\boldsymbol{\Psi}_{0}(\boldsymbol{x})=\boldsymbol{z}(\boldsymbol{x})\boldsymbol{z}^{\top}(\boldsymbol{x})-\boldsymbol{R}_{0}\)
and \(\boldsymbol{R}_{0}^{\ast}=\boldsymbol{0}_{p\times p}\), because
\(\boldsymbol{R}_{0}=\boldsymbol{D}\boldsymbol{R}_{0}+\boldsymbol{R}_{0}\boldsymbol{D}\)
with \(\boldsymbol{D}=\frac{1}{2}\boldsymbol{I}_{p}\), so that
\(\mathrm{IF}^{(2)}\) is one half of the homogeneous polynomial of
degree four
\(\mathrm{trace}((\boldsymbol{R}_{0}^{-1}(\boldsymbol{z}\boldsymbol{z}^{\top})^{\ast})^{2})\)
in \(\boldsymbol{z}=\boldsymbol{z}(\boldsymbol{x})\), which is not
identically zero because \(r\geq1\), and hence it is unbounded. For
\(\boldsymbol{R}_{0}=\boldsymbol{I}_{p}\),
\(\boldsymbol{M}^{\ast}=\boldsymbol{M}-\mathrm{diag}\{\boldsymbol{M}\}\)
and
\(\mathrm{trace}((\boldsymbol{M}^{\ast})^{2})=2\sum_{i<j}M_{ij}^{2}\),
and the off-diagonal elements of
\(\boldsymbol{\Psi}_{\beta}(\boldsymbol{x})\) are
\(w_{\beta}(\boldsymbol{x})z_{i}(\boldsymbol{x})z_{j}(\boldsymbol{x})\),
which gives the first expression of (\ref{eq-IF2ind}). For the supremum,
for a fixed value \(t\) of \(\boldsymbol{z}^{\top}\boldsymbol{z}\) the
maximum of \(\sum_{i<j}z_{i}^{2}z_{j}^{2}\) is \(\frac{p-1}{2p}t^{2}\),
reached when all the \(z_{j}^{2}\) are equal, and the maximum of
\(t^{2}\exp\{-\beta t\}\) is \(4e^{-2}\beta^{-2}\), reached at
\(t=2/\beta\).

\subsection{\texorpdfstring{Proof of Theorem
\ref{thm-PIF}}{Proof of Theorem }}\label{sec-ProofPIF}

Under \(F_{n,\epsilon,\boldsymbol{x}}^{P}\) the observations form a
triangular array of row-wise independent and identically distributed
random vectors, and \[
\begin{aligned}
\sqrt{n}\mathrm{E}_{F_{n,\epsilon,\boldsymbol{x}}^{P}}[\boldsymbol{u}_{\beta}(\boldsymbol{X},\mathrm{vech}(\boldsymbol{\Sigma}_{0}))]
&=\left(1-\frac{\epsilon}{\sqrt{n}}\right)\sqrt{n}\mathrm{E}_{\boldsymbol{\theta}_{n}}[\boldsymbol{u}_{\beta}(\boldsymbol{X},\mathrm{vech}(\boldsymbol{\Sigma}_{0}))]+\epsilon\boldsymbol{u}_{\beta}(\boldsymbol{x},\mathrm{vech}(\boldsymbol{\Sigma}_{0}))\\
&\underset{n\rightarrow\infty}{\longrightarrow}\boldsymbol{J}_{\beta}(\mathrm{vech}(\boldsymbol{\Sigma}_{0}))\boldsymbol{h}+\epsilon\boldsymbol{u}_{\beta}(\boldsymbol{x},\mathrm{vech}(\boldsymbol{\Sigma}_{0})),
\end{aligned}
\] while the variance-covariance matrix of
\(\boldsymbol{u}_{\beta}(\boldsymbol{X},\mathrm{vech}(\boldsymbol{\Sigma}_{0}))\)
converges to
\(\boldsymbol{K}_{\beta}(\mathrm{vech}(\boldsymbol{\Sigma}_{0}))\). By
the Lindeberg-Feller central limit theorem, (\ref{eq-LeCam}) holds with
the mean replaced by this limit. Since the proportion of contamination
tends to zero, the uniform convergences used in the proof of
Theorem~\ref{thm-root} also hold under
\(F_{n,\epsilon,\boldsymbol{x}}^{P}\), so the consistency of the
restricted estimators, (\ref{eq-linear}) and (\ref{eq-projscore}) remain
valid, as in Ghosh et al.~(2016) and Basu et al.~(2022b) for
contaminated contiguous sequences. By (\ref{eq-uPsi}) and
\((2\pi)^{-\frac{\beta p}{2}}\left\vert \boldsymbol{\Sigma}_{0}\right\vert ^{-\frac{\beta}{2}}=4a_{\beta}(\beta+1)^{\frac{p}{2}+2}\),
the asymptotic mean of
\(\sqrt{n}\boldsymbol{U}_{\beta,n,2}(\widetilde{\boldsymbol{\eta}}_{\beta})\)
is
\(2a_{\beta}\boldsymbol{T}^{\top}\boldsymbol{A}_{\perp}\boldsymbol{h}(\boldsymbol{\Delta}_{\beta}(\epsilon,\boldsymbol{x}))\),
and the proof of Theorem~\ref{thm-contig} with \(\boldsymbol{\Delta}\)
replaced by \(\boldsymbol{\Delta}_{\beta}(\epsilon,\boldsymbol{x})\)
gives (\ref{eq-deltastar}), and (\ref{eq-deltaexp}) follows by expanding
the quadratic form, since \(\boldsymbol{M}\mapsto\boldsymbol{M}^{\ast}\)
is linear, and from (\ref{eq-IF2}). Since
\(\frac{\partial}{\partial\delta}P(\chi_{r}^{2}(\delta)>c)=f_{\chi_{r+2}^{2}(\delta)}(c)\)
and
\(\left.\frac{\partial}{\partial\epsilon}\delta_{\beta}^{\ast}(\epsilon,\boldsymbol{x})\right\vert _{\epsilon=0}=\mathrm{ARE}(\beta,p)(1+\beta)^{\frac{p}{2}+2}\mathrm{trace}(\boldsymbol{R}_{0}^{-1}\boldsymbol{\Delta}^{\ast}\boldsymbol{R}_{0}^{-1}\boldsymbol{\Psi}_{\beta}^{\ast}(\boldsymbol{x}))\),
we obtain (\ref{eq-PIF}). With
\(\boldsymbol{\Delta}=\boldsymbol{0}_{p\times p}\),
\(\delta_{\beta}^{\ast}(\epsilon,\boldsymbol{x})=\epsilon^{2}\mathrm{IF}^{(2)}(\boldsymbol{x},\widetilde{R}_{\beta},F_{\boldsymbol{\theta}_{0}})\),
which gives (\ref{eq-levelcont}) and a level influence function equal to
zero. The last statements follow from the boundedness of
\(\boldsymbol{\Psi}_{\beta}(\boldsymbol{x})\) for \(\beta>0\) and, for
\(\beta=0\), from
\(\mathrm{trace}(\boldsymbol{R}_{0}^{-1}\boldsymbol{\Delta}^{\ast}\boldsymbol{R}_{0}^{-1}\boldsymbol{\Psi}_{0}^{\ast}(\boldsymbol{x}))=\mathrm{trace}(\boldsymbol{R}_{0}^{-1}\boldsymbol{\Delta}^{\ast}\boldsymbol{R}_{0}^{-1}\boldsymbol{z}(\boldsymbol{x})\boldsymbol{z}^{\top}(\boldsymbol{x}))\),
a quadratic form in \(\boldsymbol{z}(\boldsymbol{x})\) which is not
identically zero when
\(\boldsymbol{\Delta}^{\ast}\neq\boldsymbol{0}_{p\times p}\), and from
the unboundedness of \(\mathrm{IF}^{(2)}\).

\end{document}